\documentclass[12pt,oneside,letterpaper]{amsbook}
\newcommand\myfootnote[1]
	{ 
		\footnote{ #1 }
		
   	}

\usepackage{amsmath}
\usepackage{amsthm}
\usepackage{amscd}
\usepackage{amssymb}
\usepackage{ifthen}
\usepackage{multicol}

\usepackage[dvips]{graphics}
\usepackage{floatflt}
\usepackage{pstricks}

\newtheorem{Theorem}{Theorem}
\newtheorem{theorem}{Theorem}[section]

\newtheorem{proposition}[theorem]{Proposition}
\newtheorem{lemma}[theorem]{Lemma}
\newtheorem{corollary}[theorem]{Corollary}

\theoremstyle{remark}

\newtheorem*{note}{Note}
\newtheorem*{notation}{Notation}

\newtheorem*{definition}{Definition}

\theoremstyle{theorem}
\newtheorem*{theorem*}{Theorem}

\newcommand\R{{\mathbb R}}

\newcommand\Z{{\mathbb Z}}

\newcommand\N{{\mathbb N}}

\newcommand\Q{{\mathbb Q}}

\newcommand\hyp{{\mathbb H}}

\newcommand\union\cup
\newcommand\infinity\infty
\newcommand\from\colon

\newcommand\bdy{\partial}

\newcommand\length{ \text{len} } 
\newcommand\area{\text {area} } 
\newcommand\depth{ \text{depth} }

\renewcommand\ker{ \mathrm { ker } } 
\newcommand\im{ \mathrm { im } } 
\newcommand\rank{ \mathrm { rank } } 

\usepackage{hyperref}

\begin{document}
 \newboolean{verbose}		%
\setboolean{verbose}{true}	
 \begin{titlepage}

 \begin{center}
\renewcommand{\baselinestretch}{2}
	{\huge Algorithmic Properties of Relatively Hyperbolic Groups \\ }
\quad \\ 
	by Donovan Yves Rebbechi \\
\quad \\
	A dissertation submitted to the Graduate School-Newark \\
\quad \\
	Rutgers, The State University of New Jersey \\
\quad \\
	in partial fulfillment of the requirements  \\
\quad \\
	for the degree of \\
\quad \\
	Doctor of Philosophy \\
\quad \\
	Graduate Program in Mathematical Sciences \\
\quad \\
	Written under the direction of  \\
\quad \\
	Professor Lee Mosher \\
\quad \\
	and approved by \\
\quad \\
	\noindent
	\begin{tabular}{lr}
	Professor Lee Mosher & \underline{\qquad \qquad \qquad \qquad} \\
\quad \\	Professor Mark Feighn & \underline{\qquad \qquad \qquad \qquad}
\\
\quad \\	Professor Walter Neumann & \underline{\qquad \qquad \qquad
\qquad}\\
\quad \\
	Professor Ulrich Oertel & \underline{\qquad \qquad \qquad
\qquad} \\
	\end{tabular}
	\end{center}
\quad \\
	\begin{center}
	Date: April 2001
	\end{center}

\end{titlepage}
\pagenumbering{roman}
\setcounter{page}{2}

\renewcommand{\baselinestretch}{1}

{\center
	Algorithmic Properties of Relatively Hyperbolic Groups\\
	By Donovan Rebbechi\\
	Thesis Director: Professor Lee Mosher\\
}

\section*{Abstract}

The following discourse is inspired by the works on hyperbolic groups of 
Epstein, and
Neumann/Reeves. In \cite{epstein1}, it is shown that geometrically finite
hyperbolic groups are biautomatic. In \cite{neumann1}, it is shown that
virtually central extensions of word hyperbolic groups are
biautomatic. We prove the following generalisation:

\begin{Theorem}\label{maintheorem1}
Let $\mathcal{H}$ be a geometrically finite hyperbolic group. Let
$\sigma\in H^2 (\mathcal{H})$ and suppose that $\sigma\vert P= 0$ for
any 
parabolic subgroup $P$ of $\mathcal{H}$. Then the extension of
$\mathcal{H}$ by $\sigma$ is biautomatic
\end{Theorem}

We also prove another generalisation of the result in \cite{epstein1}.

\begin{Theorem}\label{maintheorem2}
Let $G$ be hyperbolic relative to $H$, with the bounded coset penetration
property. Let $H$ be a biautomatic group with a prefix-closed normal form.
Then $G$ is biautomatic.
\end{Theorem}

Based on these two results, it seems reasonable to conjecture the following
(which the author believes can be proven with a simple generalisation of 
 the argument in
 \ref{maintheorem1}):
Let $G$ be hyperbolic relative to $H$, where $H$ has a prefixed closed
biautomatic structure. $\sigma\in H^2(G)$ and suppose that 
$\sigma\vert_H = 0$. Then the extension of $G$ by $\sigma$ is biautomatic.

\section*{Acknowledgements}

        First and foremost, I'd like to express thanks to 
        thesis advisor, Dr. Lee Mosher, for his valuable guidance and
        inspiration, and for being there for his graduate students 
        during his sabbatical.

        I'd also like to thank Dr. Walter Neumann for his insights and 
        corrections were very helpful,
        especially his insights on the falsification by fellow traveler 
        property. Not only was Dr. Neumann helpful as a committee member,
        he inspired me to work on geometric group theory during my time at
        The University of Melbourne, where he gave a fascinating class on
        geometric group theory.

        Finally, I'd like to thank the other members of the committee, and 
        the Rutgers faculty for their friendly and helpful attitude.

 \tableofcontents
 \listoffigures
\pagebreak
\pagenumbering{arabic}

\section {Introduction}

Here is a general overview of the paper:
First, we review background material.
\begin{itemize}
\item We start by presenting several motivating examples for theorem \ref{maintheorem1}. (chapter 2)
\item In chapter 3, we do a review of essential combinatorial group theory. We
        also introduce a generalized fellow traveler property, used in the 
        proof of theorem \ref{maintheorem1}.
\item In chapter 4, we review CAT(0) and CAT($-k$) spaces, and develop a
        CAT($-k$) theory for ideal triangles.
\end{itemize}

We then prove theorem \ref{maintheorem1}. The proof uses the following essential 
steps:

\begin{enumerate}
\item First, (chapter 5) we prove a theorem of Farb (see \cite{farb1}). Our proof requires
        us to use the CAT($-k$) theory developed earlier.
        The theorem proves that geometrically finite hyperbolic groups are 
        relatively hyperbolic (in the strong sense -- they also have the bounded
        coset penetration property).
\item We then show (chapter 6) that relatively hyperbolic groups have a linear 
        \emph{electric isoperimetric inequality}.
\item We apply (in chapter 7) an argument of Neumann-Reeves (\cite{neumann3}). This article
        proves that central extensions of word-hyperbolic groups are 
        biautomatic. The main properties of word-hyperbolic groups 
        used in the proof are the linear isoperimetric
        inequality and the falsification by fellow traveler property.
        The electric isoperimetric inequality and our version of the
        falsification by fellow traveler property turn out to be sufficient 
        replacements.
\end{enumerate}

The interesting part of this theorem is step 2. The other results, steps 1 and 3
are based on the work of others, and it is step 2 that makes these old arguments
work in a new context.

We then prove theorem \ref{maintheorem2}. The proof goes as follows:
\begin{enumerate}
\item   First we have to deal with the following problem: we are using an
argument similar to that offered by Epstein in \cite{epstein1}. This argument takes
advantage of the fact that the group in question embeds into a Gromov hyperbolic metric
space. (namely, $\hyp^n$) So the first part of the problem involves finding a Gromov-hyperbolic
space to embed the Cayley graph in. Hence the first step is to construct this space, the cusped- off
Cayley complex (see 8.1-2). This space is a 2-complex with a weighted area and length function.

\item
The next step is to show that this space is indeed Gromov-hyperbolic. We do this
by first showing that it has the property that there exists a linear function $f$
such that the area of a loop with length $n$ in this space is bound above $f(n)$. (see 8.3)

\item
We then show that in this context, the linear area function implies Gromov hyperbolicity. (see 8.4)

\item 
Given this, it is a fairly simple matter to apply an argument similar to
(if somewhat simpler than)  that used by Epstein, to show biautomaticity. (see chapter 9)
\end{enumerate}

 \part{Geometrically Finite Hyperbolic Groups}

%

\section {Examples}

In order to convince the sceptical reader that this result is worth
proving, we present numerous examples of groups that satisfy the
hypothesis of the conjecture. We attain examples by two methods: 

\begin{enumerate}
\item Dehn filling on hyperbolic link complements
\item knot complements in manifolds $M$ that have second cohomology of
sufficiently high rank (rank 3 is sufficient) 
\end{enumerate}

The latter is in some sense an easier way of obtaining a rich
collection of examples, since it readily provides us with spaces that
have high rank cohomology and low rank cohomology on the boundary. 

The former is interesting as all 3-manifolds arise by way of Dehn
filling on $S^3$, so it illustrates how examples can be constructed
``from scratch''

\subsection {Group Cohomology on Manifolds}

It is a well known result (for example, see \cite{brown1} for a proof) 
	that a $K(\pi , 1 )$ 
space ( in particular, any complete hyperbolic manifold ) $M$ satisfies 
$$ \forall i\in \N , \ H^i(M)\cong H^i (\pi_1 ( M ) ) $$
We use this result to construct examples. 

We need to define a notion of ``adding a boundary'' to a complete,
non-compact, finite volume hyperbolic manifold $M$. We do this as follows: a cusp 
is homeomorphic to the Cartesian product  of a half open interval
$[0,1)$ and a quotient $E$ of $R^n$. The cusp $[0,1) \times E $ embeds 
into $[0,1] \times E$ via the inclusion map. We call 
$[0,1] \times E$ the \emph{compactified cusp}, we call the manifold 
obtained from $M$ by replacing all its cusps with the corresponding 
compactified cusps the \emph{augmentation} of $M$. 

We also need to define the notion of ``link complement''.
A closed solid torus is a space homeomorphic to 
the Cartesian product of a closed 2-ball and $S^1$. 
Similarly, an open solid torus is the product of an open
2-ball and $S^1$.
For the purposes of our discussion, a link complement will 
mean the complement in $S^3$ of disjoint embedded closed solid
tori.   Sometimes, we will need to work with compact spaces. 
In this case, we use the topological closure of the link complement.
By topological closure, we simply mean the closure of the link complement
as a subset of the space $S^3$. This space is the space that we would have
obtained by using the complement of open tori instead of closed tori.

Our examples need to satisfy the hypothesis of the following lemma: 

\begin{lemma}\label{manifolds}

Let $M$ be  the augmentation of 
a non-compact, finite volume, complete
hyperbolic manifold with finitely many cusps, and let $\sigma\in
H^2(M)$ with $\sigma\vert \bdy M  = 0$. Then $\pi_1 ( M ) $ satisfies
the hypothesis of \ref{maintheorem1}.

\end{lemma}

\begin{proof}
This is a well known result. For a reference, see \cite{brown1}.
\end{proof}

\subsection {Examples That Arise by Way of Dehn Surgery on $S^3$}
So we search for manifolds $M$ which satisfy the hypothesis of
\ref{manifolds} -- that is, manifolds $M$ and non-zero cohomology 
classes $\sigma \in H^2 (M)$ with $\sigma\vert_{\bdy M } = 0$. 
Our first candidates are Dehn filled link
complements. In particular, we consider the topological closure $U$ of 
the complement of a link
of $n$ components, and the manifold $M$ obtained by Dehn-filling $n-m$
of the boundary tori $T_i$ of $U$. 
We study the homology of $M$ via the exact Mayer Vietoris sequence
$$H^1 ( \bdy M ) \xrightarrow{\theta^*} H^2 ( M, \bdy M ) \xrightarrow{\phi^*}
H^2(M) \xrightarrow{\psi^*} H^2 ( \bdy M ) $$

First, we briefly discuss Dehn fillings.

\begin{definition}
        Let $V$ be an open solid torus embedded in $S^3$. 
        A \emph{meridian} in $V$ is a simple closed curve in 
        $\bdy { \overline{V} }$ that bounds a 
        disk in $V$.
        A \emph{longitude} of $V$ is a simple closed curve in $\bdy \overline{V}$ that 
        intersects a meridian at exactly one point, and is null-homologous 
        in $S^3 - V$.
\end{definition}

\begin{definition}
        Let $H$ be an embedded
        solid torus in $\R^{3}$ or $S^{3}$
        Let $p,q$ be coprime integers.
        Let $\alpha$ be a meridian of $H$ and $\beta$ be a longitude.
        We can obtain a manifold $M$ by gluing a solid torus 
        $V$ to $M - H $ by a map that takes the meridian of $V$ to
        the curve $\alpha^p \beta^q$. The resulting manifold $M$ 
        is said to be obtained from $\R^3$ ( or $S^3$ ) by $(p,q)$ Dehn surgery
        on $H$.

\end{definition}
        Let $U\subset S^3$ be the closure of a link complement of $n$ components.
        Let $T_i$ be the boundary tori of $U$.
        Let $M$ be the manifold obtained by filling $n-m$ of the $T_i$ 
        with a $(p_i, q_i)$ filling along meridian $\alpha_i$ and longitude 
        $\beta_i$ with a solid torus $V_i$. So $M$ has $m$ boundary components.
        Let $u_i , v_i $ be coprime integers such that the
        meridian of $V_i$ is identified with $u_i \alpha_i + v_i \beta_i$. 

We proceed as follows: $\ker( \psi^* )\cong \im( \phi^* )$, so our goal is to
show that $\im \phi^* $ has rank $> 0$,  because 
$\im (\phi^* ) \subset H^2(M)$ and cocycles in $\im(\phi^* )$ evaluate to
$0$ on $\bdy M$.
Now observe that :
$$\rank ( \im \phi^* ) = 
\rank ( H^2 ( M , \bdy M ) ) - \rank ( \ker ( \phi^* ) ) = 
\rank ( H^2 ( M , \bdy M ) ) - \rank ( \im ( \theta^* ) )  $$
Using Poincar\'e duality and the above equation, we obtain
$$\rank( \im (\phi^*) ) 
 = \rank ( H_1(M) ) - \rank ( \im (\theta_*)  ) $$
where $\theta_*: H_1 (\bdy M ) \rightarrow H_1 (M) $
is the homomorphism induced by inclusion.

So we proceed to explicitly calculate $H_1 ( M ) $ and study the map
$\theta_* : H_1(\bdy M ) \rightarrow H_1 ( M ) $. Observe that since $\bdy
M $ is the disjoint union of $m$ distinct tori, $H_1(\bdy M ) =
\Z^{2m} 
$.  Observe that $\alpha_i$
generate $H_1(U)\cong \Z^n$. 
We need a preliminary definition.
\begin{definition}
Given a set $\{ v_i\}$ of $n$ vectors in a Noetherian 
ring $\mathcal{R}$, 
the \emph{rank} of the set
is the number of vectors in a maximal linearly independent subset of 
$\{ v_i \}$ and the \emph{nullity} of $\{ v_i\}$ is equal to
the difference of $n$ and the rank of $\{ v_i \}$.
\end{definition}
\begin{lemma}\label{filling} 
Suppose that the set $\{ u_i \alpha_i + v_i \beta_i \} $ are linearly
dependent elements of $H_1(U)$. 
Then 
 the kernel of the map $\psi^* : H^2(M)\rightarrow H^2(\bdy M )$ has
rank greater than or
equal to the nullity of the set $\{ u_i \alpha_i + v_i \beta_i \}$. 

\end{lemma}
\begin{proof}
The proof follows easily from the following sub-lemma:

\begin{lemma}[Sub-lemma] 
Suppose that the set $\{ u_i \alpha_i + v_i \beta_i \} $ are linearly
dependent elements of $H_1(U)$ . Then $\rank(H_1(M)) - m \geq \mathrm{nullity}
( \{ u_i \alpha_i + v_i \beta_i \} )$
\end{lemma}
We defer the proof of the sub-lemma.

We now prove that $\rank ( \im ( \theta_* ) ) = m$. 
Applying this fact and the sub-lemma
will complete the proof.

We study the map $\theta_* : H_1(\bdy M) \rightarrow H_1 (M)$. 
We study this via the map $H_1( \bdy U) \rightarrow H_1(U)$. 
Notice that the image of this map has rank $n$ , and moreover, 
each component of $\bdy U $ contains a cycle 
corresponding to a generator $\alpha_i$ of $H_1(U)$.
Therefore each component of $\bdy U$ contains {\em exactly} one generator
of $H_1(U)$.
So for each component of $\bdy U $, there is a cycle 
that represents a cohomology class of infinite order in $\bdy U$ and is a 
boundary in $U$, hence a boundary in $M$ ( since $U\subset M$ ).
 So every component of $\bdy M$ contains exactly 
one generator that is 
a boundary in $M$. It follows that the
 rank of $\ker(\theta_*)$ is equal to $m$, the 
number of boundary components of $\bdy M$. Hence
$$\rank(\im(\theta_*) ) = \rank ( H_1 ( \bdy M ) ) -m = 2m -m = m$$

\end{proof}

We now prove the sub-lemma.

\begin{proof}[Proof of sub-lemma]
Since 
the set $\{ u_i \alpha_i + v_i \beta_i \} $ is linearly dependent, we
have that the rank of the group $\mathcal{K}$ generated by $\{ u_i \alpha_i + v_i \beta_i \} $
is less than $n-m$. Using the Mayer Vietoris sequence 
$H_1 ( U\cap V ) 
\xrightarrow{\rho_*} H_1(U)\oplus H_1(V) 
\xrightarrow{\xi_*} H_1 ( U \cup V ) $ , we obtain an exact sequence 
$$\Z^{2(n-m)} \xrightarrow{\rho_*} \Z^n \oplus \Z^{n-m} \xrightarrow{\xi_*} H_1(U\cup V)$$

For each solid torus $V_i$, choose a meridian $\mu_i$ and a longitude
$\lambda_i$. Identify $\bdy V_i$ with a component of $\bdy U $ via a
map which sends $\lambda_i$ to $p_i\alpha_i + q_i \beta_i \in H_1(U)$ and
sends $\mu_i$ to $u_i\alpha_i + v_i\beta_i \in H_1(U)$ where $u_i,v_i$ are
coprime solutions to the equation $p_i v_i + q_i u_i = 1$ . 
We view elements of $H_1(U)\oplus H_1(V) \cong \Z^n \oplus \Z^{n-m} $
as ordered pairs $(x,y)$ where $x\in H_1(U) $ and $y\in H_1(V)$.
For each component $V_i$ of $V$, we denote the corresponding generator
of $H_1(V)$ by $e_i$.
We write each $\beta_i = \sum_j k_{ij} \alpha_j $ where $k_{ij}\in \Z$.

Observe that $k_{i,i} = 0$ for all $i$. This is because $\alpha_i$ is a meridian
of $V_i$ and $\beta_i$ is a longitude. So $\beta_i$ is null-homologous 
in $M-V_i$, and so are $\alpha_j$ for all $j\neq i$. So taking images under
the map induced by the inclusion $ M \hookrightarrow S^3 - V_i$, we get 
$ 0 = k_{i,i} \alpha_i $.

Writing the vectors $\rho_* ( \lambda_i )$ in terms of the basis
$\{ \alpha_1, \dots, \alpha_n, e_1, \dots, e_{n-m} \}$,
we get $\rho_* ( \lambda_i ) = e_i + x $ where $x$ is a sum of the
$\alpha_i$. So the vectors $\rho_*(\lambda_i)$ are linearly independent,
so they span a space of dimension $n-m$.
The vectors $\rho_*(\mu_i )$ are given by $\sum  \delta_{i,j} u_{i,j} + v_{i,j} k_{i,j}$
where $\delta_{i,j} = 0 $ for $i\neq j$ and $1$ for $i=j$. The space 
spanned by the vectors $\rho_*( \mu_i )$ trivially intersects that of the 
space spanned by $\rho_*( \lambda_i )$. So

$$
\rank(\im ( \rho_* ) ) = \\
\rank(\im ( \lambda_i ))  + \rank ( \im \rho_* (\mu_i ) )   = \\
( n - m )  + \rank ( \im \rho_* (\mu_i ) )  
$$

So 
$$\rank(\ker(\xi_*)) =  \rank ( \im ( \rho_* ) )  = n - m + \rank ( \im \rho_* (\mu_i ) ) $$ and 

\begin{gather*}
\rank( H_1(M) ) - m \geq \\
\rank ( \im ( \xi_*)) -m = \\
  \rank ( H_1 ( U ) \oplus H_1 (V) ) - \rank (\ker(\xi_*))  -m =  \\ 
  (2n-m) - ( n - m + \rank ( \im ( \rho_* ( \mu_i ) ) ) )  -m = \\
 n - \rank ( \im ( \rho_* ( \mu_i ) ) )     -m=\\
 n - ( n - m - \text{nullity} ( \{ u_i \alpha_i + v_i \beta_i \} ) ) -m =\\
 m  + \text{nullity} ( \{ u_i \alpha_i + v_i \beta_i \} ) -m =\\
  \text{nullity} ( \{ u_i \alpha_i + v_i \beta_i \} )  \\
\end{gather*}

\end{proof}[Proof of sub-lemma]

This result raises the following question: how easy is it to Dehn-fill a hyperbolic link complement in such a way that
\begin{enumerate}
        \item the resulting manifold is hyperbolic and
        \item the hypothesis of \ref{filling} is satisfied
\end{enumerate}
To answer this question, we need to discuss some background material in \cite{thurston1} and use one of the key results. 

\begin{definition} The operation of Dehn surgery on a manifold (possibly with boundary)
$M$ 
 is parametrized by an ordered pair of coprime integers  
$(u_i , v_i )$. There is an embedding between these surgeries and the complex
plane, and this can be extended to a an embedding to the complex sphere (the
``trivial surgery'' where no solid torus is glued in corresponds to the point at
infinity).  For convenience, we sometimes use a rational number
instead of a pair, to denote a pair of coprime surgery coefficients, and a
vector of rational numbers to denote a particular surgery.
\end{definition}

\begin{theorem}[Thurston]
Given a set of disjoint embedded solid tori in a complete 
hyperbolic manifold $M$,
if $M=M_{\infty,\dots , \infty }$ admits a hyperbolic structure, then 
the result of some Dehn surgery on $M$ on those tori
admits a hyperbolic structure for all but finitely many possible fillings.
\end{theorem}

Note that this means that if either $p_i$ or $q_i$ is very large, 
the surgery will always result in a hyperbolic structure.,

We briefly discuss linking numbers, linking matrices, and linking graphs.

\begin{definition} 
A {\em link of n components} is a disjoint union of solid tori embedded into $S^3$. The 
{\em link complement } associated with the link $\cup_i V_i$  is $S^3-\cup_{i}V_i$. The solid tori are the {\em link components}
\end{definition}

\begin{definition}
Given an oriented link of $L=\cup_i V_i$ of $n$ components embedded into an oriented $S^3$, the {\em linking number} $k_{i,j}$ of the $i_{\mathrm{th}} $ component
with respect to the $j_{\mathrm{th}}$ component is defined as follows: 
first, choose a generating set for $H_1(S^3 - L ) $. Each generator is a 
meridian of some $T_i=\bdy V_i$. Denote this meridian by $\alpha_i$.  (ie given  $i_* : S^3-L \rightarrow S^3 - T$.) 
Let $\iota_{i*} : H_1(S^3-L )\rightarrow H_1(S^3 - V_i)$ be the homomorphism induced by inclusion. 
Using the right hand orientation on $T_i$ induced by the orientation on $S^3$, 
choose a longitude $\beta_i$ of $T_i$ with $\iota_{i*} ( \beta_i ) = 0$ . 
Then we define $k_{i,j}$ to be the number such that $\iota_{i*} ( \beta_j ) = k_{i,j} \alpha_i $ in $H_1 ( S^3 - L ) $.
Note that $k_{i,j} = - k_{j,i}$.
\end{definition}

\begin{definition}
The skew symmetric matrix $\{ k_{i,j} \}$ is called the {\em linking matrix}
\end{definition}
\begin{definition}
The {\em linking graph} associated with an oriented link complement
is a labelled directed graph with a
vertex $x_i$ for each link component $T_i$ , and for each positive $k_{i,j}$, 
an edge labelled $k_{i,j}$ from $x_i$ to $x_j$. 
\end{definition}

To satisfy the hypothesis of \ref{filling}, we need 
the elements $u_i \alpha_i + v_i \beta_i$ to be linearly independent 
in $H_1(U)$.
Recall that $u_i \alpha_i + v_i \beta_i = u_i \alpha_i + v_i \sum_{j} k_{i,j} \alpha_j$.
We can divide each of these vectors by the scalar $v_i$ without having any effect on the space
they span, so we have have vectors $( u_i / v_i  ) \alpha_i + \sum_{j} k_{i,j} \alpha_j$.
This is equivalent to the proposition that the matrix
$$
\begin{pmatrix}
u_1 /v_1 &  k_{1,2} &  k_{1,3} &  k_{1,4} & \dots &  k_{1,n}  \\
 k_{2,1} & u_2/v_2 &  k_{2,3} & k_{2,4} & \dots &  k_{2,n}  \\
\hdotsfor{6} \\
\end{pmatrix}
$$
has rank less than $n-m$. 

We approach this problem with the assumption that the $K_{i,j}$ are determined, and
we look for solutions for the vectors $\bar{\tau}=(\tau_1, \dots, \tau _{n-m} )$  
and the rational numbers $u_i/ v_i$
satisfying $\overline{\tau} C  = \overline{0}$
where $u_i^2 + v_i^2 $ is close to $\infty$ for each $i$.
In the case where $C$ is a square matrix with at least two  
rows and no rows or columns of zeros, there is a fairly simple solution. 
However, this in itself is not much help because it corresponds with the 
scenario where all link components are filled, so the resulting manifold 
is compact.
However, there are similar solutions in more difficult cases.
If the linking graph $\Gamma$ contains at least two components,
one of which contains an edge, we can obtain such a matrix by choosing our
fillings carefully , and possibly renumbering:

\begin{lemma}\label{wishful}
Let $\Gamma$ be the linking graph of a link that contains two 
components, one of which contains at
least one edge. 
 Then for all hyperbolic link complements $U$ whose linking graph is $\Gamma$, there is 
 a hyperbolic Dehn filling of the closure of $U$ such that the filled manifold $M$ has 
 the property that the map $H^2(M) \rightarrow H^2(\bdy M )$ has a kernel of nonzero rank. 
\end{lemma}
\begin{proof}
The components that we will fill are those corresponding to the vertices in the
component of the linking graph that contains at least one edge. After
renumbering, we assume that the filled components are $1, \dots ,n-m$. We need
this matrix to have linearly dependent rows. Note that the columns 
$n-m+1, \dots, n$ in the matrix $K'_{i,j} $  consist entirely of zeroes. 

The choice of $(\tau_1,\dots ,\tau_{n-m} )$ determines the $u_i / v_i$ as
follows: 
$$u_i / v_i  = \sum_{1\leq j\leq n-m} k_{i,j} \tau_j / \tau_i $$ 
It is important to note that each column vector of the 
matrix $K'_{i,j}$ is nonzero 
to ensure that the above does not give us $u_i / v_i =0 $. This is true 
because we chose to fill the link components corresponding to a single 
component of
the linking graph.
However, we have the subtle problem of ensuring that $u_i^2 + v_i^2$ can be made close to $\infty$ for each $i$. We do this by requiring the following:
\begin{enumerate}
        \item for each $1\leq i \leq n-m$, there is a $j$ for which $k_{i,j}\neq 0 $ 
        \item $\tau_{i+1} / \tau_i > t $ where $1\leq i < n-m $ and  $t$ is some large constant.

\end{enumerate}

In the following lemma, it is shown that by adjusting the value of $t$, we can ensure that each $1 / v_i < u_i / v_i <\epsilon $ 
or $1 / u_i < v_i / u_i < \epsilon $ , so either $u_i > 1 / \epsilon $ or $v_i > 1 / \epsilon $. So by choosing $\epsilon $ suitably small, we attain a hyperbolic manifold.

\end{proof}

\begin{proposition}
        Let $L$ be a link complement and let $K_{i,j}$ be an $n\times n$
        diagonal block
        of the linking matrix of $L$ (that is, $K_{i,j}$ is a sub-matrix
                        corresponding to a component of the linking graph).
        Let $A$ be the $(n-1)\times n$ matrix consisting of the first
        $(n-1)$ rows of $K_{i,j}$. Let $r$ be the column rank of 
        $M$. Then it is possible to choose $r-1$ of the $u_i / v_i$ 
        independently 
        in such a way that the manifold resulting from the fillings
        $(u_1,,v_1) , \dots , (u_{n-1}, v_{n-1} ) $ satisfies 
        the hypothesis of \ref{filling}.
        
\end{proposition}
\begin{proof}
        The goal is to choose fillings to make the matrix
$$
B = 
\begin{pmatrix}
u_1 /v_1 &  k_{1,2} &  k_{1,3} &  k_{1,4} & \dots &  k_{1,n}  \\
 k_{2,1} & u_2/v_2 &  k_{2,3} & k_{2,4} & \dots &  k_{2,n}  \\
\hdotsfor{6} \\
\end{pmatrix}
$$
        linearly dependent.
        Denote by $\vec{x_j}$ the column vectors of $A$. 
        Denote by $S$ the column span of $A$.
        The column rank of $A$ is $r$, so choose a basis for the column
        space of $A$ which has $r$ column vectors. Moreover, we can choose 
        the basis to contain the vector $\vec{x_n}$ -- this vector must 
        be non zero because of the hypothesis that we are in a component
        of the linking graph. Denote these vectors 
        by $\vec{x_{j_1}}, \dots , \vec{x_{j_{r-1}}}$, and $\vec{x_n}$.

        Given some choice of $(u_i, v_i)$, the proposition that the rows 
        of $B$ are linearly dependent is equivalent to the proposition
        that there exists a non-zero vector $\alpha \in \mathbb{Q}^{n-1},
        \alpha = (\alpha_1 , \dots , \alpha_{n-1})$
        such that $\alpha$
        is perpendicular to each column of $B$. 
        This is in turn equivalent to the proposition that $\alpha$
        is perpendicular to $S$, which is equivalent to $\alpha$
        being perpendicular to each of the vectors 
        $\vec{x_{j_1}} , \dots , \vec{x_{j_{r-1}}} $ and $\vec{x_n}$.

        Now each $u_{j_i}$ is determined by the vector $\alpha$ as 
        follows: 
        $$u_{j_i} = - \alpha \cdot \vec{x_{j_i}}$$
        So the map $f: \mathbb{Q}^{n-1} \rightarrow \mathbb{Q}^{r-1} $
        given by 
        $$f(\alpha) = ( -\alpha \cdot \vec{x_{j_1}} , \dots  , -\alpha \cdot \vec{x_{j_{r-1}}} )$$
        is linear, and the kernel of $f$ is clearly the orthogonal complement of $S$.
        So the restriction of $f$ to $S$ is injective. Moreover 
        the images $f(\vec{x_{j_i}})$ form a basis for the image of $f$. 
        Taking linear combinations of these vectors, we  can
        choose $(r-1)$ of the coordinates of the resulting sum independently.

\end{proof}

This result gives rise to a lot of fillings, because it says that it's not that
difficult to attain the linear dependency condition in the matrix $B$, and that we have
a fair degree of freedom ($r-1$ degrees of freedom, to be precise) in choosing
the pairs $(u_i, v_i)$, and all but finitely many choices will result in a hyperbolic
filling.

\subsection{Specific Examples}
\begin{figure}[!htp]
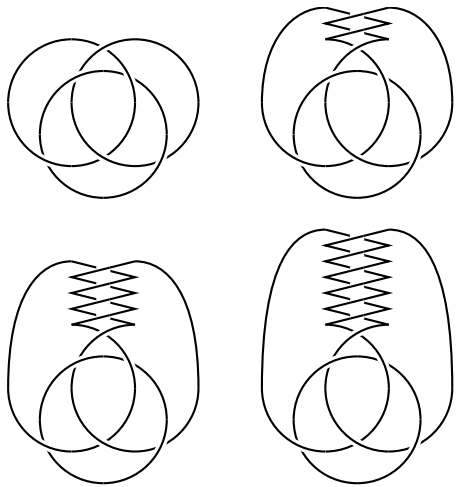
\label{rings}
\caption{The first four links $B_0 \dots B_3$ in the sequence $B_n$ of tangled Borromean rings.}
\begin{center}

\input borromean.tex
\end{center}

\end{figure}

We now give a family of examples, corresponding to \ref{wishful}. The family is:
the knotted Borromean rings , which are a series of links 
$B_0 , B_1 , B_2 , \dots , B_i , \dots $ where $B_i$ is obtained by
intertwining two of the link components $i$ times.  Fig \ref{rings} 
shows $B_0 , B_1$ and $B_2$. The author has verified that $B_0 \dots B_2$ are
hyperbolic, using snappea. However, to show that it is in fact hyperbolic for
all $B_i$, we need some further arguments.

There is a general method of looking for and constructing examples. 
This relies on a theorem of \cite{menasco}.

\begin{theorem}[Menasco]
        Let $L\subset \R^3\subset S^3$ be a link. Let $\pi:\R^3\rightarrow \R^2$
        be the projection map. Suppose further that $L$ is alternating with
        respect to $\pi$ and has no {\em trivial} crossings ( see diagram ).
        Then the following are true:
        \begin{enumerate}
                \item If $\pi(L)$ is connected, then $S^3-L$ is irreducible.
                \item Suppose that $S^3 - L $ is irreducible. Suppose that 
                        for each 
                        disk $D\subset \R^2$ where $\R^2$ is the projection plane
                        such that $\bdy D$ intersects $\pi(L)$ transversely at
                        exactly two points, neither of which are crossings, then
                        $\pi(L)\cap D$ is an embedded arc or $\pi(L)\cap \R^2-D$ 
                        is an embedded arc. ( by ``embedded arc'', we mean that it
                        contains no crossing points ). 
                        Then $L$ is prime.
        \end{enumerate}
\end{theorem}

        This gives us a strategy for finding links that satisfy the necessary 
        criteria for the hypothesis of our main theorem. Given a link $L$, we 
        perform the following checks:
        \begin{enumerate}
                \item Verify that the linking numbers satisfy the sufficient 
                        conditions.
                \item Check that $L$ satisfies the hypothesis for Menasco's 
                        theorem. If it does, then $L$ either has a complete 
                        hyperbolic structure of finite volume or $L$ is a torus link.
                        However, torus links have non-zero linking numbers, so we 
                        only need to check the hypothesis of Menasco's theorem.
        \end{enumerate}

\begin{figure}[!htp]
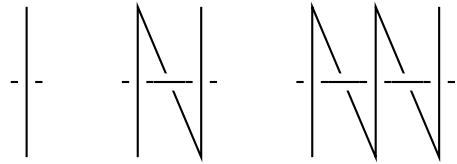
\label{tanglefig}
\begin{center}
\caption{The process of tangling}
\psset{xunit=.5pc}
\pspicture(-5,-2)(30,2)
\psline(-1,0)(1,0)
\psline[border=.1](0,-1)(0,1)

\rput{0}(7,0)
{
\psline(-1,0)(1,0)
\psline(3,0)(5,0)

\psline[border=.1](0,-1)(0,1)(4,-1)(4,1)

\psline[border=.1](1,0)(3,0)

}

\rput{0}(18,0)
{
\psline(-1,0)(1,0)
\psline(3,0)(5,0)
\psline(7,0)(9,0)

\psline[border=.1](0,-1)(0,1)(4,-1)(4,1)(8,-1)(8,1)

\psline[border=.1](1,0)(3,0)
\psline[border=.1](5,0)(7,0)

}
\endpspicture
\end{center}
\end{figure}

        Once we find a link that satisfies the appropriate conditions, we 
        can derive several links. We do this via a process called ``tangling''.
        This process is simple: given a crossing, where the components 
        at the crossing are either non-distinct, or have nonzero linking 
        number, we modify the link in a neighborhood of the crossing.
        Let the two link components be denoted by paths $\zeta$ and $\psi$. 
        Let $\zeta$ be the component that crosses over $\psi$. After reparametrisation,
        we can assume the crossing point corresponds with $\psi(0)$ and $\zeta(0)$. 
        Consider the restriction of $\zeta$ to $(-\epsilon, \epsilon)$ where $\epsilon$
        is small enough that $\zeta\vert_{( -\epsilon , \epsilon ) }$ doesn't encounter
        any other crossings. Similarly, choose $\epsilon'$ so that 
        $\psi\vert_{(-\epsilon' , \epsilon' ) }$
        doesn't encounter any other crossings.
        Take a small neighborhood $N$ of 
        $\pi(\psi\vert_{(-\epsilon' , \epsilon' )} \cup \zeta\vert_{ (-\epsilon, \epsilon ) } )$
        in $\R^2$. Taking the Cartesian product of $N$ with a sufficiently large open interval, 
        we obtain an open ball $B$ containing the crossing. Let $\zeta'$ be the  simple loop based at 
        $\zeta( \epsilon / 2  )$, such that 
        $\zeta'' = \zeta\vert_{ (-\epsilon, \epsilon /2 ) } \zeta'  \zeta\vert_{ (\epsilon/2 , \epsilon ) }$ 
        begins with an over-crossing ( after isotoping to eliminate trivial crossings and the
        self intersection  ) and 
        $\zeta'$ is a generator of $\pi_1 ( B - \psi ) $.
        We call the operation of replacing $\zeta$ by a curve isotopic to $\zeta'' $ {\em tangling}.
        The operation of tangling is illustrated in figure \ref{tanglefig}
        diagram.  The diagram also illustrates how we can repeat the 
        tangling operation near a crossing. This produces an infinite family of links.

\begin{lemma}
        Iterating the tangling operation preserves both hyperbolicity
        and the linking number hypothesis for all but finitely many iterations.
\end{lemma}
\begin{proof}
        First, we need to check that the following linking number hypotheses 
        are invariant under tangling for all but finitely many iterations:
        recall the hypotheses are that there exist $i,j_1, j_2$ such that
        \begin{enumerate}
                \item $K_{i, j_1 } \neq 0 , K_{i,j_2} \neq 0  K_{j_1 , j_2 } = 0 $
                \item The vertex $x_i$ of the linking graph is non-separating.
                \item If $c\in \Q, 0\leq j_3 \leq n$ such that $K_{j,j_3 } = cK_{j,i}$ 
                for all $j\leq n, j\neq i$, then $j_3 = i$
        \end{enumerate}
        (1) follows because iterating tangling will only produce a 0 linking number for
        a unique number of tangles. (2) is true for the same reason. (3) is similarly true.

        We also need to check that the hypothesis to Menasco's theorem is unaffected 
        by the operation of tangling. First, observe that the link projection $\pi(L)$ 
        cuts $\R^2$ into polygons. The hypotheses for Menasco's theorem is satisfied 
        if an only if $\pi(L)$ is connected and any two polygons intersect in at most one
        side. We observe that the operation of tangling introduces a new bigon that 
        intersects the adjacent polygons in a unique side. The addition of the new bigon 
        does not interfere with the adjacency relations of the polygons that existed 
        prior to the addition of the bigon. So adding the new bigon does not ``ruin'' 
        the hypothesis for Menasco's theorem. Since we insisted that the tangling be performed
        on two components with nonzero linking number, the resulting link is not a torus 
        link, hence it is hyperbolic.

\end{proof}

%
%
%
%
%
%

\section{A Review of Combinatorial Group Theory} 

We start with some preliminary definitions and notation.

\subsection{Group Presentations}

        Let $X$ be a finite set. 
        Let $X^{-1}$ be another (disjoint) set 
        in one to one correspondence with $X$ 
        and 
        for each $x_i \in X$ , 
        denote by $x_i^{-1}$ the corresponding element of $X^{-1}$.
        A { \em word on $X$ } is an element of 
        the free monoid on $\mathcal{M}(X)$ on $X\cup X^{-1}$. 
        We denote by $F(X)$ the free group on the set $X$, and given a word $w$ , denote
        by $\pi_{F(X)}w$ the image of $w$ under the projection to $F(X)$. Given two 
        words $u,v$, we say $u\sim v$ if $\pi_{F(X)}u = \pi_{F(X)}v$. It is easy to see that
        $\sim$ is an equivalence on $\mathcal{M}(X)$. Given a group
        $G$ and a surjective homomorphism $\pi_G: F(X)\rightarrow G$, we call $X$ {
        \em a generating set for $G$. } Elements of $X$ are called 
        \emph{generators}. Conversely, given a generating set, $X$ for a group
        $G$, we always denote by $\pi_G$ the projection homomorphism $\pi_G :
        F(X)\rightarrow G$. There is also a natural projection $\mathcal{M}(X)\rightarrow G$,
        this is called the {\em evaluation map} , and given $w\in \mathcal{M}(X)$, 
        its image under the
        evaluation map (also called its evaluation) is denoted by $\overline{w}$.
        A {\em group presentation} for a group $G$ is a pair $
        \langle X\vert R \rangle $ where $X$ is a generating set for $G$ and $R$ is a
        set of elements of $F(X)$ with the property that 
        $\ker \pi_G =_{F(X)} \{ prp^{-1} \vert p\in F(X), r\in R \}$.

\subsection{ The Cayley Graph and Metric Spaces }
        
        Given a finitely generated group $G$ with generating set $X$, the
        {\em Cayley graph of $G$ with respect to $X$, $\Gamma_{(G,X)} $  } is 
        a labelled directed graph with one vertex $v_g$ for each $g\in G $ and for each
        pair $g , ga $ where $a\in X,g\in G$, an edge from $v_{g} $ to $v_{ga}$. 
        We usually denote this by $\Gamma_G$  (the geometric properties of $\Gamma_G$
        that we are interested in do not depend on the choice of generating set, 
        and neither do most of our arguments.) 
        We define a length function on $(G,X)$ by defining 
        $|g|=\inf\{ \ell \vert g = a_1 , \dots , a_{\ell } , a_i\in X \forall i \}$,
        and given $g,h\in G $, we define a metric on $G$ by $d(g,h) = |g^{-1}h |$.
        Given a word $w=a_1 a_2 \dots a_i \in X^*$, denote by $w(t)$ the word
        $a_1 a_2 \dots a_t$. Now we can construct a corresponding path $\gamma$ 
        in $\Gamma_G$ such that $\gamma(t)$ maps to the vertex corresponding to $w(t)$ 
        for each $t\in \N$. Then the map can be extended to the edges of 
        $\Gamma_G$ by requiring that the path be an isometry on the open
        intervals $(n,n+1)$ where $n\in\N$.
        The homotopy classes of paths naturally correspond
        with elements of $F(X)$. 
        We can also define a corresponding metric on $\Gamma_G$ by taking the path
        metric induced by assigning a length of $1$ to each edge. Notice that the 
        $d_{G,X}(g,h) = d_{\Gamma(G,X)}(v_g , v_h ) $
        
        A geodesic metric space $M$ is said to be {\em $\delta-$hyperbolic } 
        if given any geodesic triangle $T$ in $M$ , any side of $T$ is contained in the 
        $\delta$-neighborhood of the union of the other two sides. We are interested
        in groups that act on geometric structures that are $\delta$-hyperbolic,
        or are closely related to $\delta$-hyperbolic spaces, because these groups
        enjoy nice geometric and algorithmic properties.
        
        A {\em pseudo-metric} space $(X,D)$ is a set $X$ with a 
        symmetric , non-negative distance function $D$ satisfying the triangle inequality.
        We will use pseudo metric spaces to 
        deal with the points in ``bad'' sets in a space (the idea is that
        we just define the distance between two points in a connected ``bad'' 
        set to be zero) 
        The pseudo metric does not uniquely define a ``good'' topology ,
        and using the technique used to get a topology from a metric space 
        does not produce a Hausdorff space.     
        A better topology is any Hausdorff topology where the open metric balls correspond
        to open sets and the sets of the form $\{ y : D(x,y)=0\}$ (where $x$ is an arbitrary 
        constant point of $X$) are closed.
        A geodesic in a pseudo-metric space is a path $\gamma$ such that 
        $\gamma\vert_{[t_1 , t_2 ] }  $ is no longer than $\alpha$ where $\alpha$
        is any path with $\alpha[0]=\gamma(t_1)$ and $\alpha(n) =\gamma(t_2) $

\subsection{ 2-Complexes, Groups and Groupoids }

        Any finitely presented group $G$ together with a presentation  $< X \vert R > $
        can be realized as the fundamental group
        $\pi_1(K)$ of some 2-complex $K$. This is done as follows: 
        first, we take a base vertex. For each generator
        of $G$, we attach a 1-cell (ie a loop) to this vertex.
        The resulting graph has fundamental group $F(X)$ where $X$ is the generating
        set for $G$. We attach 2-cells as follows: each $r\in R$, is a word in 
        $(X\cup X^{-1} )^*$,
        so it defines a loop $\gamma$ in the graph.     
        We glue a 2-cell to the graph 
        by identifying its boundary with $\gamma$. This gives rise to an alternative 
        definition of the Cayley graph of $G$: we can define it as the 1-skeleton
        of the universal covering space of this 2-complex.

        Now we move on to the more general notion of groupoids. 
        First, we need to define categories. 
        A {\em category} consists of the following structure:
        \begin{enumerate}
                \item A set of {\em objects }
                \item For any pair of objects $A,B$, a set of {\em morphisms }
                (maps $f: A\rightarrow B $) $hom(A,B)$ 
                \item An associative operation {\em composition} 
                $hom(B,C) \times hom(A,B) \rightarrow hom(A,C)$
                \item For each object $B$, there is a morphism $1_B$ such that 
                for any $f\in hom (A,B) $, $1_B \circ f  = f $ and if $h\in hom (B,C)$, then
                $h \circ  1_B  = h$.
        \end{enumerate}

        A group is a category with one object and invertible morphisms.
        The morphisms in a group are maps of the form 
        $\phi_g : G\rightarrow G , \phi_g(h)= gh$.
        A groupoid has invertible morphisms, but may have several objects.
        A nontrivial example of a groupoid is the groupoid of homotopy 
        classes of paths in a
        CW-complex. The objects are points of the space, and the morphisms
        are paths in the space. Notice that if we restrict ourselves to one 
        object in considering paths that begin and end at some point $p$, we
        have a group. The main difference between a groupoid and a group is 
        that in a groupoid setting , it is not always true that two morphisms are 
        composable.

        A groupoid is useful in situations where we have an action of a group 
        on a set with finitely many orbits, and we need an algebraic structure that
        preserves the structure of the group action. If there is only one orbit, the
        group is adequate. If there are several orbits, a groupoid is 
        sometimes more appropriate. In particular, each orbit 
        naturally corresponds with a unique groupoid object.

        We now discuss some properties of groupoids. 
        First, we need some definitions:
        \begin{definition}
        A {\em homomorphism of groupoids } is a map $\phi:G\rightarrow H$ 
        satisfying the
        following: 
                Let $g_1, g_2\in G$. 
                        If $g_1 g_2 $ exists, then $\phi(g_1) \phi(g_2 ) $ exists,
                        and is equal to $\phi(g_1 g_2 )$.
        \end{definition}
        \begin{definition}
        A {\em generating set}
        of a groupoid $H$ is a set of morphisms $X$ such that any morphism 
        in $H$ is some product of morphisms in $X$. A groupoid is said 
        to be {\em finitely generated} if it has some finite generating set.
        \end{definition}
        A generating set is insufficient to fully describe the generators 
        of a groupoid because it does not encode information regarding the
        composability of different generators.
        This motivates the following definition:
        \begin{definition}
        A \emph{generating graph}, $(X,\Gamma)$ is a labelled directed graph
        $\Gamma$ whose edges are in one to one correspondence with $X$. 
        There is a \emph{composability relation} 
        determined by $\Gamma$:  $g$ is composable with $h$ if the 
        path $gh$ can be traced out in $X$ (ie the initial vertex of $h$ 
                        coincides with the destination vertex of $g$)
        \end{definition}
        \begin{definition}
        A generating graph $(X,\Gamma)$ is a {\em free basis } for $H$ if for 
        any groupoid $H'$ and any function
        $f: X\rightarrow H'$ that has the property that $f(h)\circ f(g)$ 
        is well defined for all composable pairs $(g,h)\in X\times X$,  
        there is a unique homomorphism of groupoids 
        $\tilde{f}: H\rightarrow H'$.
        A groupoid $H$ is said to be a \emph{free groupoid} if it has a free
        basis.
        \end{definition}

        There is a general method to construct a free groupoid, similar to
        that used to construct a free group. First, we have a set of objects 
        $\{ A_i \}$
        for the groupoid. We also have a set of symbols. Each symbol is
        a morphism $ A_i \rightarrow A_j $ between two objects. Given a symbol
        $x : A_i \rightarrow A_j $, we have a symbol $x^{-1} : A_j \rightarrow A_i $.
        A word $w = x_1 x_2\dots x_n \in X^*$   is \emph{admissible} if
        the morphism $x_{i+1} \circ x_i$ is defined for each $i$. 
        If $x_1: A_i \rightarrow A_j$ and $x_n: A_j\rightarrow A_\ell$, 
        then we describe $w$ with the notation $w: A_i \rightarrow A_\ell$
        (this is a subtle abuse of notation, because $w$ is not a groupoid
        element, it is a word associated with a unique groupoid element).

        On the set of admissible words in $X^*$, we have 
        the operation of {\em free reduction}, 
        which corresponds to replacing a word of the form $u x x^{-1} v $ or 
        $u x^{-1} x v $ with $uv$ where $u,v\in X^* , x\in X$. A word in $X^*$
        is said to be {\em freely reduced} if it admits no free reduction. 
        It is a theorem (see \cite{mks}) that given a word in $X^*$, the outcome 
        of free cancellation is unique (this is usually invoked in group 
        theory, but
        observe that it is also applicable to the groupoid setting. 
        It is a fact about
        words, not groups or groupoids). 
        It is easy to see, but worth mentioning, that admissibility 
        and inadmissibility are both invariant under free cancellation.

        Now we are in a position to say what ``free groupoid'' means in this
        new sense. The free groupoid is the set of equivalence classes 
        of admissible words modulo the equivalence relation of free equality.
        The operation is word concatenation. The concatenation of two words
        is only defined if the result is admissible. (So two words 
        $w: A_i \rightarrow A_j, v: A_k \rightarrow A_\ell$ have a product
        if and only if $j = k$.)

        We now justify the name ``free groupoid''.
        Let $\phi$ be a map of a set of symbols $X$ into a groupoid such that
        for each $i,j$, there is a $\phi(A_i), \phi(A_j)$ such that if 
        $x : A_i \rightarrow A_j $, then 
        $\phi(x) : \phi ( A_i ) \rightarrow \phi ( A_j )$.
        Then we can define a map $\tilde{\phi}$ by the rule 
        $\phi ( x_1 x_2 \dots x_n ) = \phi ( x_1) \phi(x_2 ) \dots \phi(x_n)$. 
        This map is well defined by the uniqueness of free reduction. 

        Given a generating graph $(X,\Gamma)$, a relator is a morphism in the free 
        groupoid $F(X)$. The groupoid $< X\vert R > $ is defined by the set of 
        equivalence classes given by the rule $uv \sim urv $ where $ r\in R ; u, v \in F(X)$.

        We observe some properties of generating graphs. Corresponding
        with each object $A_i$, we have a vertex $v_i$. For each generating morphism 
        $x\in X: A_i\rightarrow A_j$,
        we have a labelled directed edge $e_x$ from $v_i $ to $v_j$. Note that the groupoid 
        of homotopy classes of edge paths beginning and ending at vertices is isomorphic 
        to the corresponding free groupoid.

        For a groupoid $< X\vert R> $ , we can construct a 2-complex.
        First, construct the graph corresponding to the free groupoid on $X$.
        Then for each relator $A_i\rightarrow A_j$, identify the 
        vertices $v_i , v_j$, then
        attach a 2-disk in the same manner as for a group.
        We will call this space the {\em associated 2-complex} of the groupoid.
        Note that the groupoid of homotopy classes starting and ending at the 
        vertices $v_i$ is isomorphic to $<X\vert R>$.
        We will call the 1-skeleton of this space the {\em associated graph}. 
        The {\em Cayley graph} of a groupoid is the 1-skeleton of the universal covering
        space of its associated 2-complex.

\subsection{ Regular Languages and Finite State Automata}
        A {\em language} , $L$ on an alphabet $X$ is a subset of $X^*$.
        A {\em deterministic finite state automaton} $M$ is a $5$-tuple 
        $$
                (
                \mathcal{S}, 
                s_0\in\mathcal{S}, 
                X, 
                f: X\times\mathcal{S} \rightarrow \mathcal{S}, 
                g:\mathcal{S}\rightarrow \{0,1\}
                )
        $$
        where $X$ is a collection of symbols, and $\mathcal{S}$ 
        is a finite set of objects called {\em states}. $s_0$ is called the
        \emph{initial state}, $f$ is referred to as the 
        {\em state function}. A state $s\in\mathcal{S}$ is referred to as an {\em accept state}
        if $g(s)=1$ and a {\em fail state} if $g(s)=0$. A state $s$ is said to be {\em terminal}
        if $f(x,s)=s$ for all $x\in X$. For each element $w=x_0 x_1 \dots x_n \in X^*$, 
        we assign a recursively defined state $s(w)$ as follows:
        $s(x_0 x_1 \dots x_n ) = f( x_n , s( x_0 x_1 \dots x_{n-1} )) $ and $s(\epsilon) = s_0$ where
        $\epsilon\in X^*$ is the empty word. A state $s\in\mathcal{S}$ is said to be {\em inaccessible} if 
        $w\in X^*$ implies that $s(w) \neq s$.
        The {\em language accepted by $M$ }, $L_M$ is the set of words
        $w\in X^*$ such that $s(w)$ is an accept state. A language, $L$ is {\em regular}
        if it is the language accepted by some finite state automaton $M$. 
        Because the most interesting aspect of an automaton (for our intents and purposes) 
        is the language accepted by it, we are not interested in inaccessible states.
        A {\em non deterministic finite state automaton} $M'$ is a $5$-tuple
        $$
        (
                \mathcal{S}, 
                s_0\in 2^\mathcal{S}, 
                X, 
                f: X\times\mathcal{S} \rightarrow 2^\mathcal{S}, 
                g:\mathcal{S}\rightarrow \{0,1\}
        )
        $$
        where $f$ and $g$ have similar properties. 
        There are a few key differences: we extend $f$ to
        a function 
        $\tilde{f}: X\times 2^\mathcal{S}\rightarrow 2^\mathcal{S}$ 
        by imposing the condition that 
        $\tilde{f} ( A\cup B ) = \tilde{f} ( A) \cup \tilde{f} (B) $. 
        Now we can define $s(w)\in 2^\mathcal{S}$
        in a way analogous to the definition used in the deterministic case. 
        We define the language accepted 
        $L_{M'}$ by $M'$ to be the set $\{w\in X^* \vert \exists s\in s(w) : g(s) = 1 \}$. 
        We can construct from $M'$ a deterministic automaton 
        $M''$ whose state set is $2^\mathcal{S}$ such
        that $L_{M''}=L_{M'}$. 
        Hence it is true that a language is regular if it is the language accepted
        by a finite state non deterministic automaton.
        
        A language $L\subset X^*$ on a group $G$ with a generating set $X$ is a surjection 
        $L\rightarrow G$. $L$ is a {\em rational structure}
        if $L$ is a regular language. 

\subsection{Automatic Groups}

\begin{definition}
A language $L\subset X^*$ is said to have the 
\emph{$k$ fellow traveler property} if given 
any two words $w_1, w_2$ in $L$, with $d(\bar w_1 , \bar w_2 ) \leq 1 $,
then $d(w_1(t), w_2(t)) < k$.
It is said to have the \emph{$k$ two sided fellow traveler property}
if for any $w_1, w_2 \in X^*$ such that there exists 
$x_1, x_2 \in X \cup 1$, with $\overline{x_1 w_1 x_2} = w_2$, 
$d( (x_1 w_1 x_2 )(t) , w_2(t) ) < k$.
\end{definition}

\begin{definition}
If we replace $w_2(t)$ with $w_2(\rho(t))$ where $\rho : \R \rightarrow \R$ is a 
(not necessarily strict) monotone increasing surjective function with $\rho(0) = 0$,
we get a weaker property which we will call the \emph{asynchronous fellow 
traveler property}.
\end{definition}

\begin{definition}
        An \emph{automatic structure}, $L$ on a group $G$ is a regular language
$L$ on $G$ with the fellow traveler property.
A regular language $L$ on $G$ is an \emph{asynchronous automatic structure}  
if it has the asynchronous fellow traveler property.
We use the term \emph{biautomatic} (or asynchronously biautomatic) to refer 
to structures which have the two sided fellow traveler property (or
two sided asynchronous fellow traveler property).
\end{definition}

\subsection{The Falsification by Fellow Traveller Property}
	The falsification by fellow traveller property is a means
	of showing that languages (in particular, languages consisting
	of geodesic words) are regular. For example, this property can
	be used to show that the language of geodesics in a Gromov 
	hyperbolic group is regular. Epstein makes a subtle generalisation
	of this property by pointing out that you can instead consider the
	language of geodesic words within a regular language on $X^*$.
	We further explore this notion, and aim to geometric idioms 
	for showing a language is regular, using geometric properties of
	the language.

	A {\em height function} $\mathcal{H}: X^*\rightarrow \Z$ 
	on the set of words $X^*$ on a generating set $X$ 
	of a group $G$	
	is a function that satisfies
	the property that
	for each $g\in G$, 
	$ \sup \{ \mathcal{H}(v) \vert \bar{v} =g \}$
	exists and is equal to $\mathcal{H}(w)$ for some $w\in X^*$ with $\bar{w} = g$
	If $w$ has this property, we call $w$ a {\em maximising word} of $\mathcal{H}$
	
	A language $L$ on a group $G$ is said to possess the 
	$\delta$-{\em falsification by fellow traveller property},
  	where $\delta >0$ is a constant, 	
	if there is a height function  $\mathcal{H}: X^* \rightarrow \Z$ on $X^*$ 
	with the following properties:
\begin{enumerate} 
	\item $L$ is the set of maximising words of $\mathcal{H}$
	\item for all  $u\in X^* - L$, there exists $v\in X^*$  with 
	$\mathcal{H}(u)<\mathcal{H}(v)$ , $\bar{u} = \bar{v}$
	and a monotone function $x: \R\rightarrow \R , x(0 ) = 0 $ such that 
	$$d\left[ u(t) , v(x(t)) \right] < \delta $$ for all $t\in \R$.
\end{enumerate}
	Any words $u, v $ such that 
	there is a monotone function $x: \R\rightarrow \R , x(0 ) = 0 $ such that 
	$$d\left[ u(t) , v(x(t)) \right] < \delta $$ for all $t\in \R$ are
	known as {\em fellow travellers} , and are said to {\em fellow travel} each other.

	At this stage, we are dealing with a very abstract setting. This is because we wish
	to construct a falsification by fellow traveller property that is applicable to 
	many different situations. For example, 
		$\mathcal{H}(a) = - \length(a)$.
		This is how one constructs an automatic structure on a word hyperbolic
		group.
%
%
%
	For the purpose of this discussion, we will be interested in 
	height functions that are well behaved. We develop a 
	several desirable properties a height function can have.

%
\begin{definition}
\begin{enumerate}
	\item Additivity: 
		For all $u,v\in X^*$,
		$\mathcal{H}(uv) = \mathcal{H}(u) + \mathcal{H}(v) $
		This property essentially implies that the height function 
		is a weighted length function, and is what is traditionally
		used in falsification by fellow traveller arguments.

	\item Strong translation invariance:
	We say that $\mathcal{H}$ is {\em strongly translation invariant} if the following holds:
	For all $u,u',v,v'\in X^* $ such that $\bar{u} = \bar{u'}, \bar{v} = \bar{v'} $,
	$$
		\mathcal{H}(uv) - \mathcal{H}(u'v') = 
		( \mathcal{H}(u) - \mathcal{H}(u') ) + ( \mathcal{H}(v) - \mathcal{H}(v') )
	$$
	This property is implied by additivity.

	\item Translation invariance: 
			For all $u,v\in X^*$, 
			$\mathcal{H}(u) = \mathcal{H}(v) \Rightarrow \mathcal{H}(uw) = \mathcal{H}(vw)$.
			This is a weaker property than strong translation invariance. 

	\item Left order-preserving: 
		For all $x\in X$, and for all $u,v\in X^*$ with $\bar{u} = \bar{v}$, if 
	$\mathcal{H}(u) \leq \mathcal{H}(v)$, then $\mathcal{H}(xu) \leq \mathcal{H}(xv)$.
	This is a weaker property than translation invariance. This implies a \emph{suffix-closure}
	property: for any maximising word $w = uv$, the word $v$ is also maximal. In particular,
	the trivial word is maximal.

	\item Right order-preserving:
		For all $x\in X$, and for all $u,v\in X^*$ with $\bar{u} = \bar{v}$, if 
	$\mathcal{H}(u) \leq \mathcal{H}(v)$, then $\mathcal{H}(ux) \leq \mathcal{H}(vx)$.
	This is a weaker property than translation invariance. This  \emph{prefix-closure}: 
	for any maximising word $w = uv$, the word $u$ is also maximal.
	As with left order-preserving height functions, this implies maximality 
	of the trivial word.
\end{enumerate}

\end{definition}

	Note that if a height function is both left and right order-preserving, then 
	any subword of a maximising word is maximising (subword-closure).

	Currently, the main result of this section requires the strong translation invariance
	property. It would be nice to weaken the hypothesis. Where possible, we will prove 
	intermediate results assuming the weaker properties.

%

	Given a group $G$, a number $\delta\in \R$ and an element $g\in G$, 
	we will use the notation $B_{\delta, g}$ to denote the set 
	$\{ h\in G\vert d_G ( g,h) \leq \delta \}$

	To each strongly translation invariant height function $\mathcal{H}$,
	we associate a {\em state function} 
	$\Phi_{\delta} : X^* \times B_{\delta, 1 }\rightarrow \Z $ 

	A {\em $B_\delta $-word  } is a word $w\in X^*$ 
	such that $w(t) \in B_{\delta, 1 } $ for all $t$.

	Given a group $G$, for each $g\in B_{\delta, 1 }\subset G $, 
	choose a $B_\delta$-word $z_g$ such that 
	$z_g$ is a maximising word.
	
	The state function is an inf taken over certain sets. We define
	these sets first.

	Let 
	$$
		V_{\delta, u,g } = 
		\left\{ 
			v\in X^* \vert 
			\bar{v} = \bar{u}g^{-1} , 	
			v \text{ } \delta-\text{fellow travels } u  
		\right\}
	$$
	Note that all fellow travellers $v$ of $u$  with $\bar{u} = \bar{v}$ are in some
	$V_{\delta,u,g}$. Now we can define $\Phi_{\delta}$.
	
	$$
		\Phi_{\delta}(u,g) = 
			\inf_{v\in V_{\delta,u,g}} 
			\left\{ 
				\mathcal{H}(u) - \mathcal{H}(vz_g)
			\right\}
	$$
	We now need to demonstrate the worthiness of $\Phi_{\delta}(u,g)$ as
	a state function. 

\begin{lemma}
	If $G$ has the $\delta$-falsification by fellow traveller property 
with respect to  a left order-preserving height function
	$\mathcal{H}$, $\Phi_{\delta}$ has the following property: 
	{\em 
		$u$ is maximising if and only if $\Phi_{\delta}( u, g ) \geq 0 $
		for all $g\in B_{\delta, 1 }$. 
	}
\end{lemma}
\begin{proof}
	First, one of the implications is clear: if $u$ is maximising, then
	$\Phi_{\delta}(u,g) \geq 0 $. This is immediate from the definition
	of $\Phi_{\delta}$.

	The converse however is nontrivial. We suppose that 
	$\Phi_{\delta}(u,g) \geq 0$.

	Let $v'$ be a maximising fellow traveller of $u$ with $\bar{v'} = \bar{u}$.
	Then $v'$ enters $B_{\delta , u } $ at some point. 
			This implies that there exists $t_0\in\Z$ such that
	$v'(t_0) = \bar{u}g^{-1}$ 
	and $v'(t) \in B_{\delta , \bar{u} } $ for all $t>t_0 $.
	We define words $v$ and $w$ with $vw = v'$ as follows:
	$v(t) = ( v(0), \dots , v(t_0 ) )$ and 
	$w(t) = ( v(t_0) , \dots , v(\length (t) ))$

	To complete the proof, we need to compare $u$ with $v'$.
	First we show that $vz_g$ is a maximising word.
	We do this by comparing $v' = vw$ with $vz_g$.
	Note that
	$  \mathcal{H}(z_g) \geq \mathcal{H}(w)$, so by 
	left order-preservation, $\mathcal{H}(vz_g) \geq \mathcal{H}(vw)$.

	The hypothesis that $\Phi_{\delta}(u,g) \geq 0$ implies that 
	$u$ has height greater than or equal to all fellow travellers
	$v'$ that end in $z_g$, and the above argument shows that there
	exists a maximising that ends in $z_g$. Therefore, $u$ is maximising.

\end{proof}
	$\Phi_{\delta}$ tells us how well $u$ `measures up' to it's competitors. 
	We introduce some notation: we define $\Phi_{\delta}(u)$ as the
	function
	$$\Phi_{\delta}(u) : B_{\delta,1 } \rightarrow \Z \qquad
		\Phi_{\delta}(u)(g) = \Phi_{\delta}(u,g)$$

	Given $u\in X^* , x\in X$, the critical question is this: does the state function 
	$\Phi_{\delta}(ux ) $ depend only on $\Phi_{\delta}( u ) $ and $x$ ?  We will see that 
	this is important in proving that the language of maximal words with regard
	to the height function $\mathcal{H}$ is regular.

\begin{lemma}\label{fft1}
	Let $G$ be a group generated by a set $X$, 
	and let $L\subset X^*$ be a language on $G$ that has the 
	$\delta$-falsification by fellow traveller property with 
	respect to a height function $\mathcal{H}$ on
	$X^*$.
	Suppose that $\mathcal{H}$ has the following properties:
\begin{enumerate}
	\item $\mathcal{H}$ is right order-preserving.
	Note that this  implies
	maximality of the trivial word.
	\item Bounded difference: 
	$\mathcal{H}$ has the {\em bounded difference} property if there exists $K\in \N$ such that 
	for all $w\in X^*, x\in X$, $$\vert \mathcal{H}(w) - \mathcal{H}(wx) \vert < K$$
	\item The function $\Phi_{\delta}(ux)$ is uniquely determined by the pair 
	$( \Phi_{\delta}(u) , x  )$ , ie is not dependent on $u$.

\end{enumerate}
	Then $L$ is a prefix-closed regular language.
\end{lemma}

	These conditions might sound somewhat contrived, but in practice, ``most''
	partial orders that one would want to define will satisfy this hypothesis.
	The bounded difference condition is almost certainly necessary.

\begin{proof}

	If $\Phi_{\delta}(ux)$ depends only on $\Phi_{\delta}(u)$ and $x$ , 
	( ie if $\Phi_{\delta}(u)$ and $x$ uniquely 
	determine $\Phi_{\delta}(ux)$, independently of $u$ ) then 
	we can define a finite state automaton $\mathcal{M}$ that accepts the language $L$.	
	The success states of $\mathcal{M}$ will be state functions $\Phi_{\delta}(u)$ with the range
	$\{ 0, \dots , 2k \delta \}$ , and the additional terminal fail state corresponds to
	state functions $\Phi_{\delta}(u')$ whose range lies outside this set.
	There are well defined transition functions 
	$\tau ( \Phi_{\delta}(u) , x ) = \Phi_{\delta}(ux)$. The initial state is given by $\Phi_{\delta}(1)$.
	Note that $\Phi_{\delta}(1)$ is {\em not } the fail state, because of our hypothesis
	that the trivial word is maximising.

	First,  the bounded difference
	property implies that the state functions $\Phi_{\delta}(w)(g) < 2 k\delta $. 
	If $w\in X^*$ , then for any 
	$g\in B_{\delta , 1  }$, there is a $z$ with 
	$\bar{z} = g , \length(z) \leq \delta$. $wzz^{-1}$ is a fellow traveller 
	of $w$ and by the bounded difference property, 
	$\mathcal{H} ( w ) - \mathcal{H}( wzz^{-1}  ) \leq 2\delta K$

	We now justify the terminal fail state. Let $w$ be any word that leads to
	the terminal fail state. 
		
	Now
	if $\Phi_{\delta}(w_{[0, t_0]})(g) < 0 $ for some $g\in B_{\delta , 1 }$ , then 
	$\mathcal{H} ( vz_g ) > \mathcal{H}(w_{0,t_0})$ for some $t_0$. So 
	if $w=w_{[0,t_0] } w' $, then we have that 
	$\mathcal{H} ( vz_gw' ) > \mathcal{H}(w_{0,t_0} w')$ 	
	by right order preservation.
	So $w\not\in L$.

	We also need to show that if $w\not\in L$ , then $w$
	labels a path from the initial state the terminal fail state.
	Let $w\in X^* - L$. Then by the
	falsification by fellow traveller property, there is some word $u\in X^*$ such 
	that $u$ $\delta$-fellow-travels $w$, and $\mathcal{H}(u) - \mathcal{H}(w) > 0 $
	and $\bar{u} = \bar{w}$. Decompose $u$ as $u=u'x$ where $x$ is the last generator
	in the word $w$. Then 
	$$\Phi_{\delta}(w)(x)\leq \mathcal{H}(w) - \mathcal{H}(u) < 0 $$ 
	So $\Phi_{\delta}(w)$ corresponds to the terminal fail state.

\end{proof}

	This lemma in itself is a little unwieldy. The next step will be to
	find some sufficient conditions for computability of the transition 
	functions.

\begin{lemma}\label{fft2}
	If $\mathcal{H}$ is strongly translation invariant height function on $X^*$, and 
	$\mathcal{H}$ is well defined on $F(X)$, then 
	$\Phi_{\delta}(ux)$ is uniquely determined by the pair $(\Phi_{\delta}(u) , x )$

\end{lemma}

	Note: this is the example that will become important later on. 
	It is applicable to height functions on central $\Z$ extensions of 
	groups $\Z\rightarrow_{\iota} E \rightarrow_{\pi} G$. 
	Given a section $\rho: G\rightarrow E $ corresponding to the extension,
	the height $\mathcal{H}$ can be given by 
	$\mathcal{H}(w) =  \iota^{-1} ( \bar{w} \rho\circ \pi ( \bar{w}^{-1} )  )$.

	Note that the hypothesis for this result is still fairly strong, it would
	certainly be nice to find a weaker hypothesis that worked.

\begin{proof}
	Let $\mathcal{W}_{\delta, x, g }$ be the set of words $w$ satisfying the following
	conditions: 
	\begin{enumerate}
		\item $w(t_0) = w(t_1 ) \Rightarrow t_0 = t_1 $
		\item $g^{-1}w(t) \in B_{\delta, 1 } \cup B_{\delta, x }$ for all $t$ , $g^{-1}w(0) \in
			B_{\delta, 1 } $ and $g^{-1}w(\length{w} ) \in B_{\delta , x }$
	\end{enumerate}			

	Let 
	$$\mathcal{W}_{\delta, x, g,h } = 
	\{ w\in \mathcal{W}_{\delta, x,g } , \bar{w} = gxh^{-1} \}$$

	We now wish to compute $\Phi_{\delta}(ux)(h)$ for any $h\in B_{\delta, 1 }$
	Now consider a fellow traveller $v'$ of a word $ux$ where $x\in X$
	where $\overline{v'} h = \overline{ux}$. 

	We introduce some notation for writing $\inf$ functions. 
	We will denote by $\inf_{ a\in A } f(a) $  
	the $\inf$ of
	the function $f(a)$ as $a$ ranges over the set $A$.
	We will denote by 
	$\inf_{ \{a \vert p(a)\} } f(a)$	the $\inf$ of $f(a)$ as 
	$a$ ranges over the set of values that satisfy some proposition $p(a)$.
	This notation is necessary to clarify the set over which the $\inf$ is being 
	evaluated.

%

$$		
		\Phi_{\delta}(ux)(h) = 
		 \inf_{ v'\in V_{\delta,ux,h }  } 
			\mathcal{H}(ux) - \mathcal{H}(v'z_h) 
$$

	Observe the following : any $v' \in V_{\delta, ux, h }$ can be decomposed
	in the following manner: let $g \in B_{ \delta , 1 }$ such that 
	$v'(t_0)g = \bar{u} $ for some $t_0$
	and $v'(t)\in B_{\delta , \bar{u} } \cup B_{\delta, \bar{ux}}$
	for all $t>t_0$ and  $v'_{ [0, t_0 ] } $ fellow travels $u$.
	The existence of such a 
	$t_0$ is guaranteed by the fact that $v'$ is a fellow traveller 
	of $ux$. We can write $v = v'_{ [ 0, t_0 ] }$ and 
	$w = v'_{ [ t_0, \length ( v' ) ] }$, hence $v'=vw$.
	So for any $v'\in V_{\delta, ux,h}$, we can find a triple 
	$(g, v, w )$ such that 
	$g\in B_{ \delta, 1 }$ , and $v\in V_{\delta, u,g}$ 
	satisfies $\bar{v} = \bar{u}g^{-1}$, and 
	$w\in\mathcal{W}_{\delta, x , g, h }$, which is constrained
	by $g$ , but not by $v$ , satisfies $\bar{w} = g\bar{x}h^{-1}$.

	Conversely, any triple $(g, v, w )$ such that $g\in B_{\delta, 1 }$,
	$v\in V_{\delta, u,g }$ and $w\in \mathcal{W}_{\delta, x, g,h  }$ uniquely determines 
	a fellow traveller of $ux$, namely $v'$.

	Consequently, given a function $f(v')$, it follows that 

	\begin{align*}	
		& \inf_{ 
					\{ v'\in V_{\delta, ux, h  } \} 
				}
				f(v' )  =\\
		& \inf_{ 
			\{      
				g,v,w \vert
				g\in B_{ \delta, 1 } ;
				v\in V_{\delta, u,g} ;
				w\in \mathcal{W}_{\delta, x, g,h  }
			\}
		}
		f(vw ) 
		=\\
		& \inf_{g \in B_{\delta , 1 } } 
			\inf_{v\in V_{\delta, u,g} } 
			\inf_{w\in \mathcal{W}_{\delta, x, g, h } }
			f(vw)
	\end{align*}	

	Note that the last two nested infs are independent, ie the choice of $v$ and 
	$w$ are independent of each other ( both depend only on $g$ ).
	So
	\begin{align*}
		& \Phi_{\delta}(ux)(h) = \\
		 & \inf_{ v'\in V_{\delta, ux,h }  } 
			\mathcal{H}(ux) - \mathcal{H}(v'z_h)  =\\
		& \inf_{g \in B_{\delta , 1 } } 
			\inf_{v\in V_{\delta, u,g} } 
			\inf_{w\in \mathcal{W}_{\delta, x, g, h  } }
			\mathcal{H}(ux) - \mathcal{H}(vwz_h)  =\\
		& \inf_{g \in B_{\delta , 1 } } 
			\inf_{v\in V_{\delta, u,g} } 
			\inf_{w\in \mathcal{W}_{\delta, x, g, h } }
			\mathcal{H}(ux) - \mathcal{H}(vz_gz_g^{-1}wz_h) + \mathcal{H}(z_gz_g^{-1}) =\\
		& \inf_{g \in B_{\delta , 1 } } 
			\inf_{v\in V_{\delta, u,g} } 
			\inf_{w\in \mathcal{W}_{\delta, x, g, h } }
			\left[ \mathcal{H}(u) - \mathcal{H}(vz_g)   \right] + 
			\left[ \mathcal{H}(x) - \mathcal{H}(z_g^{-1}wz_h)\right]  
		 + \mathcal{H}(z_gz_g^{-1}) 	=\\
		& \inf_{g \in B_{\delta , 1 } } 
		\left[
			\inf_{v\in V_{\delta, u,g} } 
			(
				 \mathcal{H}(u) - \mathcal{H}(vz_g)    
			)
			+ 
			\inf_{w\in \mathcal{W}_{\delta, x, g, h } }
			( 
				\mathcal{H}(x) - \mathcal{H}(z_g^{-1}wz_h) + \mathcal{H}(z_gz_g^{-1}) 
			)  
		\right]
		=\\
		& \inf_{g \in B_{\delta , 1 } } 
		\left[
			\Phi_{\delta}(u)(g)	
			+ 
			\inf_{w\in \mathcal{W}_{\delta, x, g, h } }
			( 
				\mathcal{H}(x) - \mathcal{H}(z_g^{-1}wz_h) + \mathcal{H}(z_gz_g^{-1})
			)  
		\right]
	\end{align*}
	So $\Phi_{\delta}(ux)(h)$ depends only on $\Phi_{\delta}(u)$ and the function 
			$\inf_{w\in \mathcal{W}_{\delta, x, g, h } }
			( 
				\mathcal{H}(x) - \mathcal{H}(z_g^{-1}wz_h)
			)$
	which has finite domain and range and only depends on $g$ , $h$ and $x$.

\end{proof}

\section{CAT 0 Spaces}

\subsection{Background Material}

        We review a geometric notion of non-positive and negative 
        curvature based on ``comparison triangles''. The idea is that 
        triangles in a CAT($k$) space should be at least as ``thin'' as
        their counterparts in the Riemannian manifold of constant 
        sectional curvature $k$. The 
        case we are interested in is where $k=0$ (hence we compare with 
        Euclidean triangles) or $k<0$ (ie we compare with triangles in
        $\hyp^n$, or equivalently, $\hyp^2$.)

        The basics of CAT(0) geometry, including the results presented
        are widely understood by geometers, and recently, a comprehensive 
        text \cite{bridson1} has been published by Bridson and Haefliger,
        and it is recommended reading for anyone wishing to pursue CAT(0)
        geometry.  This section uses several results of that book.
\begin{definition}
        A {\em geodesic triangle} $T$ is a triple of points $(a,b,c)$ and a geodesic
        arc between each pair of points, parametrized by arc-length. The geodesic
        arcs are called {\em sides} of $T$.
\end{definition}

\begin{figure}
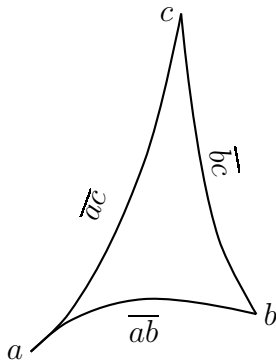
\label{trianglenotation}

\caption{Notation for Triangles}
\begin{center}
        \input triangle_notation.tex
\end{center}

\end{figure}

\begin{definition}
        A metric space $M$ is said to be {\em CAT(k)} where $k\in \R , k\leq 0 $
        \myfootnote
        {
                There is a more general definition of a CAT($k$) space which allows
                the possibility that 
                $k>0$. In this more general definition, 
                we only require comparison triangles to exist if the distance
                between any two vertices of the triangle is no more than
                $\pi/k$.
        }
        if the following is true: for any geodesic triangle $T\subset M$, with 
        vertices $a, b, c$, and sides $\alpha, \beta, \gamma$, there is a 
        {\em comparison triangle $\overline{T}$ }
        in the space of constant curvature $k$ such that 
        \begin{enumerate}
                \item There is a map $\phi: T \rightarrow \overline{T}$ such 
                that the restriction of $\phi$  to each side of $T$
                is an isometry with respect to the subspace metrics on $T$ 
                and $\overline{T}$
                \item $\phi$ takes vertices to vertices
                \item For any points $\alpha(s),\beta(t) $, 
                $$d( \phi\circ\alpha(s) ,\phi\circ\beta(t)) \leq d(\alpha(s), \beta(t))$$
                Similarly for $( \alpha(s), \gamma(t) )$ and $( \beta(s),\gamma(t) )$
        \end{enumerate}
        
\end{definition}

\begin{definition}
        Angles in CAT(0) spaces are defined using the law of cosines:
        the {\em approximating angle} of a triangle $(a,b,c)$ is given by
        $$A^2 = B^2 + C^2 - 2BC \cos\alpha$$
        where $A,B$ and $C$ are the lengths of the sides
        opposite $a,b$ and $c$ respectively, and $\alpha$ is the angle opposite the side $A$.
        Given $\epsilon > 0 $, we can choose a triangle $T_\epsilon ( a )$
        whose vertices are
        $a, b'_\epsilon$ and $c'_\epsilon$ where 
        $b'_\epsilon$ and $c'_\epsilon$ are (resp)  
        points on $\overline{ab}$
        and $\overline{ac}$ within distance $\epsilon$ of $a$. The angle at $a$ is defined
        as the limit as $\epsilon$ approaches $0$ of the approximating angle at $a$ 
        of triangle $T_\epsilon$. This is well defined if $M$ is sufficiently nice 
        ( a Riemannian manifold is more than sufficient ). This angle is known as
        the {\em Alexandrov} angle.
\end{definition}

For convenience, given a set/point/path/arc $A$ in a $CAT(k)$ space, $M$, we 
refer to its counterpart in the comparison space as $\overline{A}$.

\begin{definition}
        The {\em boundary} of a CAT(0) space $M$ is the set of all (infinite) geodesic  
        rays modulo the relation that $r \sim r'$  if there exists $K\in \R$ 
        such that $d(r(t), r'(t) )< K \ \forall t\in\R$. It is easy to check that
        $\sim$ is indeed an equivalence.  We call these equivalence
        classes ``points at infinity''
        or ``boundary points''.

\end{definition}

We now state without proof some lemmas of \cite{bridson1}
\begin{lemma}
        Let $M$ be a complete CAT(0) space, and let 
        $r: \left[ 0, \infty \right) \rightarrow M$ be 
        a geodesic ray with $x=r(0)$. Then for any $y\in M$, there is a unique 
        geodesic ray $r'$ such that $r'(0) = y$ and $r' \sim r$. 
\end{lemma}     

We also obtain a useful 
intermediate result from the proof of the above:
\begin{lemma}[Sub-lemma]\label{cat0_sublemma}
For large $t$, the $r_t$ are ``close'' in the following sense:
        Given $a > 0 ,\epsilon > 0 , s > 0 \in \R $, there exists a 
        $K(\epsilon, a , s ) \in \R$ 
        such that if $t>K , t' > 0 $, then 
        $d( r_t ( s ), r_{t + t' } ( s ) ) < \epsilon$.
\end{lemma}

\begin{definition}
        Let $M$ be a CAT(0) metric space. Let 
        $r: \left[ 0 , \infty \right) \rightarrow M$ 
        be a geodesic ray. The {\em Buseman function} associated with 
        $r$ is given by 
        $b_r ( x ) = \lim_{t\rightarrow\infty} d( x, r(t) ) - t $.
\end{definition}
        It remains to show that this definition makes sense.
        The following result of \cite{bridson1} demonstrates that it does.
\begin{lemma}
        The function $f(t) = d( x , r(t) ) - t $ is non-increasing and
        bounded below by $-d( x, r(0 ) )$, hence the limit used 
        to define Buseman function exists.
\end{lemma}

We need another result of \cite{bridson1}, which says that Buseman functions 
only change by a constant when the base-point is changed:

\begin{lemma}
        Let $M$ be a CAT($-k$) space. Let $x,y\in M$. Let $\beta, \gamma$ be 
                        geodesic rays from ( resp. ) $x$ and $y$ such that 
                        $\beta\sim\gamma$.      Then 
                The Buseman functions of $\gamma$ and $\beta$ differ by some constant.
\end{lemma}

\begin{definition}
        Let $M$ be a CAT(0) metric space. A {\em horosphere} centered at a 
        boundary point given by a geodesic ray 
        $r: \left[ 0, \infty \right) \rightarrow M$
        is a level set of the Buseman function $b_r(x): M\rightarrow \R$.
        A  {\em horoball} is a set of the form $\{ x : b_r(x) <    K \} $ where 
        $K\in\R$.
\end{definition}

\ifthenelse{\boolean{verbose}}
{
\begin{note}
        In $\hyp^n$, horospheres are very simple constructions: consider the upper
        half space model with coordinates $ (x_1 , \dots, x_n ) $ where $x_n >0$.
        In the case where the boundary point is $\infty$, a 
        horosphere is given by the set 
        $\{ x\in \hyp^n : x_n = k \}$ where $k$ is some positive real number, 
        and the corresponding horoball is 
        given by the set of points $\{ (x_1, \dots ,x_n) | x_n > k \}$.
        In the cases of the other boundary points $p = ( x_1 , \dots , x_{n-1} , 0 )$, 
        let $S$ be the Euclidean sphere 
        \myfootnote{
                By ``the Euclidean sphere'', we mean this set:
                $\{ ( y_1 , \dots , y_k )  | (k-y_k)^2 + \sum_1^{k-1} ( y_i - x_i )^2 = k^2 \}$
                }
        of radius $ k $ about the point 
        $ (x_1, \dots , x_{n-1} , k )$, so $S$ is tangential to $p$. Then $S - \{p\}$
        is a horosphere about $p$, and the interior of $S$ is a horoball.
        
\end{note}

\subsection{Convexity Results}

We need some basic results on convexity. We will prove that horoballs 
        and balls in CAT(0) spaces are convex.

\begin{lemma}[Hypotenuse of a right angled triangle ]
        Let $T = ( a,b,c )$ with sides $\alpha, \beta, \gamma$ 
        (see figure \ref{trianglenotation})
        be a right-angled triangle in a CAT(0) space
        $M$, where the angle at the vertex $a$ is a right angle. Then the longest
        side of $T$ is $\overline{bc}$.
\end{lemma}
\begin{proof}
        Given a triangle $T=(a,b,c)$, with a right angle at 
        the vertex $a$, we can take a 
        $T_\epsilon (a) = (a,b_\epsilon, c_\epsilon)$ (with $b_\epsilon \in 
                        \overline{ab}$
        and $c_\epsilon \in \overline{ac}$) for very small
        $\epsilon$. 
        Investigating the image $(a', b'_\epsilon, c'_\epsilon)$
        of $T_\epsilon$
        in the comparison triangle $(a',b',c')$, 
        we observe that $d( b'_\epsilon , c'_\epsilon ) \geq d ( b_\epsilon , c_\epsilon ) 
        \geq B^2 + C^2 - \epsilon'$ 
        where 
        $$
                B = d(a, b_\epsilon ) = d(a', b'_\epsilon ), 
                C = d( a, c_\epsilon )= d( a', c'_\epsilon )$$
                and $\epsilon'(\epsilon)$
        is a small constant.
        conclude that the 
        comparison triangle also has an angle of at least $\pi / 2 $  
        at the vertex $a'$. 
        So the comparison triangle's
        longest side is $\overline{b'c'}$, hence the longest 
        side of $T$ is $\overline{bc}$.
\end{proof}

                The following lemma is a result of \cite{bridson1}
\begin{lemma}[ $\tilde{M}$ balls are convex ]\label{horoballs_convex}
        Let $\tilde{M}$ be a contractible $CAT(0)$ space. Then any sphere in
        $\tilde{M}$ is convex, ie if $\gamma\subset\tilde{M}$ is a geodesic 
        such that
        the endpoints of $\gamma$ are in some ball $B\subset \tilde{M}$ , 
        then $\gamma\subset B$.
\end{lemma}

        The result fairly easily generalizes to horoballs:
\begin{lemma}[Horoballs Are Convex]
        Let $B$ be a horoball in a CAT(0) space $M$. Then $B$ is a convex set.
\end{lemma}
\begin{proof}
        Let $x, y \subset B$, where $b_\gamma(x) \geq b_\gamma (y)$ for any geodesic ray
        $\gamma$ corresponding with $B$. 
        Let $r$ be a geodesic ray such that $r(0) = x $ 
        and $r$ corresponds to the horosphere $H = \bdy B$. Consider the 
        triangle $T$ with vertices $x,y,r(t)$. 
        Let $\alpha$ be a geodesic arc from $x$ to $y$.
        Let $\overline{T}$ be a comparison triangle
        and let $\phi$ be the comparison map. Let $z$ be some point on the side
        $\alpha$. Then either 
        $$
                d(z, r(t) ) \leq d(\overline{z} , \overline{ r(t) } ) \leq d( \overline{x}, \overline{ r(t) } )
                = d ( x, r(t) ) 
        $$ 
                or the above holds after replacing $x$ with $y$. 
                Since $b_r (x ) \geq b_r(y) $, we have that 
                $d( r(t) , x ) - t \geq d( r(t) , y ) -t $ for sufficiently 
                large $t$, hence $d(r(t) , x ) \geq d(r(t) , y ) $.
                So if $t$ is large enough, we have that 
                $d( z, r(t) ) \leq d(x,r(t) ) = t $. It follows that
                $$b_r(z) = \lim_{ t \rightarrow \infty } \left(d(z,r(t) ) - t < 0\right) = b_r(x) $$ 
                so 
                the restriction of 
                $b_r$ to  $ \alpha$ realizes its maximum on an endpoint, ie the 
                geodesic arc $\alpha$ is contained in $B$.
        
\end{proof}

}

\subsection{Triangles With an Ideal Vertex}
        We now investigate triangles with an ideal vertex. 
        An {\em generalized triangle} in
        a CAT($-k$) space  is a set of three 
        geodesic sides, which may be arcs, rays or lines; such that 
        any two sides $\alpha, \beta$
        either meet at an endpoint, or 
        $\alpha\vert_{ \left[0 , \infty \right) } \sim
        \beta\vert_{ \left[0 , \infty \right) } $ ( possibly after 
        re-orienting  $\alpha$ and $\beta$ )
        Choosing comparison triangles for ideal triangles is nontrivial,
        because ideal triangles in $\hyp^2$ depend on more data than the side 
        lengths alone. 
        We will study ideal triangles with one ideal vertex, ie a geodesic 
        arc $\alpha$ between vertices $x$ and $y$ with rays $r_x, r_y$
        such that $r_x ( 0 ) =x, r_y(0)=y, r_x\sim r_y$.

        We define a comparison triangle in the following manner:

\begin{definition}
        Let $T$ be a triangle in a CAT($-k$) space with vertices $x$ and $y$. 
        Then $\overline{T}$ is said to be a comparison triangle for $T$ if :
        \begin{enumerate}
                \item   There is a map $\phi: T\rightarrow \overline{T}$  such
that the restriction of $\phi$ to each side of $T$ is an isometry.
                \item   Let $x$ be a vertex of $T$. 
                                Let $\overline{\beta}$ be the geodesic ray of $\overline{T}$ 
                                originating at $\overline{x}$. Let $\beta$ be the geodesic
                                ray of $T$ originating at $x$. Let
                                $b_{\overline{\beta}}$ be the Buseman function associated to 
                                $\beta$. Then $b_{\overline{\beta}}(\overline{y}) = b_{\beta}(y)$.
        \end{enumerate}

\end{definition}

We now need to show that these ``comparison triangles'' really do have the desired 
comparison properties.

\begin{lemma}
        Let $M$ be a CAT($-k$) space and let $T\subset M$ be a triangle with one 
        ideal vertex. Let $\alpha$ be the geodesic arc in $T$ from vertex $p$ 
        to vertex $q$. Let $\beta \sim \gamma$ be geodesic rays from 
        $p $ and $q$ respectively. Let $x, y\in T$ with $y\in\gamma$. 
        Let $T'$ be a comparison triangle for $T$. Let $x'$ and $y'$ be 
        the points of $T'$ corresponding to $x,y \in T$.
        Then $d( x,y ) \leq d( x', y' )$.
\end{lemma}
\begin{figure}
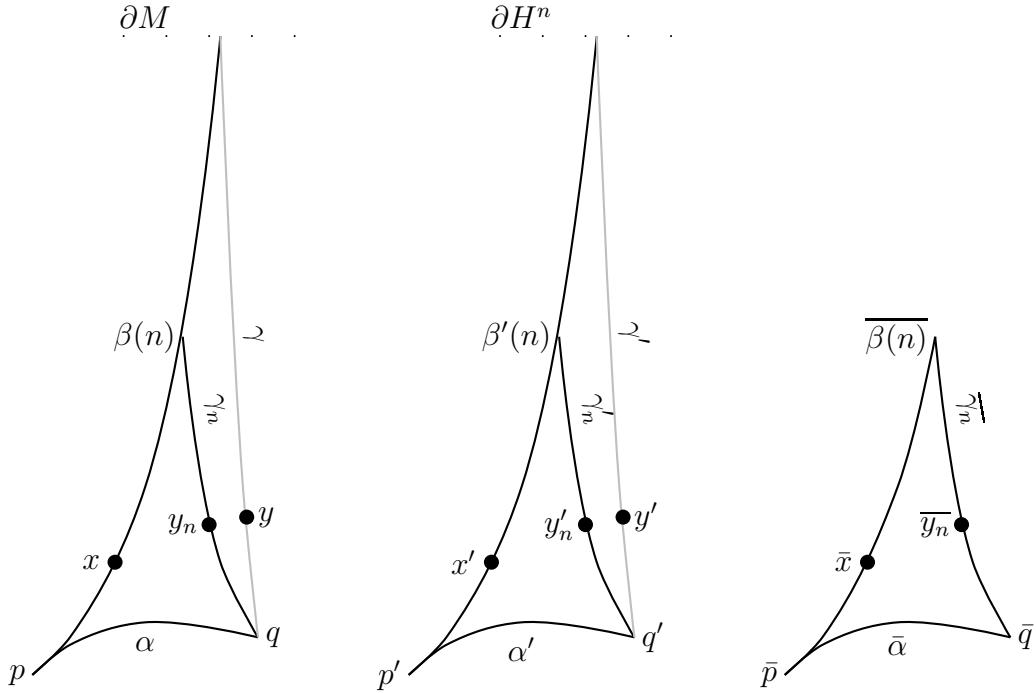

\caption
{
Illustration of the CAT($-k$) property for triangles with an ideal vertex.
The sequence $T_n$ converges to the triangle $T$ (left). There is a sequence
of triangles $T_n'$ converging to the comparison triangle $T'$ for $T$ (center).
The triangle $T_n'$ can be approximated by $\overline{T_n}$ (right), 
which is the comparison triangle for $T_n$.
}
\begin{center}

\input three_triangles.tex
\end{center}
\end{figure}

\begin{proof}
        We construct some sequences of triangles. Let $T_n$ be a triangle with 
        sides $\alpha , \beta\vert_{ [ 0,n ] } $ and the arc $\gamma_n$ from
         $q$ to $\beta(n)$. Let $d_n$ = $\length(\gamma_n)$.
         For sufficiently large $n$, we can define a point $y_n\in T_n$
         as the unique point on $\gamma_n$ such that $d(y_n , q ) = d (y, q )$.
         Let $\ell$ be the length of $\alpha$.
         Let $\overline{T_n} \subset \hyp^2$ be a comparison triangle for $T_n$.
         Let $T_n '\subset \hyp^2$ be a triangle derived from $T_n$ in the following 
         sense: choose a geodesic ray $\beta'\in \hyp^2$ originating at a point 
         $p'$, and a geodesic arc 
         $\alpha'$ from $p' $ to $q'$ of length $\ell$ such that 
         $b_{\beta'}(q') =b_{\beta}(q)$. Let $\gamma'_n$ be a geodesic arc 
         from $q' $ to $\beta'(n)$. For large enough $n$, we can define $y_n'$
         as the unique point on $\gamma'(n)$ 
         satisfying $d(y_n' , q' ) = d( y_n ,q ) = d(y,q )$.
         Let $d' = \length ( \gamma'_n )$.
         Observe that $T_n'$ converges point-wise to $T'$ ( hence the unconventional
         notation for the comparison triangle $T'$ ) 

         We will prove the following facts:
                For any $\epsilon>0, \delta > 0 $ , 
                for any points $x\in T, y\in \gamma$ , there exists 
                and $N_1(x,y,\epsilon)$ and $N_2(x,y,\delta)$ 
                such that the following are true:
         \begin{enumerate}
                \item for all $n>N_2$, 
                        $| d_{\hyp^2} ( \overline{x} , \overline{y_n} ) - d_{\hyp^2}( x' , y'_n ) | < \delta    $
                \item for all $n>N_1$, $d_{\hyp^2}( y' , y_n' )  < \epsilon     $
                \item for all $n>N_1$, $d_{M}( y , y_n )  < \epsilon    $
         \end{enumerate}
         Assuming this is true, the result follows:
                \begin{align*}
                        d( x,y )        &\leq d_M ( x,y_n ) + d_M ( y_n,y ) \\
                                                &\leq d_M ( x,y_n ) + \epsilon \\
                                                &\leq d_{\hyp^2} ( \overline{x}, \overline{y_n} ) + \epsilon  \\
                                                &\leq d_{\hyp^2}( x', y_n' ) + \epsilon + \delta \\
                                                &\leq d_{\hyp^2}( x', y' ) + d_{\hyp^2}( y', y_n' )+ \epsilon + \delta \\
                                                &\leq d_{\hyp^2}( x', y' ) + 2\epsilon + \delta \\
                \end{align*}
                Since $\delta$ and $\epsilon$ can be made arbitrarily small, 
                we see that $d_M (x,y) \leq d_{\hyp^2}( x', y' )$.
        So it remains to show that we can indeed find such an $N_1$ and $N_2$.
        The existence of $N_1(\epsilon)$ was established in \ref{cat0_sublemma}.
        So we need to show that there is an $N_2(\delta)$. 
        We argue this as follows: observe that 
        $\length{ \overline{ \gamma_n } } = \length( \gamma_n )  = d_n $.
        By the definition of the Buseman function $b_\beta$, it follows that
        $d_n = n + b_\beta ( q ) + a_n $ where $a_n$ is a constant that 
        goes to $0$ as $n$ approaches $\infty$. Similarly, 
        $d_n' = n + b_{\beta'} ( q' ) + a'_n  = n + b_\beta( q ) + a'_n$.
        So $d_n' -d_n = a'_n - a_n$, ie $d_n'$ approximates $d_n$ for large enough $n$.
        This means that $\overline{T_n}$ and $T_n'$ can be made arbitrarily close to
        identical by choosing $n$ large enough. 
\end{proof}

We now need to prove some basic facts about hyperbolic trigonometry
\begin{lemma}\label{cat0_rightangle}
        Let $T$ be a hyperbolic right angled triangle. Let $t$ be the length 
        of the side opposite the right angle, let $v$ be a vertex, $u>1$ the length 
        of the side opposite $v$ and $\theta $ the angle at $v$. Then 
        $$t-u \geq \log \left( \frac{1}{ 2\sin \theta } \right) $$
\end{lemma}
\begin{proof}
        By the hyperbolic law of sines, 
        $$ \frac{ \sinh u }{ \sin \theta } = \frac{ \sinh  t }{ 1 } $$
        Since $u > 1 $, we have $e^u / 2 > e^{-u}$, so 
        \begin{align*}
                \frac{ e^u }{ 4 }
                        \leq \frac{ e^u - e^{ -u } }{ 2 }   
                        = \frac{ e^t - e^{ -t } } { 2 } \sin \theta 
                 \leq \frac{ e^t  } { 2 } \sin \theta 
        \end{align*}
        So we conclude that 
        $$e^u \leq 2e^t \sin \theta$$
        Rearranging, we get the required inequality.
\end{proof}

\begin{lemma}\label{ideal_rightangle}
        For any  $C\in \R$, there exists $L(C)\in \R$  with the following property:
        Let $T\subset\hyp^n$ be a triangle with one ideal vertex. Let $x$ and $y$ be the 
        vertices of $T$. Let $\alpha$ be a geodesic arc of length $\ell > L$ 
        from $x$ to $y$. Let $\beta$ and $\gamma$ be sides of $T$
        from ( resp. ) $x$ and $y$. Suppose further that $b_{\beta}(y) = 0$.
        Let $z$ be the midpoint of $\alpha$. Then 
        $$b_{\beta}(z) \leq -C$$
\end{lemma}

\begin{proof}

        Observe that the two vertex angles of $T$ are identical. Let $\theta$
        be this angle. There is a right angled triangle $T'$ whose vertices are 
        $x$, $z$ , and the ideal vertex of $T$. Let $\xi$ be a ray from 
        $z$ such that $\xi\sim \beta$. 

        \begin{gather*}
                \cosh( \ell ) = \frac{ 1 + \cos^2 \theta }{ \sin^2 \theta }\\
                e^\ell / 2 \leq 
                \frac{ 1 + \cos^2 \theta }{ \sin^2 \theta } \leq 
                \frac{ 2 }{ \sin^2 \theta } \\
                \Rightarrow \sin^2 \theta \leq 4e^{-\ell }
        \end{gather*}

        \begin{figure}[h]
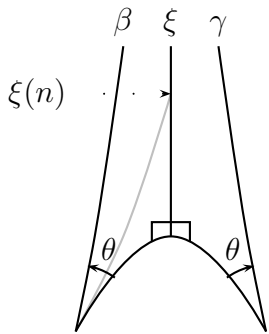

                \begin{center}
                \caption{An Ideal ``Isosceles'' triangle}
                \input isosceles.tex
                \end{center}
        \end{figure}

        This shows that we can make $\theta$ arbitrarily small by choosing 
        $L$ large enough. 

        Now consider the triangle $T_n$ whose vertices are $x,z$ and $\xi(n)$.
        $T_n$ converges point-wise to $T'$ ( see \ref{cat0_sublemma} ). 
        Moreover, by \ref{cat0_rightangle}, 
        $d( x, \xi(n) ) - n \geq  \log \left( \frac{1}{ 2\sin\theta } \right) $.
        So it also follows that 
        $$
                b_{\xi}(x) 
                = 
                        \lim_{n\rightarrow\infty } d( x , \xi(n) ) - x 
                \geq 
                        \log \left( \frac{1}{ 2\sin\theta } \right )
        $$
        But 
        $
          b_{\beta}(z) = - b_{\xi}(x)
                \leq 
        - \log \left( \frac{1}{ 2\sin\theta } \right )
        $.
        So the result holds for $L = 2C + \log 16 $

\end{proof}

\begin{lemma}\label{buseman_comparison}
        CAT($-k$) space.
        Let $T\subset M$ be a triangle with exactly one ideal vertex, 
        and $\overline{T}$ a comparison 
        triangle for $T$.
        Let $\alpha$ be the geodesic arc between the vertices 
        of $T$. Let $\beta$ be one of the geodesic rays of $T$. 
        Then for any point $x$ on $\alpha$, 
        $$ b_\beta (x) \leq b_{\overline{\beta} } (\overline{x} ) $$
\end{lemma}
\begin{proof}
        It follows from the CAT($-k$) inequality that 
        $d(\beta(t) , x )        \leq d( \overline{ \beta(t)} , \overline{x}  ) $
        so 
        $$d(\beta(t) , x ) -t    \leq 
        d( \overline{ \beta(t) }, \overline{x}  ) - t $$
        Take limits of both sides as $t\rightarrow \infty$, 
        and we have the desired result. 
\end{proof}

\begin{proposition}[ Penetration Depth ]\label{penetration_depth}
        Let $D> 0 \in R$ be given. Then there exists a constant $a=a(D)>0 \in \R $
        such that the following is true:
        Let $\alpha: [ 0, a ]\rightarrow M$ be a geodesic arc such that 
        $b_\gamma(\alpha ( 0 ) ) = b_\gamma(\alpha ( a ) ) = 0 $ for some geodesic
        ray $\gamma$ corresponding to a horoball $B$. Then there exists 
        a point $z\in\alpha$ such that $b_\gamma( z ) \leq -D$
\end{proposition}
\begin{proof}
        This result follows easily from previous lemmas. 
        The \ref{ideal_rightangle}
        proves the result for $\hyp^n$. In the general case, 
        a comparison triangle ( with one ideal vertex ) in $\hyp^n$  is used.
        The result follows immediately from \ref{buseman_comparison}.
\end{proof}

%
%

\section{Bounded Coset Penetration}

\subsection{Introduction}

        This section is based on the work of Farb in \cite{farb1}, 
        in fact it's essentially a rework of section 4 of \cite{farb1}. The goal
        is to provide a slight generalization of (and corrections to) the
        proof of the bounded coset penetration (BCP) property. 
        This property plays an important role in our proof.

        We begin with some basic definitions

\begin{definition}

        A {\em Hadamard manifold }, $\tilde{M}$ 
        is a complete, simply connected Riemannian manifold
        with non-positive sectional curvatures. 
        We typically denote such manifolds 
        by $\tilde{M}$ because for our applications, the Hadamard manifold 
        $\tilde{M}$ we are interested in is a universal covers of some complete
        Riemannian manifold $M$, such that the sectional curvature $\kappa(M)$ 
        satisfies $- a^2 \leq \kappa(M) \leq - b^2 $, where $a,b\in \R$.
        We call such a manifold a {\em pinched Hadamard manifold }

%


        Let $\tilde{M}$ be a Hadamard manifold and suppose that a 
        group $G$ acts on $\tilde{M}$ and $M = \tilde{M}/G$
        is a complete, non-compact, finite volume manifold.
        We can choose 
        $G$-invariant set of horoballs such that any $\tilde{M}$ geodesic 
        intersecting two distinct horoballs has length at least 
        $\max( 1, \delta )$, where
        $\delta$ is  the Gromov hyperbolicity constant for $\tilde{M}$ and the 
        $G$ action on the horoballs has finitely many orbits. 
        Remove the interiors of the horoballs corresponding to those 
        horospheres to obtain a space $\Psi$ on which $G$ acts cocompactly.
        When we speak of horospheres in $\Psi$,
        it should be understood that we refer to the {\em boundary horospheres} 
        of $\Psi$. 
        
        We define a Riemannian pseudo metric $h_{i,j}$
        on $\Psi$ as follows: let $g_{i,j}$ be the Riemannian metric on $\tilde{M}$. Then 
        on $\text{interior} ( \Psi ) $, we define $h_{i,j} = g_{i,j}$ and on $\bdy{\Psi}$, $h_{i,j}=0$. 
        In other words, we are locally inheriting the $\tilde{M}$ metric on the 
        interior of $\Psi$, but adding the condition that $d(x,y)=0$ for any two points 
        in a horosphere $S\subset \bdy{\Psi}$.
        We call this metric the {\em electric metric}, and we call the pair 
        $( \Psi, h_{i,j} )$  the electric space. We denote the electric space by $\widehat{\Psi}$. 
        A distance function 
        $d_{\widehat{\Psi}}: \widehat{\Psi} \times \widehat{\Psi}\rightarrow [0,\infty)$ 
        can be induced from $h_{i,j}$ as follows: we compute $d_{\widehat{\Psi}}(x,y)$ 
        by taking the inf of the path length with respect to $h_{i,j}$ of all paths
        from $x$ to $y$. 
        We call the path length induced by $h_{i,j}$ \emph{electric length}.
        It is not clear that there is a path that realizes 
        this inf, and we won't need to use such a thing, we will be
        more interested in ``electric quasi-geodesics''
        (when we want to use geodesics, we will use geodesics in $\tilde{M}$).
        A {\em geodesic} is a locally length-minimizing path. Note that geodesics 
        in $\widehat{\Psi}$ can behave arbitrarily on the horospheres in $\Psi$. 
 
        A $(\lambda,\epsilon)$ quasi-geodesic in $\widehat{\Psi}$ 
        is a path $\gamma$ that has the property that any sub-path
        $\gamma'$ of $\gamma$ whose endpoints are electric distance 
        $k$ apart is of electric length at most $\lambda k + \epsilon$.
        We will refer to such paths as $(\lambda,\epsilon)$ 
        \emph{electric quasi-geodesics}.

        Let$\gamma$ be a path in $\Psi$. We will parameterize paths in the
        electric space by the path metric on $\Psi$ 
        (since we can hardly parameterize by arc-length!). We say that a 
        {\em strongly non-horospherical segment} of $\gamma$ is a segment 
        $\gamma_{ [t_1 , t_2 ] }$ such that $\gamma(t_1)\in S , \gamma(t_2) \in T$
        for some horospheres $S\neq T $, and $\gamma_{ (t_1 , t_2 ) } \subset 
        \Psi - \bdy {\Psi} $. A non-horospherical segment is defined in the same way, but
        we allow $S=T$. A {\em horospherical segment} of $\gamma$ is a maximal segment 
        that is contained in a horosophere $S$. a {\em quasi-horospherical} segment
        is a maximal segment that starts and ends inside some horosphere $S$, and 
        intersects no other horospheres. 
        Strongly non-horospherical segments are complementary ( modulo endpoints )
        to quasi-horospherical segments, and horospherical segments are complementary
        (again, modulo endpoints) to non-horospherical segments.

        We say that a path $\gamma$ {\em penetrates} $S$ at $x\in S$ if $x$ is the 
        initial point of a horospherical segment. We say that $\gamma$ 
        {\em greets} $S$ at $x$ if $x$ is the initial point of a quasi-horospherical 
        segment. If $\gamma$ penetrates $S$ at $x$ , but doesn't greet $S$ at $x$,
        we say that $\gamma$ {\em re-penetrates } $S$ at $x$. Note that this terminology 
        is slightly subtle: it is possible that a path could greet $S$ more than once.
        We do {\em not} call the second greeting a re-penetration. We say that 
        $\gamma$ {\em first greets } $S$ at $\gamma(t_0)$ if $\gamma(t_0)\not\in S $ 
        for $t<t_0$, $t_0 > 0 $ 
        and $\gamma(t_0)\in S $. We say that $\gamma$ 
        {\em permanently leaves } $S$ at $\gamma(t_0)$ if $\gamma(t_0)\in S $, and
        $\gamma(t) \not \in S $ for all $t>t_0$

\end{definition}

\begin{notation}

        In this subsection, we shift gears in terms of notation, because we
        are dealing entirely with geometry ( as opposed to group theory ). 
        We adopt the following conventions: $G=\pi_1 (M) $ is a group,
        $M$ is a negatively curved manifold 
        with $n$ cusp(s).  $\tilde{M}$ is a Hadamard ( ie simply connected,
        complete and non-positively curved ) Riemannian manifold on which 
        $G$ acts freely, with $n$ orbits of parabolic fixed points.
        $S$ will be used to denote a horosphere in $\Psi$.

\end{notation}

\begin{proposition}[\cite{heintze}]\label{heintze}

        Let $\tilde{M}$ be a Hadamard manifold, such that 
        $-b^2 \leq \kappa(\tilde{M} ) \leq -a^2 < 0$, $a,b\in \R$. Let $\gamma$
        be a geodesic tangent to a horosphere $S$ and let $p$ and $q$ 
        be the projections of the endpoints of $\gamma$ onto $S$.
        Then 
                $$2/b \leq d_S ( p,q ) \leq 2/a $$ 
        where $d_S$ denotes the path metric on $S$ induced by the metric 
        on $\tilde{M}$.

\end{proposition}

The proposition implicitly states that horospheres have tangent planes.
This deserves some clarification. Heintze's argument shows that Buseman
functions in $\tilde{M}$ are $\mathcal{C}^2$ and have $\mathcal{C}^1$ 
gradient. So it makes sense to speak of ``tangency''.

Note that proposition \ref{heintze} admits some generalization:
\begin{corollary}\label{generalised_heintze}
        Let $\tilde{M}$ be a Hadamard manifold, such that 
        $-b^2 \leq \kappa(\tilde{M} ) \leq -a^2 < 0$, $a,b\in \R$. Let $\gamma$
        be a geodesic that doesn't intersect a horosphere $S$ and let $p$ and $q$ 
        be the projections of the endpoints of $\gamma$ onto $S$.
        Then 
                $$d_S ( p,q ) \leq 2/a $$ 
        where $d_S$ denotes the path metric on $S$ induced by the metric 
        on $\tilde{M}$.
\end{corollary}
\begin{proof}
        Choose a horosphere $S'$ corresponding to the same boundary point as $S$ such that
        $\gamma$ is tangent to $S'$. Then the result applies to $S'$. We can then project $S'$
        onto $S$ by flowing along the projection lines between $\gamma$ and $S'$. This 
        projection is length decreasing, so 
        $\length ( \pi_S ( \gamma ) ) \leq \length ( \pi_{S'} ( \gamma ) ) \leq 2/a $
\end{proof}

\begin{corollary}[Projections of horospheres on horospheres are bounded \cite{farb1}]\label{farb4.2}

        Let $\tilde{M}$ be a Hadamard manifold such that 
        $-b^2 \leq \kappa(\tilde{M} ) \leq -a^2 < 0$, $a,b\in \R$. 
        Let $S_1$ and $S_2$ be disjoint horospheres based at different points
        of $\partial \tilde{M}$.
        Then the projection of $S_2$ onto $S_1$ has diameter at most $4/a + 2\delta$ 
        with respect to the metric $d_{S_1}$, where $\delta$ is the Gromov-hyperbolicity constant for $\tilde{M}$.

\end{corollary}

\begin{proof}
        Let $b:\tilde{M}\rightarrow \R$ be a Buseman function
        based at $S_1$ such that $b(S_1) = 0 $, and let $S_1'$ = 
        $b^{-1}( \delta )$ where $\delta$ is the Gromov hyperbolic constant
        for $\tilde{M}$ ( so the horoball corresponding to $S_1' $ contains $S_1 $ ).
        Let $\xi$ be  the geodesic  between
        the boundary points of $\tilde{M}$ corresponding to $S_1'$ and $S_2$.
        Let $w=\xi\cap S_1'$.
        It is easy to see that for any point $z$ on $S_2$, the line between 
        $z$ and $w$ intersects $S_1'$ only at $w$ ( inspect the triangle formed 
        by $\xi\cap S_2$, $w=\xi\cap S_1'$ and $z$. The angle at $\xi\cap S_2$
        is at least $\pi / 2 $, so the angle at $w$ is less than $\pi / 2 $ ).

        Let $x$ and $y$ be two points on $S_2$. 
 Let $\alpha$ be the geodesic containing 
        $w$ and $x$, parametrized so that $\alpha(0) = x $. 
        Let $\beta$ be a geodesic between $w$ and $y$ with $\beta(0) = y $. 
        If $\alpha$ intersects $S_1$ transversely, let $\alpha'$ be some geodesic 
        containing $x$ that is tangent to $S_1$, otherwise let $\alpha' = \alpha$.
        Similarly, define $\beta'$. 
        
        We need to bound $\pi_{S_1} ( \overline{xw}  )$.
        If $\alpha' = \alpha$, then it is bounded by $2/a$ as an immediate 
        consequence of \ref{heintze} or \ref{generalised_heintze}.
        If $\alpha' \neq \alpha$, then let $p=\alpha' \cap S_1$, and let 
        $q$ be the first point of $S_1$ lying on $\alpha$.
        The triangle $(x,p,q)$ is $\delta$-thin,  but the side $\overline{pq}$
        lies in the horoball corresponding to $S_1$ since horoballs are convex
        ( \ref{horoballs_convex} ). So any point on the arc $\overline{xw}$
        must lie within distance $\delta$ of some point on $\alpha'$. 
        Since horospherical projection is distance decreasing , 
        $\pi_{S_1} ( \overline{xw} ) $ stays within distance $\delta$ of 
        $\pi_{S_1} ( \alpha' )$. So $\pi_{S_1} ( \overline{xw} ) $ 
        has diameter at most $diam ( \pi_{S_1} ( \alpha' ) ) + \delta = 2/a + \delta$.
        
        The same argument shows that the diameter of $\pi_{S_1} ( \overline{yw} ) $ 
        is at most $2/a + \delta$. So the diameter of their union is at most 
        $4/a + 2\delta$, and in particular, the distance between $\pi_{S_1} ( y ) $ and $\pi_{S_1} ( x ) $ is at most $4/a + 2\delta$.
        
\end{proof}

\begin{definition}[Visual Size of a Horosphere]
        
        Let $\tilde{M}$ be a pinched Hadamard manifold, and let $S\subset\tilde{M}$ 
        be a horosphere, and let $B$ be the corresponding horoball.
        Let $\gamma$ be a bi-infinite geodesic in $\tilde{M} - B$.
        Given $s\in S$, we say that $\gamma$ {\em can be seen from } $s$ 
        if $ \overline{s\gamma(t)} \cap S = \{ s \} $ for some $t$.
        The {\em visual size}
        of $S$ with respect to $\gamma$ is the diameter of the set
        of points $s\in S$ that $\gamma$ can be seen from.      

        The {\em visual size} of the Horosphere $S$ is the supremum over all geodesics 
        $\gamma$ of the visual size of $S$ with respect to $\gamma$. 

\end{definition}

        It is clear that the diameter of $\pi_S (\gamma)$ is no more than the 
        visual size of $S$ provided that $\gamma\cap S = \emptyset $.

        Less obvious, but also true is that the diameter of the 
        projection of $S$ onto $\gamma$ is no more
        than the visual size of $S$ with respect to $\gamma$. This is because projection 
        onto $\gamma$ is length decreasing and the projection of $S$ onto $\gamma$
        is the same as the projection of the set of points on $S$ from which $\gamma$ 
        can be seen. 

\begin{lemma}[Horospheres are visually bounded]\label{visual}
        
        Let $\tilde{M}$ be a pinched Hadamard manifold with 
        $-a^2 \leq \kappa ( \tilde{M} ) \leq -b^2$
        where $a,b\in \R$. Then there exists $D>0$ that depends 
        only on $a$ and $b$ such
        that the visual size of any horosphere 
        $S \subset \tilde{M}$ is less than $D$.

\end{lemma}

\begin{proof}

        Let $S\subset \tilde{M}$ be a horosphere, and let $\gamma\subset \tilde{M} - S$
        be a geodesic.  Let $x= \gamma(t)$, and let $y\in S$ be a point such that
        the geodesic arc $\beta$ from  $x$ to $y$ has the property that
        $\beta\cap S = y$ ( ie $\gamma$ can be seen from $y$ ). 
        Let $z = \pi_S (x)$. Let $\alpha$ be
        the geodesic segment between $y$ and $z$.

        We use the geometry of the triangle $xyz$ to bound the length of $\alpha$.
        Observe that $\alpha=\overline{yz}$ lies in one component of $\tilde{M} -S $ while 
        $\overline{xz}$ and $\overline{xy}$ lie in the other component. Since $\tilde{M}$ is 
        Gromov-hyperbolic, 
        $\alpha$ lies in a $\delta$ neighborhood of $S$ where $\delta(\tilde{M})$ is the 
        hyperbolicity constant. By \ref{penetration_depth}, there is some $C(\delta)$ such that the any
        geodesic segment of length greater than $C$ between two points on $S$ does not
        stay in a $\delta$ neighborhood of $S$.  So $\alpha$ has length at most $C$.
        So $d( y, \pi_S ( \gamma )  ) \leq d ( y, \pi_S( x ) ) = d ( y, z ) \leq C$.

        We have shown that for any point $y\in S$ such that $\gamma$ can be seen from
        $y$, $d(y, \pi_S(\gamma ) ) \leq C$.
        By \ref{generalised_heintze}, $\pi_S ( \gamma ) $ has diameter at most $2/a$.
        So $\gamma$ has visual size at most $D=2/a + C(\delta(\tilde{M}))$. 
        Note that $D$ depends only on $\tilde{M}$.

\end{proof}

We need to prove the following:

\begin{lemma}[Projection of a Horosphere onto a penetrating geodesic]
        Let $\gamma$ be an $\tilde{M}$ geodesic that intersects some 
        horosphere $H$.  Let $B$ be the horoball corresponding to $H$.
        Let $\pi_\gamma : \tilde{M} \rightarrow \gamma$ be the 
        projection onto $\gamma$.
        Then there exists a $D\in\R$ depending only on $\tilde{M}$
        such that $\pi_\gamma (H) - B $ is in a $D$ 
        neighborhood of $H$.
\end{lemma}

\begin{proof}
        Since $\tilde{M}$ is a pinched Hadamard manifold, $\tilde{M}$ is
        $\delta$-hyperbolic for some $\delta\in \R$. Let $a\in\gamma$ be an element
        of the image of $\pi_\gamma ( H )$. Let $b\in H $ such that 
        $a = \pi_\gamma ( b )$, and let $c$ be the point of  $ \gamma\cap H $
        closest to $a$. Consider triangle 
        $(a,b,c)$. Since $\tilde{M}$ is $\delta$-hyperbolic, the side $\overline{ac}$
        is contained in a $\delta$ neighborhood of $\overline{ab} \cup \overline{bc}$.
        Let $p \in \overline {ac}$. If $p$ is in a $\delta$ neighborhood of 
        $\overline{bc}$, then $p$ is in a $\delta$-neighborhood of $B$ since 
        $\overline{bc}\subset B$. 

        If $p$ is not in a $\delta$-neighborhood of $\overline{bc}$, then 
        $p$ is in a $\delta$-neighborhood of $\overline{ab}$. 
        Let $\alpha$ be a geodesic line segment of length less than $\delta$ from $p$ 
        to $\overline{ab}$ and let $q = \overline{ab} \cap \alpha$. Since $\overline{ap}$
        meets $\overline{aq}$ at right angles, it follows that $\overline{pq}$ is 
        opposite a right angle, and hence is the longest side of triangle 
        $(a,p,q)$. 
        
\end{proof}

\begin{definition}
        Given an $\tilde{M}$  geodesic  $\gamma$, the {\em projection onto $\gamma$ relative
        to $\widehat{\Psi}$ } , denoted by $\pi_{\gamma, \widehat{\Psi} } $ is given by the following 
        rule: let $x\in \widehat{\Psi}$. Let $y$ be the projection of $x$ onto $\gamma$. Let 
        $\alpha$ be the $\tilde{M}$-geodesic from $x$ to $y$. If $y$ is in the 
        interior of some horosphere $S$, then define 
        $\pi_{\gamma,\widehat{\Psi}} (x) = \alpha\cap S $  otherwise, define $\pi_{\gamma,\widehat{\Psi}}(x) = y $
\end{definition}

\begin{lemma}[ Electric Quasi-geodesics electrically track $\tilde{M}$ geodesics ]\label{tracking}

        Given a pinched Hadamard manifold $\tilde{M}$, $\lambda \in \R, \epsilon \in \R$,
        there exists real constants $K = K(\tilde{M}, \lambda, \epsilon )$ and
        $L = L ( \tilde{M}, \lambda, \epsilon )$ such that for any
         electric $\lambda,\epsilon$ quasi-geodesic ray $\beta$ in $\tilde{M}$,
        if $\gamma$ is the $\tilde{M}$ geodesic with the same endpoints as $\beta$,
         then any subsegment of $\beta$ 
        outside $nbhd_{\widehat{\Psi}}(\gamma,K)$ has $\tilde{M}$-length at most $L$. 

\end{lemma}

\begin{proof}

        Choose $K$ so that $K\geq 1/a \log (2\lambda(D+1))$ where $D$ is the constant given by
        \ref{visual}.
        If $\beta$ stays inside $N=nbhd_{\widehat{\Psi}}(\gamma, K) $, we are done. So assume that
        $\beta$ leaves $N$. Let $\beta'$ be a maximal subsegment of $\beta$ contained in 
        $N'=\Psi - N$.
        Let $x=\beta'(0)$ and let $y$ be the endpoint of $\beta'$. 

        Suppose that $\beta'$ greets $q$ horospheres. 
        Then $q+1 \leq \length_{\widehat{\Psi}}\beta' $ since any path between distinct horospheres 
        has length at least 1.

        Let $x' = \pi_{\gamma, \widehat{\Psi}} (x) $ and $y' = \pi_{\gamma, \widehat{\Psi}} ( y ) $

        Now we have that
        \begin{align}\label{track_eqn}
                \length_{\widehat{\Psi}}\beta' 
                & \leq  \lambda \left( 
                        d_{\widehat{\Psi}}(x,x') + d_{\widehat{\Psi}}(x',y') + d_{\widehat{\Psi}}(y',y)
                        \right) + \epsilon 
        \end{align}

        The next step is to establish a bound on $d_{\widehat{\Psi}}(x,x')$ , and
        $d_{\widehat{\Psi}}(y,y')$. By symmetry, it suffices to obtain a bound on the former.
        The $\widehat{\Psi}$ distance between $x$ and $\gamma$ is at most $K$. This distance 
        is realized by some $\widehat{\Psi}$-geodesic $\alpha$. Let $z=\alpha\cap\gamma$. 
        Note that $z\in \Psi$.
        $$
                d_{\widehat{\Psi}}( x,x' ) 
                \leq d_{\widehat{\Psi}}( x , z ) + d_{\widehat{\Psi}}( z , x' )  
                \leq  K + d_{\widehat{\Psi}}(z, x' ) 
        $$
        So it suffices to bound $d_{\widehat{\Psi}}(z,x') $. We do this by bounding the
        $\widehat{\Psi} $-length of the  projection $\pi_{\gamma, \widehat{\Psi}}(\alpha)$ of 
        the path $\alpha$.
        First, observe that if $\alpha$ penetrates any horosphere $T$ that is penetrated
        by $\gamma$, then $\alpha$ does not leave $T$ ( because if $w\in \alpha$, and $w\in T$,
        then $w$ is $\hat{\Psi}$-distance $0$ from $\gamma$. ) Let $\alpha'$ be a maximal 
        initial segment of $\alpha$ that doesn't penetrate any such horosphere, and
        let $\alpha''$ be the remaining segment of $\alpha$.

        $ \alpha'' \subset T $, so 
        by ``projection of a horosphere onto a penetrating geodesic'', 
        $\pi_{\gamma,\widehat{\Psi}}(\alpha'')$ stays in a $2\delta$ 
        neighborhood of  $T$ where $\delta$ is the Gromov hyperbolicity constant of
        $\tilde{M}$. Hence the $\widehat{\Psi}$ length of $\pi_{ \gamma, \widehat{\Psi} }( \alpha'' )$
        is at most $2\delta$.

        So it remains to bound $\length_{\widehat{\Psi}} ( \pi_{\gamma, \widehat{\Psi} } ( \alpha' ) ) $.
        Let $U$ be a horosphere penetrated by a horospherical segment of $\alpha'$.
        Then $U$ is not penetrated by $\gamma$. 
        So $\pi_\gamma (U) $ has diameter at most $D$.
        Let $\eta$ be the (maximal) quasi-horospherical segment of $\alpha'$ 
        corresponding with $U$. Recall that $\eta$ has $\widehat{\Psi}$-length of 
        at most $\epsilon$. Since $\pi_{\gamma,\widehat{\Psi}}$ is length decreasing,
        $\pi_{\gamma,\widehat{\Psi}} ( \eta ) $ stays within an $\epsilon$ radius
        of $\pi_{\gamma,\widehat{\Psi}}( U ) $. 
        So $\pi_{\gamma , \widehat{\Psi}} ( \eta )  $ has a diameter of at most $D + \epsilon$.

        The projections onto $\gamma$ of the strongly non-horospherical
        segments have lengths totalling no more than $K$. There are at most $K$ strongly 
        non-horospherical segments and $K$ quasi horospherical segments, so the projection
        of these onto $\gamma$ has length at most $K( 1 + D + \epsilon ) $.

        We have shown that $\length ( \pi_{\gamma, \widehat{\Psi}} ( \alpha''  ) ) \leq 2\delta$
        and $\length ( \pi_{\gamma, \widehat{\Psi}} ( \alpha' ) ) \leq K ( 1 + D + \epsilon )$.
        So 
        $$d_{\widehat{\Psi}} ( z, x' ) \leq K ( 1 + D + \epsilon ) + 2\delta $$
        Since $d_{\widehat{\Psi}}( x,z) \leq K $, it follows that 
        \begin{align}\label{proj_eqn}
        d_{\widehat{\Psi}} ( x,x' ) \leq K ( 2 + D + \epsilon ) + 2 \delta
        \end{align}
        Let $K' = K ( 2 + D + \epsilon ) + 2 \delta$.
        It follows from \ref{track_eqn} and \ref{proj_eqn} that

        \begin{align*}
                \length_{\widehat{\Psi}}(\beta' ) 
                & \leq \lambda \left(
                        d_{\widehat{\Psi}}(x,x') + d_{\widehat{\Psi}}(x',y') + d_{\widehat{\Psi}}(y',y)
                        \right) + \epsilon \\
                & \leq \lambda \left( 
                        K' + 
                        e^{-aK} ( \length_{\widehat{\Psi}} (\beta') + Dq )
                        + K' 
                        \right)+\epsilon \\
                & \leq \lambda \left(
                        2K' +
                        e^{-a ( 1/a \log ( 2\lambda(D+1) ) ) }
                        ( \length_{\widehat{\Psi}} (\beta') + D ( \length_{\widehat{\Psi}}(\beta') - 1 ) ) 
                        \right )+\epsilon \\
                & = \lambda \left(
                        2K' +
                        \frac{1}{2\lambda(D+1)}
                        ( \length_{\widehat{\Psi}} (\beta') + D ( \length_{\widehat{\Psi}}(\beta') - 1 ) ) 
                        \right )+\epsilon \\
                & = \lambda \left(
                        2K' + 
                        \frac { 1 } {2\lambda(D+1)}
                        ( ( D + 1 ) \length_{\widehat{\Psi}} ( \beta' )  - D ) 
                        \right )+\epsilon \\
                & \leq \lambda \left(
                        2K' + \frac{ \length_{\widehat{\Psi}} (\beta' ) }{2\lambda }
                        \right )+\epsilon \\
                & = 2K' \lambda + \length_{\widehat{\Psi}}(\beta' ) / 2+\epsilon \\
        \end{align*}

        Moving the $\length_{\widehat{\Psi}}$ term to the left hand side, we get
        $$
                \frac
                {
                        \length_{\widehat{\Psi}}(\beta') 
                }{
                        2
                }
                \leq 2K' \lambda + \epsilon
        $$      
        so  it follows that 
        $$\length_{\widehat{\Psi}}(\beta') \leq 4K' \lambda + 2\epsilon = 4K\lambda( 2 + D + \epsilon ) + 8\lambda\delta + 2\epsilon$$
        This completes the proof, with $L(K,\lambda,\epsilon,\delta) = 4K\lambda( 2 + D + \epsilon ) + 8\lambda\delta + 2\epsilon$.

\end{proof}

\begin{lemma}\label{neighbourhoods}
        There exists $D(\lambda,\epsilon) \in \R$  with the following property:
        Let $\beta$ be an electric $(\lambda,\epsilon) $ quasi-geodesic  and let $S$ be
        a horosphere where $\beta\cap S = \emptyset$. Then $\pi_S(\beta)$ has diameter
        (in the $\tilde{M}$-metric on $S$) of at most $D\length_{\widehat{\Psi}}(\beta)$.
\end{lemma}

\begin{proof}
        Recall that all horospheres in $\widehat{\Psi}$ are a distance of at least 1 apart.
        Each strongly non-horospherical segment of $\beta$ is a 
        $(\lambda,\epsilon)$-quasi-geodesic in $\tilde{M}$. 
        
        Let $\beta'$ be such a segment.
        Since $\pi_S$ is length decreasing on $\beta'$, 
        it follows that $\length( \pi_S(\beta') ) \leq \length ( \beta' ) $.

        So we investigate the quasi-horospherical segments of $\beta$. 
        Let $\beta''$ be such a segment, and let $T$ be the corresponding horosphere. 
        $\beta''$ has $\widehat{\Psi}$ length of 
        no more than $\epsilon$. Since $\pi_S$ is length decreasing, 
        $\pi_S( \beta'' ) $  stays within and $\epsilon$ neighborhood 
        of $\pi_S ( T) $.
        So the diameter  of $\pi_S(\beta'')$ is at most $E + \epsilon$ 
        where $E$ is the constant given by \ref{visual}

        So $\sum ( diam ( \pi_S ( \beta' ) ) ) + \sum ( diam ( \pi_S ( \beta'' ) ) ) $
        where $\beta'$ and $\beta''$ range over the strongly non-horospherical and
        quasi-horospherical segments of $\beta$ (respectively) is bounded by 
        $\length_{\widehat{\Psi}}(\beta) + q ( E + \epsilon ) $ where $q$ is the number of 
        horospheres that $\beta$ greets. But $q\leq \length_{\widehat{\Psi}}(\beta ) $.
        So we conclude that $\pi_S ( \beta )$ has diameter at most 
        $\length_{\widehat{\Psi}}(\beta) ( 1 + E + \epsilon ) $.

\end{proof}

\begin{lemma}\label{quasigeodesic projections}
        Let $\beta$ be a $(\lambda,\epsilon)$-quasi-geodesic in $\widehat{\Psi}$ from $x$ to $y$. Let 
        $S$ be a horosphere such that $\beta(0)\not\in S$. If $\beta\cap S\neq \emptyset $,
        let $t_0 = \min \{ t: \beta(t) \in S \} $ and $t_1 = \max \{ t: \beta(t) \in S \} $.
        If $\beta\cap S = \emptyset$, let $t_0 = t_1 = \infty $.
        Let $x'=\beta(t_0)$ and $y'=\beta(t_1)$. Then there exists a constant $D\in\R$ such 
        that the following are true: 
        
        \begin{enumerate}
                 \item $ \pi_S ( \beta_{ [ 0,t_0 ] } ) $ has diameter at most $D$
                 \item $ \pi_S ( \beta_{ [ t_1, \infty )   
                 } )  $ 
                 has diameter at most $D$
        \end{enumerate}

\end{lemma}

\begin{proof}
        First, we prove (1). 
        By \ref{tracking}, there exists a constant $D$ such that any quasi-geodesic
        stays within a $\hat{\Psi}$-distance of at most $D$ from the $\tilde{M}$
        geodesic with the same endpoints.
        By \ref{neighbourhoods}, there is a constant $E$ such that for any arc $\eta$ that
        doesn't intersect $S$, $\pi_S(\eta )$ has diameter no more than $E( \length(\eta) )$. 
        Let $\delta$ be the Gromov hyperbolicity constant for $\tilde{M}$. If $\delta > E$,
        we will choose $E=\delta$.
        Let $\gamma$ be the $\tilde{M}$ geodesic with the same endpoints as $\beta$.
        Let $\gamma'$ be a $\tilde{M}$ geodesic such that $\gamma'(0) = \gamma(0)$ and 
        $\gamma'$ is tangential to $S$.
        Let $B$ be the $E$ neighborhood of $S$ in $\hat{\Psi}$.
        Let $\beta'$ be a maximal connected 
        subset of $\beta$ with endpoints in the closure
        of $B$.  Since $\beta'$ is a $(\lambda,\epsilon)$ quasi-geodesic, $\beta'$ has 
        electric length at most $\lambda D + \epsilon $.
        So $\pi_S ( \beta' )  $ has a diameter of at most
        $E ( \lambda D + \epsilon ) $.

        We show that the rest of $\beta$ is in a bounded radius of $\pi_S(\gamma ) $.
        Let $y$ by a point on $\beta$. Then $y$ is an $\hat{\Psi}$-distance of no more 
        than $D$ from $\gamma$. This distance is realized by some electric
        geodesic $\xi_y$. The length of $\pi_S ( \xi_y )$ is no more than 
        $E D$.
        Let $q$ be the first point on $\gamma$ that intersects $S$.
        Any point $p$ on the arc between $\gamma(0)$ and $q$  is within distance
        $\delta$ of some point on $\gamma'$ : consider the delta thin triangle whose
        vertices are $q' = \gamma' \cap S$, $\gamma(0)$, and $q$. $p$ is within distance
        $\delta$ of either $\gamma'$ or the arc between $q$ and $q'$. But this arc is 
        contained in the horoball corresponding to $S$ since horoballs in $\tilde{M}$ are 
        convex.
        We have shown that any point on $\beta_{[0,t_0]}$ 
        is within distance at most $D+\delta$ of 
        $\gamma'$, and the diameter of $\pi_S ( \gamma' )$ is at most $2/a$ ( recall that $-a^2$
        is the upper bound for the curvature of $\tilde{M}$ ).

        Part (2) uses exactly the same argument as part (1) (you can apply precisely 
        the same argument to the path $\beta^{-1}$.)

\end{proof}

\begin{lemma}[Bounded Horosphere Penetration]
        Let $G$ be a group acting on a pinched Hadamard manifold $\tilde{M}$.
        Let $\Psi$ be a defined as previously, recall that
        $\Psi$ is obtained by removing the interiors of a set of $G$-invariant 
        horoballs of $\tilde{M}$ such that any two horoballs are distance
        $\max(1,\delta)$ apart, where $\delta$ is the Gromov-hyperbolicity
        constant for $\tilde{M}$. Let $\widehat{\Psi}$ consist of the points
        in $\Psi$ with the electric metric.
        Let $\alpha$ and $\beta$ be $\lambda$ quasi-geodesics from $x$ to $y$ 
        in $\widehat{\Psi}$. 
        Let $S$ be a horosphere in $\Psi$. Then there exists a 
        constant $E\in \R$ such that
        the following hold: 
        \begin{enumerate}
                \item   Suppose $\alpha$ first greets $S$ at $\alpha (s_0)$ and $\beta$
                first greets $S$ at $\beta (t_0)$. Suppose that $\alpha$ and $\beta$ permanently
                leave $S$ at $\alpha(s_1)$ and $\beta(t_1)$
                Then $d_S ( \alpha (s_0) , \beta(t_0) ) < E$ and 
                $d_S ( \alpha (s_1) , \beta(t_1) ) < E$
                \item   Suppose $\alpha$ greets $S$ at $s_0$ and permanently leaves $S$ 
                at $\alpha(s_1)$. Suppose that $\beta$ doesn't greet $S$. Then 
                $d_{\tilde{M} } ( \alpha(s_0), \alpha(s_1) ) < E$.
        \end{enumerate}
\end{lemma}

\begin{proof}
        First, we prove (1). If $x\in S$, we are done. So assume that $x\not\in S$.
        Let $D$ be the constant given by \ref{quasigeodesic projections}. So 
        $diam_S ( \pi_S ( \alpha_{ [ 0, s_0 ] } ) ) \leq D $ and 
        $diam_S ( \pi_S ( \beta_{ [ 0, t_0 ] } ) ) \leq D $. 
        Since $N=\alpha_{ [0,s_0 ] } \cup \beta_{ [ 0,t_0 ] } $ is a connected set,
        so $\pi_S (N ) $ is connected and has diameter at most $2D$. This completes 
        the proof of the first assertion in (1).

        The second assertion in (1) follows by applying the above argument to 
        $\alpha^{-1} $ and $\beta^{-1}$.

        Now we prove (2). By \ref{quasigeodesic projections}, $\pi_S(\beta)$ has 
        $S$-diameter at most $D $. Similarly, $\pi_S( \alpha_{ [0,s_0 ] } )$ 
        and $\pi_S ( \alpha_{ [s_1, \infty 
        )  } )$ 
        also have diameter at most $D$
        Since the set 
        $$\pi_S( \alpha_{ [0,s_0 ] } ) \cup \pi_S ( \alpha_{ [s_1, \infty ) } )\cup \pi_S(\beta) $$
        is connected, it has diameter at most $3D$. It also includes the points 
        $\pi_S(\alpha(s_0)) $ and $\pi_S(\alpha(s_1))$, so the $S$-distance between these 
        points is at most $3D$.
        This proves part (2).

        Hence the result is true for $E=3D$.

\end{proof}

\subsection{Relatively Hyperbolic Groups}
Given a group $G = \langle X\vert R\rangle$,  with Cayley graph $\Gamma_G $ , and a 
collection of subgroups $P_i = \langle  y_{i,j} \rangle \subset G$ , the 
{\em coned off Cayley graph  over $P$ }  is the graph obtained by adding a
$v_{g,i}$
vertex for each 
coset $gP_i$ , and adding an edge of length $1/2$ from $v_{g,i}$ to the vertex corresponding to $gP_i$.
For convenience, we usually omit $P$ from the notation.
A group $G$ is {\em hyperbolic relative to a collection of groups 
        $\{ P_i  = \langle y_{i,j} \rangle : i\in I \}$ } 
if the coned off Cayley graph 
is a delta-hyperbolic metric space. 

For most of this discourse, we discard this coned off Cayley graph. 
The main problem with this space is that the metric is ``bad'' ,
 our subsequent results depend on a pseudo-metric $d$ 
that gives $d(gp_1 , gp_2) = 0 $ where $p_i\in P$ and $P$ is a parabolic subgroup.

We define a notion of \emph{electric length}:
\begin{definition}
        Given a collection of groups subgroups of $G=\langle X\vert Y\rangle $, 
        $P_i = \langle y_{i,j} $,  and assume that $y_{i,j} \in X$.
        Let $\mathcal{Y}$ be the union of all the $y_{i,j}$.
        We define the \emph{electric word length} with respect to the 
        groups $P_i$ as
        $$\ell : X^* \rightarrow \Z , \ell(x_1, \dots, x_n) = 
        \text{cardinality} \{ x_j \vert x_j \not\in \mathcal{Y} \}
        $$
        In other words, electric length is a length function that assigns a length
        of $0$ to the generators $y_{i,j}$.
\end{definition}

The  pseudo-metric space we use is the {\em electric Cayley graph} $\widehat{\Gamma_G} $ 
which consists of the Cayley
graph with a pseudo-metric induced by electric length on $G$.

Given a group $G$ acting on a pinched Hadamard manifold,
with some choice $y_{i,j}$ of the 
cusp subgroups $P_i$, for any word $w$, a {\em parabolic segment} 
of $w$ is a maximal sub-word of the form $y_{i,j_1} y_{i,j_2} \dots y_{i,j_k} $
A word $w$ is said to {\em penetrate } a coset, $gP_i$ if $w(t_0)\in gP_i$ 
for some $t_0$. $w$ {\em leaves} a coset $gP_i$ at $t_0\in \N$ if 
$w(t_0) \in gP_i$ and $w(t_0  + 1 ) \not\in gP_i$ 
$w$ is said to {\em backtrack} if it penetrates some coset
more than once.

\begin{definition}{Bounded Coset Penetration} \\
        A pair consisting of a groupoid $G$ and a finite collection of subgroups $\{
        P_i : i = 1 , 2 , \dots , n\} $ satisfies the {\em bounded coset
        penetration property } if for any $k\geq 1$, there is a constant $c(k)>0$
        such that if $u$ and $v$ are $(\lambda,\epsilon)$ electric quasi-geodesics 
        and $d_{\Gamma}(\bar{u} ,\bar{v}) \leq 1 $ , then the
        following are true
        
        \begin{enumerate}
                \item   if $u$ penetrates a coset $gP_i$ and $v$ does not penetrate $gP_i$, then $u$ 
                                travels a $\Gamma$-distance of at most $c$ in $gP_i$. 
                \item   If both $u$ and $v$ penetrate a coset $gP_i$ , then the 
                                vertices of $\Gamma$ at which $u$ and $v$ first enter 
                                $gP$ lie a $\Gamma$-distance of at most $c$ from each other. The 
                                same is true for the vertices of $u$ and $v$ where $u$ 
                                and $v$ leave $gP_i$.
        \end{enumerate}

\end{definition}

\begin{theorem}[Bounded Coset Penetration]\label{BCP}
        Let $G$ be a group acting on a pinched Hadamard manifold 
   $\tilde{M}$ such that $M=\tilde{M}/G$ is a complete,
non-compact, finite volume manifold.
        Let $\{P_i\}$ be the cusp subgroups. 
        Then $G$ has the bounded coset penetration property with respect 
        to $P_1 , \dots , P_n$.
\end{theorem}

\begin{proof}
        First, we address the (relatively easy) one cusp case.
        
        We use a quasi isometric embedding $f: \hat{\Gamma_G}\rightarrow \widehat{\Psi}$.
        In the case where the action of $G$ on $\Psi$ is cocompact, this map will in fact
        be a quasi-isometry.  Let $G=< X\vert R >$. 
        Let $P=<Y>$ where $Y\subset X$.
        
        Choose a point $p$ on 
        a horosphere $S$. For each $g\in G$, denote by $v_g$ the corresponding vertex 
        of $\Gamma_G$. We define $f(v_g) = g\cdot p$. 
        For each edge in the Cayley graph $(1,a)$ where $a\in X$, 
         we can join the points $p$ and $a\cdot p$ by a path $\beta_a$ that 
         \begin{enumerate}
                \item  intersects $\bdy \Psi $ only at its endpoints if $a\not\in Y$ 
                \item lies entirely in $\bdy \Psi$ if $a\in Y$.
        \end{enumerate}
        We can then translate
        these edges around the Cayley graph (ie the edge $(g,ga)$ is mapped to the edge
        $g \beta_a $.) This is a quasi-isometric map (since the $G$-action is by isometries) 
        If the action of $G$ on $\widehat{\Psi}$ is cocompact 
        it follows that some $\epsilon$ neighborhood with respect to the $\tilde{M}$ metric 
        contains $\Psi$.    
        
        $f$ maps  all of the vertices of $\Gamma_G$ 
        to  points on horospheres. Moreover, it maps the generators $Y$ of $P$ 
        to the horospheres. In particular, if $g_1 , g_2 $ are in the same coset 
        $g_P$, then their images can be joined by a path lying in the horosphere 
        $gS$. So 
        This implies that it is a quasi-isometric embedding of the electric Cayley graph 
        $\widehat{\Gamma_G}$ into the space $\widehat{\Psi}$.
        
        Let $\gamma$ be some $(\lambda, \epsilon )$ electric quasi-geodesic that penetrates 
        a coset $gP$ and
        let $\beta$ be a $(\lambda, \epsilon) $ electric quasi-geodesic that doesn't.
        Suppose $\gamma(t_1)$ is the point at which $\gamma$ first enters $P$ and
        $\gamma(t_2)$ is the last point at which $\gamma$ leaves $P$. Then 
        $f(\gamma(t_1) ) $ is the first 
        point where $f(\gamma)$ enters the horosphere corresponding with 
        the coset $gP$, and $f(\gamma(t_2))$ is the last point at which $f(\gamma)$ leaves
        this horosphere. By the bounded horosphere penetration lemma, $f(\gamma(t_1))$ 
        and $f(\gamma(t_2))$ are distance at most $D$ where $D$ is the constant given
        by the result. So $d(\gamma(t_1) ,\gamma(t_2 ) ) \leq \lambda D + \epsilon $ since 
        $f$ is a quasi-isometric map. If $\beta$ and $\gamma$ both penetrate a given coset,
        we can use an analogous argument to show that the initial penetration points
        $\beta(s_1) $ and $\gamma(t_1)$ are a bounded distance apart as are the points 
        $\beta(s_2) $ and $\gamma(t_2)$ which are the first points where  $\beta^{-1}$ 
        and $\gamma^{-1}$  penetrate $gP$.

        This completes the argument for the single subgroup case. The case  where there are 
        several subgroups is a little more complex. First, there is an apparent difficulty
        embedding the Cayley graph: embedding the Cayley graph requires one to choose 
        a $G$ orbit of the horosphere. Each such embedding maps the cosets $gP_i$ of 
        a ``preferred'' parabolic subgroup $P_i$ to horospheres. The problem is that 
        any given a different subgroup $P_j$, $f$ does not map $P_j$ to
        a single horosphere. 

        So the solution is to ``choose'' all of them. 
        Let $G=<X\vert R>$ and let $Y_i$ be the generating set of $P_i$. 
        Let $Y_i \subset X$.
        For each horosphere orbit represented by a horosphere $S_i$, 
        we construct an embedding in the same manner
        as for the single coset case.
        For each $G$ orbit, we choose a base-point $p_i\in S_i$, 
        and we connect each two 
        points $p_i, p_j$ by an path $e_{i,j}$. 

        The  result is the Cayley graph of a 
        groupoid $\tilde{G}$
        whose morphisms are of the 
        form 
        $$(g, i, j )\in G \times \{1 , \dots n \} \times \{ 1, \dots , n \} $$
        The morphisms are composed as follows:
        $$
                (g, j_1, k ) \circ ( h, i, j_2 ) = ( gh ,  i,k  ) $$
                if $j_1 = j_2 $, otherwise $g$ and $h$ aren't composable.

        We denote morphisms of the form $(g, i, i) $ by $(g, i )$  and morphisms 
        of the form $( 1_G , i ,j ) $ by $m_{i,j}$. Note that these two 
        classes of morphisms generate the groupoid $\tilde{G}$.
        There are some $\tilde{G}$ cosets that are important: 
        the coset of morphisms of the form $( g, i )( p, i ) $ where 
        $g\in G , p \in P_i$. We will denote this by $(g, i ) ( P_i , i ) $.
                
        An application of the argument used in the single coset case 
        shows that $\tilde{G}$
        enjoys a property analogous to
        the bounded penetration property with respect to the groups 
        $( P_i , i )$. 
        First, we define an electric pseudo-metric on the Cayley 
        graph $\Gamma_{\tilde{G}}$ by using the edge path metric, but counting
        any edge that lies on a horosphere as length 0.

        If $\beta$ and $\gamma$ are electric quasi-geodesics in $\Gamma_{\tilde{G}}$,
        we can use the same arguments as above to show the following:
                \begin{enumerate}
                        \item If $\beta$ penetrates some groupoid coset $( g, i ) ( P_i , i )$, and
                        $\gamma$ does not, then $\beta$ travels a distance of at most $K$ in 
                        $(g, i ) ( P_i , i )$
                        \item If $\beta$ and $\gamma$ both penetrate a coset $(g,i) ( P_i , i )$,
                        then the $P_i$-distance between the points at where $\beta$ and $\gamma$
                        first enter $(g,i) ( P_i , i ) $ is less than $K$. The same is true for 
                        the points where $\beta$ and $\gamma$ last leave $(g,i)(P_i , i )$.
                \end{enumerate}

        It is easy enough to map groupoid elements to group elements. 
        We do this via a homomorphism of groupoids. We do this by
        taking a maximal tree in the base graph from which the groupoid is formed
        and contract it (ie all groupoid generators corresponding to edges in that 
        tree are mapped to the identity element.) 
        The image of this homomorphism is a groupoid with one object, which is a group.
        It is a well known result 
        that the group obtained is independent (up to isomorphism) of the 
        choice of maximal tree.  We choose the edges $e_{1,j}$. These do indeed 
        form a maximal tree because they form a tree that includes all vertices
        of the base graph.

        For our argument, we will also need a method to ``lift'' a path $\gamma$ in $\Gamma_G$ to 
        a path in $\Gamma_{\tilde{G}}$. We do this as follows: 
        Let $\tilde{\gamma}(0) = p_1 $.         Let $a\in X $ be 
        the generator such that $\gamma(t)a = \gamma(t+1)$. Let $\tilde{\gamma}(t')$ be the 
        lift of $\gamma_{[0,t] } $. Then if $a\not\in Y_i$ for any $i$, we set 
        $\tilde{\gamma}( t' + 1 )  = \tilde{\gamma}(t') \tilde{a} $.
        If $a\in P_i$ for some $i$, there are two cases: either $\tilde{\gamma}(t')$ 
        is in $(G , i )$ or it is in $(G,j)$ for some $j\neq i $. 
        \begin{itemize}
        \item
                If $\tilde{\gamma}(t')$ is in $(G,i)$ , then we set $\tilde{\gamma ( t' + 1 ) } = 
                \tilde{\gamma}(t) \tilde{a} $ 
        \item
                otherwise, we set 
                        \begin{itemize}
                                \item[] 
                                        $\tilde{\gamma} ( t' + 1 )  = \tilde{\gamma } ( t' ) e_{j,i}$ 
                                \item[] 
                                        $\tilde{\gamma }( t' + 2 )   = \tilde{\gamma } ( t' ) e_{j,i} \tilde{a} $ 
                        \end{itemize}
        \end{itemize}

        Note that since the $P_i$ are malnormal, it is not possible that 
        $a\in P_i , a\in P_j, i\neq j$. So $\tilde{\gamma}$ is well defined. $\tilde{\gamma}$ 
        is the preimage under a quasi-isometric embedding of $\gamma$. It follows that the bounded
        coset penetration property is true for $G$.

\end{proof}

\section{Electric Isoperimetric Inequalities}

\subsection{Preliminary Definitions}

To proceed, we need to consider the notion of ``electric isoperimetric functions''.
An isoperimetric function $f(n)$ for a group $G$ is a bound for the area of a loop in the 
Cayley graph $\Gamma_G$ of length at most $n$. An \emph{electric} isoperimetric function of $G$ with respect 
to a group $P<G$ is also a  bound for the area of a loop in $\Gamma_G$, but we assign area 
0 to all loops in $P$. We clarify this in the following definition: 

\begin{definition} 
Let $G$ be a finitely presented group and $\{P_i \}   $  be a collection of 
finitely presented  
subgroups of $G$ indexed by some set $I$.
Fix finite presentations $\langle X \vert R \rangle = G$, and $\langle B_i \rangle = P_i $ 
with $B_i \subset X$.
Let $S_i$ be sets of defining relators for the $P_i$. Let $S=\bigcup_iS_i$
Given $g\in F(X)$ such that  $g$   is in the normal closure of $R$.
 define the {\em electric area of $g$ with respect  to $\cup_i P_i$} to be
$$\area_{\cup_iP_i} (g) = 
\inf \{ n: g = \Pi_{i=1}^n ( p_i r_i p_i^{-1}  \Pi_{j=1}^{b_i} p_{i,j} s_{i,j} p_{i,j}^{-1} ) \}
$$
where $r_i\in R$ and $s_{i,j}\in S$.

Define the {\em electric length of $g$ with respect to $\cup_i P_i$ } to be 
$$\text{card} \{ j : x_j\not\in B_i \forall i \} $$
We say that $f: \N \rightarrow \N $ is an {\em electric isoperimetric function of $\langle X\vert R\rangle
$ with respect to $\{ \langle B_i \rangle \} $ }   if for any $g\in F(X)$ such that $\pi_G(g) = 1_G $ 
and the electric length of $g$ is less than $n\in\N$, 
$ \area_{\cup_i P_i}(g) < f(n)$

A \emph {$k$-local geodesic path/word} is a path/word that has the property that all sub-paths/sub-words of length no more than $k$ are geodesic.
We can also apply this definition to pseudo-metrics, hence we have 
$k$-local electric geodesics.
\end{definition}

\subsection{Results}

The main goal of this section is the following result: 

\begin{theorem}
Let $G$ be a geometrically finite hyperbolic group. Let $P_1 , \dots P_m $ 
be the conjugacy classes of parabolic subgroups. 
Then $G$ has a linear electric isoperimetric function with respect to $\bigcup_i P_i$
\end{theorem}

Now we use a similar path shortening algorithm to that used for the coset
graph.

Firstly, we recall a result about $\delta$-hyperbolic metric spaces: 
\begin{lemma}
        Let $M$ be a $\delta$-hyperbolic metric space. Then there exists a $k\in\N$ such that
        $k$-local geodesics $k$-fellow travel geodesics in $M$.
\end{lemma}
Note that this implies the same is true for $\delta$-hyperbolic pseudo-metric spaces
( given a pseudo-metric space, we apply the result to the quotient metric space 
obtained by identifying all points distance $0$ apart.  )

First, we prove 
\begin{lemma}\label{rel_iso_exists}
There is a function $f: \N \rightarrow \N $ 
such that for any cycle $\eta$ in $\Gamma_G$  
where $\bar{\eta} = 1 $ and
the length of $\hat{\eta}$ is less than $k$, $\eta$ has electric 
area of no more than $f(k)$. In other words, we show that an electric  isoperimetric function
for $G$ exists. 
\end{lemma}

This argument only provides an exponential bound. However, the proof 
of the proposition requires the existence of this function since it involves 
a ``coarse'' decomposition of an arbitrary word to words of bounded length.
We need the lemma to prove that these words of bounded length really 
do have area bounded by a uniform constant.

\begin{proof}[Proof of lemma]
We introduce the idea of {\em coset reduction}. 
This works as follows: if any parabolic segment $y_i$ of $w$ is non geodesic, 
replace it with a $\Gamma_G$-geodesic.
Strictly speaking, if $w=\alpha y \gamma$ where $y$ is a non-geodesic element of $P$ , then 
replace $w$ with $(\alpha y {y'}^{-1} \alpha^{-1}) \alpha y' \gamma$
where $y'$ is a $P$-geodesic.
Coset reduction preserves electric length and electric area. So for any word $w$ 
in $\Gamma_G$, there is a corresponding
word  $w'$ whose parabolic segments are all geodesics, which has the same electric area 
and electric length as $w$. See fig \ref{coset_reduction1}.
\begin{figure}
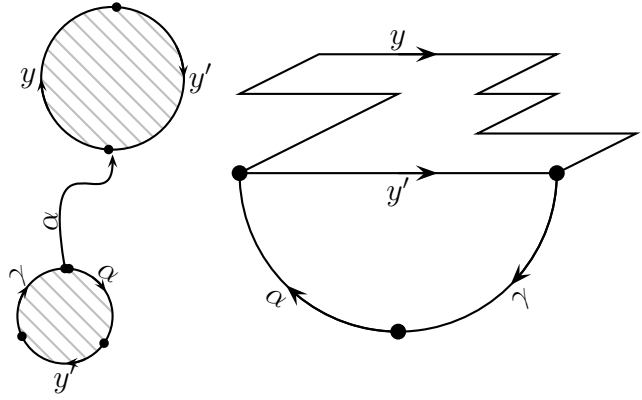
\label{coset_reduction1}

\begin{center}
\caption{Coset reduction -- A Dehn Diagram (left), and a schematic in $\Gamma_G$ (right) }

\input coset_reduction_1.tex
\end{center}
\end{figure}

First, we deal with the case where $\eta$ is without 
backtracking. We can assume that $\eta$ is coset reduced. 
Because an isoperimetric function for $G$ exists ( \cite{epstein1} ), 
it suffices to get a uniform bound on the $\Gamma_G$-length of $\eta$, which is
what we do in this case. Observe that $\eta$ is a $(0, \length (\hat{\eta}) )$
electric quasi-geodesic.
Compare $\eta$ to the trivial path in $\Gamma_G$.
Observe since that the trivial path only penetrates the costs 
$1\cdot P_1, \dots , 1 \cdot P_N $. 
cosets, 
any cosets penetrated by $\eta$ ( with the exception of $1\cdot P_1, \dots ,1\cdot P_n $ )
are not penetrated by the trivial path.
So by the bounded coset penetration property, $\eta$ travels a distance of 
at most $C$ in any of the cosets it penetrates except 
$1\cdot P_i,\dots 1\cdot P_n$, where 
$C(\length ( \hat{\eta} )\leq C(k) )$ is the constant 
given by the bounded coset penetration property.

We now deal with the cosets $1\cdot P_1 ,\dots 1\cdot P_n$. For each $P_i$, there
is some $t$ such that $\eta(t)^{-1}\not\in P_i$ (otherwise, the electric length of $\eta$
is zero and there is nothing to prove ). Consider the cyclic permutation
$\eta' = x_t x_{t+1} x_{t+2} \dots x_{k'} x_0 x_1 \dots x_{t-1} $ of the word $\eta(t)$
where $k'$ is the $\Gamma_G$-length of $\eta$.
Note that $\eta'$ penetrates the coset 
\begin{gather*}
                x_t \dots x_{k'} \cdot P_i  = \\
                \bar{\eta'}^{-1} x_t\dots x_{k'}\cdot P_i \\
                x_{t-1}^{-1} \dots x_0^{-1} x_{k'}^{-1}\dots x_t^{-1}  x_t \dots x_{k'} \cdot P_i =\\
                x_{t-1}^{-1} \dots x_0^{-1} \cdot P_i    = \\
                \eta'(t)^{-1} P_i
\end{gather*}
                and $1$ does not penetrate this coset. So $\eta'$ travels a distance of at most 
                $C$ in $\eta'(t)^{-1} P_i$. Translating $\eta'$ and $\eta'(t)^{-1} P_i$ 
                by $\eta'(t)$, we conclude that $\eta$ travels a distance of at most 
                $C$ in $P_i$.

Since the parabolic segments of $\eta$ are bounded by $C$,
 the length of $\eta$ is at most
$k(C+1)$. 

Now suppose that $\eta$ backtracks.
We will use induction.  The idea behind our induction is that if
$\eta$ backtracks, we
can either find a word $\eta'$ shorter than $\eta$ such that  
$\area(\eta) $ is bounded by a linear function of the electric length of $\eta'$ and 
$\area(\eta')$
If $\eta'$ turns out to be without backtracking, we are done by the previous argument.
If $\eta'$ backtracks, we reiterate the reduction process.
This will show that an electric isoperimetric function exists (though the exponential 
bound obtained is horribly un-optimal) 
Suppose that for all $k'<k$ ,
there is a constant $f(k')$ such that any loop of electric length less than or equal to $k'$ 
has electric combinatorial area of at most $f(k')$. Then  suppose
$\eta$ is of electric length $k$.

Decompose $\eta$ as $\eta_1 \eta_2 \eta_3$ where $\eta_2$ is a backtracking sub-word
( ie $\eta_2$ is a non-parabolic backtracking sub-word of $\eta$ ).
Let $\ell_1, \ell_2$ and $\ell_3$ be respectively the electric lengths 
of $\eta_1, \eta_2$ and $\eta_3$.
Let $y$ be a parabolic geodesic such that $\bar{y} = \bar{\eta_2} $.
Then 
$$\eta\sim  ( \eta_1\eta_2 y^{-1}\eta_1^{-1} ) ( \eta_1 y  \eta_3  ) $$
So 
$$f(k)\leq \area(\eta_1\eta_2 y^{-1}\eta_1^{-1}) + \area(\eta_1 y  \eta_3 ) \leq
f(\ell_2 - \ell_1 ) + f ( \ell_1 + \ell_3 - \ell_2 )\leq 2f(k-1 )  $$
\begin{figure}
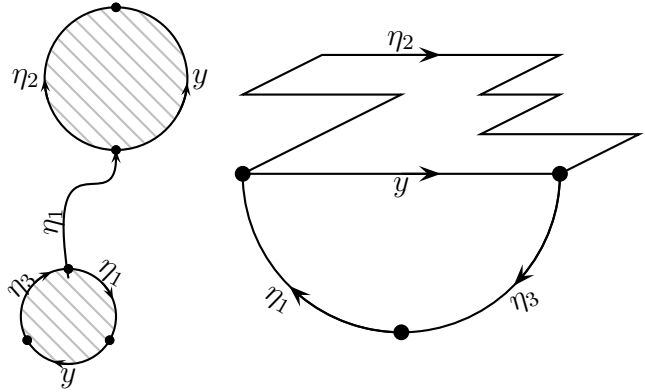

\caption{$\eta$ backtracks -- Dehn diagram (left), and schematic in $\Gamma_G$ (right)}
\begin{center}
\input coset_reduction_2.tex
\end{center}

\end{figure}

Hence we get $f ( k ) \leq 2f(k-1 ) $.
In the case where $\eta$ is without backtracking, since $G$ is biautomatic
( see \cite{epstein1} ), it has quadratic 
isoperimetric inequality and we get that $f(k) < Ek^2 + E'k + E''$ for 
some constants $E, E', E''$ that depend only on the presentation of $G$.
Hence $G$ has an electric isoperimetric function. 
\end{proof}

\begin{proof}[proof of proposition]
Let $\eta$ be a path in $\Gamma_G$, and let $\text{hat} : \Gamma_G \rightarrow
\hat{\Gamma}_G$ be the natural (identity) map. That is to say that
the spaces are the same, the metric is different.
  Then given a closed path $w\in \Gamma_G$ , there is a
corresponding path in $\hat{\Gamma}_G$ that is obtained by the formula $ \hat{w} (t)
= \text{hat}(w(t)) $. 
The strategy we will employ is to use two types of moves to reduce $w$ to a cycle such that
$\hat{w} $ is a $k$-local geodesic. The first type of move is {\em coset reduction}
as described previously. 

If $w$ is coset reduced, then we can perform an \emph{electric length reduction} as follows: 
if $\hat{w}$ is not a $k$-local geodesic, then we can 
decompose $\hat{w} = \alpha\beta\gamma = \alpha\beta\xi^{-1}\alpha^{-1}\alpha\xi\gamma$ where $\beta$ is a minimal non-geodesic segment of $\hat{w}$, and choose $\tilde{\alpha} , \tilde{\beta} , \tilde{\gamma}$ 
such that $\text{hat}(\tilde{\alpha}) = \alpha , \text{hat} ( \tilde{\beta} ) = \beta $ and
$\text{hat} ( \tilde { \gamma } ) = \gamma $ and further $w=\tilde{\alpha}\tilde{\beta}\tilde{\gamma}$. 
Choose $\tilde{\xi}$ such that each of the parabolic segments of $\tilde{\xi}$ is a 
geodesic.The length of $\beta\tilde\xi^{-1}$ is less than $2k$. So the relative combinatorial area
of $\beta\xi^{-1}$ is at most $A=f(2k)$ where $f$ is a relative isoperimetric function for $G$.
Let 
$v = \tilde{\alpha} \tilde{\beta } \tilde{\xi^{-1} } \tilde{\alpha^{-1}}$ and 
let $w' = \tilde{\alpha} \tilde{\xi} \tilde{ \gamma}$. Then $w=vw'$ , so $\area_P(w) \leq \area_P ( v) + \area_P (w' ) = A + \area_P (w' )$.
Moreover , $\length (\hat{w' }) < \length ( \hat{w} ) $. 

We iterate the procedure above inductively: given a word $w$, first we can
coset reduce it , to find a word $w'$ such that $w$ is freely equal to $w'$ and
$\area_P(w') \leq \area_P(w) $. If $\hat{w'}$ is not a $k$-local geodesic, then we
can replace it by a word $w'' $ where $\length(\hat{w''}) < \length ( \hat{w'})$
where $\area_P ( w ) \leq \area_P (w') \leq A +\area_P ( w'' )$. This process must
terminate. So inducting, we see that $\area_P(w) \leq A (\length(w) ) + \area_P ( \theta ) $ where
$\hat{\theta} $ is a  $k$-local geodesic and $\theta$ is coset reduced. Note that since $\hat{\theta}$ 
is a $k$-local geodesic, it is a quasi-geodesic, hence $\length(\hat{\theta}) \leq D$ 
for some constant $D$ that depends only on $G$.
By lemma \ref{rel_iso_exists}, this has bounded area. 

\end{proof}

\section{Proof of Theorem 1} 

We prove the following:
\begin{theorem*}[\ref{maintheorem1}]
Let $A\rightarrow_{\iota} E \rightarrow_{\pi} G  $ 
be a boundary-trivial central extension of a 
geometrically finite hyperbolic group $G$. 
Then $E$ is biautomatic.
\end{theorem*}
The proof is similar to that in \cite{neumann1}.
We will first prove the case where $A=\Z$.
Given this, the more general case follows, see \cite{neumann1} for
example.

First, we need some definitions:
\begin{definition}
        A $2$-cocycle $\sigma$ is \emph{weakly bounded} if the sets 
$\sigma(g,G)$ and $\sigma(G,g)$ are finite for
each $g\in G$. $\sigma$ is also said to be \emph{$L$-regular} 
if for each $h\in G$ and $a\in A$, the 
set $\{ g\in G \vert \sigma(g,h) = a \} $ is an $L$-rational subset of $G$.
A cohomology class is $L$-regular if one of its representative cocycles is.
\end{definition}

The theorem follows from the following result:
\begin{lemma}\label{mainlemma} 
Let $J$ be a geodesic biautomatic structure on a geometrically
finite hyperbolic group $G$
such that $J$-words do not backtrack. 
Let $A\rightarrow E \rightarrow G$ be a boundary-trivial central
extension of $G$.
Then there is a
rational structure $L$ on $G$ equivalent to $J$ such that the cocycle $\sigma$
that defines this central extension is $L$-regular.
\end{lemma}

\ref{mainlemma} can be used to show that $E$ is biautomatic.
The argument is as follows:
It is a well known result 
that given two rational structures $L_1$ and $L_2$, a set is $L_1$-rational
if and only if it is $L_2$-rational. 
Since $\sigma$ is $L-regular$, this implies that 
$\sigma$ is also $J$-regular.  
Theorem A of \cite{neumann1}  states that 
$E$ carries a biautomatic structure if and only if 
$G$ has a biautomatic structure $J$ for which $\sigma$
is given by an $J$-regular cocycle.  
Since $J$ satisfies the conditions given by theorem A
from \cite{neumann1}, $E$ is biautomatic.
So the goal for the remainder of this section is to prove \ref{mainlemma}.

Let $X$ be a finite collection of symbols  in 1-1 correspondence with a
generating set for $E$ closed under inversion.
Let $Y_{i,j}\subset X$ be a collection of parabolic generators for conjugacy 
class representatives
$P_i$ of the parabolic subgroups of $G$.
Let $\pi : E\rightarrow G $ be the projection map. 
Then $X$ maps onto a generating
set for $G$ via $\pi$. So we will use $\bar{x} $ to denote the generators of
$E$ and $\pi(\bar{x})$ to denote the generators of $G$.
Let $\bigcup S_{i,j}$ be a collection of defining relators for 
$P_i=\langle Y_i\vert S_{i,j} j\in \{1,\dots , n_i\} \rangle$
and let $R$ be a set of defining relators for $G=\langle X \vert R \rangle $
such that $\bigcup_{i,j}S_{i,j} \subset R$.
For each $g\in G $, we define an ordering $<$ on $\pi^{-1}(g)$ by the 
natural ordering on $\Z$ and we refer to a maximal element of a subset 
of $U\subset \pi^{-1}(g)$ as the maximum, denoted by $\max(U)$.
We say that the set $U$ is 
{\em bounded above} if it has a maximal element
and that it is 
{\em bounded } if it has a minimal and maximal element.
Let $\Gamma_G$ be the Cayley graph of $G$.
For the purpose of this section, for a word $w=x_1 x_2 \dots x_j \in X^*$ ,
we denote  by $\length(\hat{w})$ the 
number of $x_i$ satisfying $x_i \not\in Y_{i,j}$, ie $\length ( \hat{w} ) $ is the number of non parabolic generators 
appearing in the word $w$, and take the view of the coned off Cayley graph $\hat{\Gamma}_G$
as a pseudo-metric space where the path between any any point in $gP_i$ and it's associated 
cone point $v_i(g) $ is assigned a length of $0$.

First, we define our language $L$ which will be shown to satisfy the hypothesis of 
\ref{mainlemma}. $L$ will be the language of {\em maximizing words}. The following result
verifies the appropriateness of this definition.

\begin{lemma} 
There is a $C>0$ such that for any $g\in G$, 
such that a {\em maximizing word } 
$$\max\{ \bar{w} \iota(-C \length(\hat{w})) \vert w\in X^* , \pi(\bar{w}) = g \} $$
exists. Moreover, there exists a $\lambda$ that does not depend on $w$
such that if $w$ is a maximizing word, then $w$ is a 
$(\lambda, 0 ) $ quasi-geodesic in
$\hat{\Gamma}_G$
\end{lemma}

\begin{proof}
Let $T\in\Z $ be defined by the formula 
$\iota(T)= \max \{ \bar{r} \vert r^{\pm 1 }\in R - \bigcup_i S_{i,j} \}$. 
Let $K$ be the constant given by the isoperimetric inequality for $G$
relative to the $P_i$. Then if $\bar{w} = 1_G $ , 
then there exist a finite collection of $p_j\in F(X)$ and relators
$r_j$ such that 
\begin{enumerate}
\item $w =_{F(X)} \Pi_j p_j r_j p_j^{-1} $ 
\item there are at most $K\length(\hat{w}) $ 
$r_j $ with $r_j\not\in P_i$ for some $i$. 
\end{enumerate}
Let $K$ be the minimal 
number.
Note that $ \iota (g p r_j p^{-1}) = \iota (g ) $ if $r_j \in S_{i,j} $ for some $i$.
Let $A$ be the {\em relative } combinatorial area enclosed by $w$.
So $\iota(\bar{w}) \leq \iota ( T A ) \leq \iota ( T K \length(\hat{w}) ) $ Let $C\in \Z$, 
$C>TK$. For any $g\in G , v,w \in X^* $ with $\pi(\bar{w}) = \pi(\bar{v})  = g$,
we have 
\begin{gather}
\bar{v}^{-1}\bar{w} \leq \iota ( TK \length(\hat{v}^{-1} \hat{w} ) ) 
           = \iota  ( TK \length ( \hat{v} ) ) \iota ( TK \length ( \hat{w} ))  \label{eqn1}
\end{gather}

This implies that 
\begin{gather*}
\bar{w} \iota ( -C \length ( \hat{w} ) ) \leq \bar{w} \iota ( -TK \length( \hat{w} ) ) \\
\leq \bar{v}^{-1} \iota ( TK \length\hat{v} )  
\end{gather*}

So the set $\{ \bar{w} \iota ( -C \length( \hat {w} ) ) : w\in X^* , \pi(\bar{w})=g \} $ is 
bounded , in particular it has a maximum.

Now if $w$ is maximizing and $v$ defines a geodesic in $\hat{\Gamma}_G$,
then $\bar{v} \iota ( -C \length( \hat{v} ) ) \leq \bar{w} \iota ( -C \length( \hat{w} ) ) $.
so 
\begin{gather}
\iota ( C \length( \hat{w} ) - C \length( \hat{v} ) )  \leq ( \bar{v}^{-1} \bar{w} ) . \label{eqn2}
\end{gather}

It follows from equations 
\ref{eqn1} and \ref{eqn2} that $C( \length( \hat{w} )  - \length( \hat{v} ) ) \leq 
TK( \length( \hat{v} )  + \length( \hat{w} ) ) $.
Rearranging terms, we get
$$\length(\hat{w}  ) \leq  \frac{ \length(\hat{v} ) ( C + TK  ) }{  C-TK } $$
So let $\lambda = (C+TK ) / ( C-TK ) $. Then 
$\length( \hat{w}(t) ) \leq \lambda \length ( \hat{u}  ) $ for each $t$ , where $u_t$ is 
a geodesic such that $\overline{u_t}=\overline{w(t)}$.
So
$\hat{w} $ is a $(\lambda, 0 ) $
quasi-geodesic. 
\end{proof}

Hence we are in a position to define a language $L\subset X^*$ to be 
the set of all maximizing words 
such that any parabolic sub-word $Y_{i,j_1} Y_{i,j_2}\dots Y_{i,j_n}$ follows the 
biautomatic structure on $P_i$. Note that the parabolic sub-words of a word $w$
have no bearing on whether or not $w$ is maximizing, so changing a word 
$w$ by replacing some of its parabolic sub-words with different parabolic 
sub-words will preserve maximality (or non-maximality) of $w$.

\begin{lemma} 
The language $L$ is regular and has the asynchronous fellow traveler property.
The words in $L$ do not backtrack. Moreover, any sub-word of an $L$-word is in $L$
(that is, sub-words of maximizing words are maximizing)
\end{lemma}
\begin{proof}
First, we show that sub-words of maximizing words are maximizing.
Let $u\in X^* -  L$. Let $u_1$ be the shortest initial 
section that is not in $L$. Let $w$ be the trailing section (ie such that
$u_1 w = u$)
By minimality of $u_1$, it follows that $u_1=u_0x$ with $u_0\in L$.  Let
$v_1\in L$ with $\pi(\overline{v_1}) = \pi (\overline{u_1})$.  Then since $L$
has the (asynchronous) fellow traveler property, $u_1$ and $v_1$ fellow
travel.  Let $v$ be the word obtained from $u$ by replacing the initial segment
$u_1$ with $v_1$ ($v=v_1 w $). Then 
$$\overline{u_1}\iota( -C \length ( \hat{u_1} ) ) < \overline{v_1} \iota( -C \length ( \hat{v_1} ) ) $$
Hence , applying (1)
$$\overline{u_1 w}\iota( -C \length ( \hat{u_1} ) ) \iota( -C \length ( \hat{w} ) )< \overline{v_1 w} \iota( -C \length ( \hat{v_1} ) )\iota( -C \length ( \hat{w} ) ) $$
So 
$$\bar{u}\iota( -C \length ( \hat{u} ) ) < \bar{v} \iota( -C \length ( \hat{v} ) ) $$
This argument easily generalizes to arbitrary ( not just initial ) sub-words 
of maximizing words.

To show that $L$ is a regular language, it suffices to demonstrate that it
satisfies the falsification by fellow traveler property. 
We then appeal to the results \ref{fft1} and \ref{fft2}.

Recall that it is sufficient to prove that $L$ has a 
height function with the following properties:
\begin{itemize}
        \item Weak translation invariance
        \item Bounded difference: 
        $\mathcal{H}$ has the {\em bounded difference} property if there exists $K\in \N$ such that 
        for all $w\in X^*, x\in X$, $$\vert \mathcal{H}(w) - \mathcal{H}(wx) \vert < K$$
        \item The trivial word is maximizing.
\end{itemize}

Let 
$\mathcal{H}(w) = \iota^{-1} ( wv^{-1} ) - C \length( \hat{w} ) + C \length( \hat{v} )  $ 
where $v$ is a maximizing word with
$\pi(\overline{wv^{-1} } ) = 1 $. 
The height function measures the difference between
a word and an appropriate  maximizing word.

First, we prove weak boundedness: 
Let $u,v, u', v' \in X^* $ with 
$\pi(\bar{u} ) = \pi( \bar{u'} ) , \pi( \bar{v} ) = \pi( \bar{v'} ) $
Let $a$ be a maximizing word with 
$\pi( \bar{a} ) = \pi ( \overline{uv} ) = \pi( \overline{u'v'} ) $.
Write 
        $ a= bc$ 
        where 
        $
        \pi( \overline{b} ) = 
        \pi ( \overline{u} ) = 
        \pi( \overline{v} )
        $.
$b$ and $c$ are sub-words of a maximizing word, hence they are maximizing.
Then 
\begin{gather*}
        \iota [ \mathcal{H} ( uv ) - \mathcal{H} ( u'v' ) ]=\\
        ( uva^{-1} ) - \iota [ C \length( \hat{uv} ) ] + \iota [ C \length( \hat{a} )   ]  - 
        \left[ ( u'v'a^{-1} ) - \iota [ C \length( \hat{u'v'} ) ] + 
                \iota [ C \length( \hat{a} )   ]  \right] = \\
        ( uvc^{-1}b^{-1} ) - \iota [ C \length( \hat{uv} ) ]   - 
        ( u'v'c^{-1}b^{-1} ) + \iota [ C \length( \hat{u'v'} ) ]   = \\
        ( vc^{-1}) ( u b^{-1} ) - \iota [ C \length( \hat{u} ) + \length ( \hat{v} ) ]   - 
        ( v'c^{-1} ) ( u'b^{-1} ) + \iota [ C \length( \hat{u'} ) + \length (\hat{v'} ) ]   = \\
        \iota [ \mathcal{H} ( u ) - \mathcal{H} ( u' ) ]
        + \iota [ \mathcal{H} ( v ) - \mathcal{H} ( v' ) ]\\
\end{gather*}

Next, we prove that the words in $L$ do not backtrack. 
This implies immediately that the trivial word is maximizing.
Suppose $v\in X^*$ backtracks.
Then $v=v_1v_2v_3$ where $\pi(\overline{v_1} ) P_i = \pi ( \overline{v_1v_2} ) P_i$ 
for some $i$ and $v_2$ is some nonempty word. Let $u$ be a word in $Y_{i,j} $
such that $\pi(\overline{u}) = \pi(\overline{v_2})$.
Then  

Recall that equation \ref{eqn1} asserts that 
$$\overline{v_2 u^{-1} }\leq \iota ( TK \length ( \widehat{v_2 u^{-1} } ) ) $$ 
where $T$ and $K$ are as used in the definition of the maximizing word, and since $C>TK$, 
$$ \iota ( TK \length ( \widehat{v_2 u^{-1} } ) ) < \iota ( C \length ( \widehat { v_2 u^{-1} } ) ) $$
so
$$\overline{v_2u^{-1}} < \iota ( C \length \widehat{v_2u^{-1}} )   $$

Since $\length \hat{u} = 0 $, 
$$\overline{v_2u^{-1}} < \iota ( C \length \widehat{v_2} )   $$  which implies that
$$\iota (- C \length \widehat{v_2} )\overline{v_2} <  \overline{u} =\bar{u}-\iota ( - C \length ( \hat{u} ) ) $$
So $v$ is not a maximizing word ( since $v_1 u v_3 $ is closer to maximizing ). 
Note that this also tells us that all words in the $Y_{i,j}$ are maximizing.

Now we show the boundedness property.
Let $u\in X^*$, and let $v\in X^*$ be a maximizing word
with $\pi(\bar{v} ) = \pi(\bar{u} ) $.
Let $x\in X$, and $w\in X^*$ be a maximizing word
with $\pi(\bar{w}) = \pi( \overline{ ux } )$.
By weak translation invariance,
$$
        \mathcal{H} ( u ) - \mathcal{H}( ux ) = 
        \mathcal {H}( v ) - \mathcal{H}( vx ) =
        - \mathcal { H } ( vx ) 
$$
So we need to evaluate $\mathcal{H}(vx)$. 
Since $v$ is maximizing, $\mathcal{H}( wx^{-1} ) \leq \mathcal{H} ( v ) $.
By weak translation invariance, $-2C = \mathcal{H}( wx^{-1}x ) \leq \mathcal{H}(vx) $.
So $$\mathcal{H}(w) - \mathcal{H}(wx) = - \mathcal{H}(vx) \leq 2C$$

Now we prove the asynchronous fellow traveler property. 
Let $w_1 , w_2 \in L $ where $\pi (\overline{w_1 x} ) = \pi(\overline{w_2} ) $.
Then $\hat{w_1}$ and $\hat{w_2}$ are both quasi geodesics without backtracking. 
The bounded coset penetration property , 
combined with the fact that $w_1$ and $w_2$ follow the same biautomatic structure
inside the cosets immediately implies that 
$w_1$ and $w_2$ satisfy the asynchronous fellow traveler property.
\end{proof}

We now define a cocycle $\sigma$ with a view to proving that $\sigma$ is $L$-regular.
We define $\sigma$ via a section $\rho$.
Define $\rho : G \rightarrow E $ by setting 
$\rho ( \bar{w} ) =\bar{w} \iota ( -C \length(\hat{w} ) ) $

First, we prove weak boundedness.
Let $h_1, h_2 \in \pi^{-1} ( g)$ , $g\in G $ , and 
define $h_1 - h_2 \in Z $  by the equation $\iota ( h_1 - h_2 ) = h_1 h_2^{-1} $.
The cocycle $\sigma$ that comes from the section $\rho$ is defined by the formula
$\sigma (g_1 , g_2 ) = \rho ( g_1 ) \rho ( g_2 ) - \rho(g_1 g_2 ) $.
So $\sigma$ is weakly bounded if and only if 
$\rho(\pi(\bar{x})) \rho(g) - \rho(\pi(\bar{x})g): G\times X\rightarrow \iota(A) $  and
$\rho(g) \rho (\pi(\bar{x})) - \rho (gx):G\times X\rightarrow \iota(A) $ are bounded functions.
The following sub-lemma proves that $\sigma$ is weakly bounded:
\begin{lemma}\label{resulttech} 
$\sigma$ is weakly bounded with respect to the generating set
$\pi(\bar{X})$. Explicitly, for all $g\in G$ and $x\in X$, \\
\begin{enumerate}
\item $ \vert \rho(g) \bar{x} - \rho(g\pi(\bar{x})) \vert \leq C $
\item $ \vert \bar{x} \rho(g) - \rho ( \pi(\bar{x})g ) \vert \leq C $
\end{enumerate}
\end{lemma}

\begin{proof}[Proof of \ref{resulttech} ]
First , we prove (1). Let $w_1 , w_2$ be maximizing words with $\pi(\overline{w_1})=g$ and
$\pi(\overline{w_2}) = gx$. Then since $w_2$ is maximizing and $\pi(\overline{w_1 x} ) = \pi(\overline{w_2})$, there exists
$h\geq 0 $ such that 
\begin{gather}
\rho(gx) = \overline{w_2} \iota ( -C \length (\hat{w_2} ) ) = \overline{w_1} \iota ( -C \length(\hat{w_1} ) ) \bar{x} \iota ( -C )\iota (h) = \rho(g)\bar{x} \iota(h-C)   \label{eqn3}
\end{gather}
By symmetry, since $w_1$ is maximizing and $\pi(\overline{w_1}) = \pi ( \overline{w_2 x^{-1} } ) $, 
there exists $k\geq 0 $ such that 
\begin{gather}
\rho(g) = \overline{w_1} \iota ( -C \length (\hat{w_1} ) ) = \
\overline{w_2} \iota ( -C \length(\hat{w_2} ) ) \overline{x}^{-1} \iota ( -C )\iota (k) \
= \rho(gx) \bar{x}^{-1}\iota(k-C)   \label{eqn4}
\end{gather}
Combining equations \ref{eqn3} and \ref{eqn4}, we get 
\begin{gather}
\rho(g)^{-1}\rho(gx) = \bar{x} \iota ( h-C ) =  \iota (  C - k ) \bar{x} 
\end{gather}
Since the image of $\iota$ is central in $E$, we obtain 
$h-C = C -k $. Therefore, 
$$-C\leq h-C = C-k \leq C $$
which completes the argument. 
The proof of part 2 is analogous.
\end{proof}
We now prove that $\sigma$ is $L$-regular where $L$ is a language of quasi-geodesics
on $G$. Equivalence classes of rational structures on $G$ depend only on choices 
of rational structures on the parabolic subgroups of $G$. Fix a bijective
biautomatic structure $J$  on $G$. 

Now we lift $L$ via $\rho$ . 
The generic way to do this is to use an alphabet consisting of a set in correspondence
with $\{ \rho(\pi(\bar{x})) \iota ( -\sigma ( g,\pi(\bar{x}) )) : g\in G , x\in X \}$ 
then we can lift a word $v\in X^*$ to a word in $E$ whose initial vertices have values
$\rho ( \pi ( \overline{x_1} ) ), \rho ( \pi ( \overline{x_1 x_2} ) ) ,\rho ( \pi ( \overline{x_1x_2 x_3} ) ) ,\dots $
by using 
        $       \rho ( \pi(\overline{x_i}) ) 
                \iota ( - \sigma ( \pi(\overline{x_1 x_2 \dots x_{i-1}} , 
                \pi(\overline{x_i}) ) )) $
        as our 
$i_{th}$ generator. In fact this is how we will lift the biautomatic structure $J$ 
to get a biautomatic structure on the group $E$.
However, lifting the  rational structure $L$ is somewhat simpler -- 
because $L$ consists of maximizing 
words, 
$\sigma ( \pi(\overline{x_1 \dots x_{i-1} }), \pi(\overline{x_i}) ) = 0 $. 
So the set $\{ \rho(\pi(\bar{x})) : x\in X \} $ is good enough.
For any $x\in X-Y $, define $x' = \rho(\pi(\bar{x})) = \bar{x} \iota( -C) $ and 
for any $y_{i,j} \in Y_i$ , we use $y_{i,j}' = \bar{y_{i,j}}$.
Let $L'$ be the language with the evaluation map $x\mapsto x' , y_{i,j} \mapsto y_{i,j}'$.
Let $v=x_1 x_2 \dots x_n  \in L'$. 
Then the initial segments of $v$ have values
$\rho ( \pi ( \overline{x_1} ) ), \rho ( \pi ( \overline{x_1 x_2} ) ) ,\rho ( \pi ( \overline{x_1x_2 x_3} ) ) ,\dots $

We will complete the proof that $\sigma$ is $L$-regular with the following result;

\begin{lemma}
Let $A\rightarrow_\iota E \rightarrow_\pi G$ be a central extension. 
Let $L'$ be a regular language with an evaluation map $\phi_{L'} $
mapping bijectively onto the image of 
a section $\rho:G\rightarrow E$ , and suppose
$L'$ satisfies the asynchronous fellow traveler
property. Let $L$ be the normal form on $G$ defined by the regular language $L'$ 
with the evaluation map $\phi_L = \pi\circ\phi_{L'}$
Then the cocycle $\sigma$ determined by $\rho$
has the property that $\{ g\in G  : \sigma (g,\pi(\bar{x}))  =a \} $ is $L$-rational
for all $x\in X, a\in A$.
\end{lemma}
\begin{proof}
The asynchronous fellow traveler property for $L'$ implies that the language 
$\{ (u,v) \in L'\times L' : \bar{u} \bar{x} = \bar{v} \}$ is a regular language
accepted by an asynchronous two tape automaton. 
Note that this implies that 
$\{ (u,v) \in L'\times L' : \bar{u} \bar{x}\iota(-a+b)  = \bar{v} \}$ is also regular language
 since we can choose an automatic structure $L_A$ for $A$ and observe that $L'L_A$
has the asynchronous fellow traveler property.  
Projection onto the first factor 
in the following sense: $\{ u\in L' : \exists v\in L' , \overline{ux} =\bar{v}\iota(a-b)  \}$ 
is also regular  \cite{epstein1}  . The image of this in $G$ is 
$\{g\in G :\rho(g) \bar{x}   = \rho ( g \pi(\bar{x}) ) \iota(a-b)  \}$

If we choose ( without loss of generality ) $b$ so that $\iota(b) = \bar{x}^{-1} \rho( \pi(\bar{x} ) ) $
, then this is 
$\{g\in G :\rho(g) \rho ( \pi(\bar{x} ))  = \rho ( g \pi(\bar{x}) ) \iota(a)  \}$
which is the same as $\{ g\in G : \sigma (g,x ) =a \} $. 
\end{proof}

\begin{lemma}
The rational structure $L$ is equivalent to $J$
\end{lemma}
\begin{proof}

        Let $u$ be a $J$-word, and let $v$ be an $L$-word such that $\bar{u} = \bar{v}$.
        We need to show that $u$ and $v$ asynchronously fellow-travel. The first step
        towards doing this is to use a quasi-isometric map 
        $\phi : \Gamma_{\tilde{G}}\rightarrow \Psi$ ( where $\tilde{G}$ denotes the groupoid 
        constructed in the proof of \ref{BCP} ).
        There is also a quasi-isometric map $\psi: X^* \rightarrow \Gamma_{\tilde{G}}$ 
        as outlined in the proof of  \ref{BCP}.
        Composing this with the quasi-isometric map of $\Gamma_{\tilde{G}}$ into 
        $\Psi$ that maps parabolic subgroups to horospheres of $\hat{X}$, we obtain 
        two quasi-geodesics in $\Psi$.

        $\phi\circ\psi (u) $ and $\phi\circ\psi (v) $  asynchronously fellow travel. 
        This follows from a theorem of \cite{neumann2}: 

        \begin{lemma}[Lemma 3.2, \cite{neumann2} ]
                For any $\lambda>1, \epsilon>0 $, there exists $\ell\in\R$ 
                such that for any $(\lambda,\epsilon )$ quasi-geodesic in $\Psi$, $\gamma$
                the following holds:
                Let $\beta$ be a geodesic with the same endpoints as $\gamma$.
                Then $\beta$ asynchronously fellow travels $\gamma$ after 
                possibly modifying $\gamma$ on some of the horospheres that it 
                penetrates.
        \end{lemma}

        Since $\phi\circ\psi$ is a quasi-isometric map, it follows that $u$ and $v$ 
        also fellow travel.
\end{proof}

This proves \ref{mainlemma}. This completes the proof that $E$ is biautomatic.

 \part{Relatively Hyperbolic Groups}
 %
%

\section{The Cusped Off Cayley Graph is Hyperbolic}

\setcounter{footnote}{0}
\newcommand{\mass}{\mathcal{M}}
\subsection{The Groupoid Associated With a Relatively Hyperbolic Group}

In \cite{epstein1} , a geometric construction is used to analyze geometrically
finite hyperbolic groups. First, we consider a geometrically finite hyperbolic
group $G$ and a set of parabolic fixed points $p_i$ that lie in distinct 
$G$-orbits. Hence there are finitely many $p_i$. Centered at each $p_i$ , we
chose a horosphere $H_i$. It is possible to choose the $H_i$ small enough that
the $H_i$ orbits are disjoint.

A graph is then constructed as follows: firstly, for each horosphere
in $H_i$, we choose a vertex $x_i$. There is a unique geodesic ray 
$\gamma_i$ corresponding to $H_i$ emanating from $x_i$. 
We join each distinct pair of $x_i$ by a geodesic edge.

Let $a_j$ be a set of generators for $G$. For each $i,j$, we add 
the vertex $a_j \cdot x_i$ , and an edge joining $x_i$ to $a_j \cdot x_j$.
We extend this as follows: for each $g\in G$, $a_j$, $x_i$, 
add an edge from $g\cdot x_i$ to $ga_j \cdot x_i$.
The resulting graph is the Cayley graph of 
a groupoid whose associated group is $G$. 

We mimic this construction in the more general setting of groups
that are relatively hyperbolic in the sense of Farb.

        Let $G$ be a group that is hyperbolic relative to a 
        finitely presented subgroup $H$.
        Let $G= < X\vert R > , H = < Y\vert S >$, and assume that 
        $Y\subset X$ and $S \subset R$.

        The groupoid is constructed via a complex $K$. 
        We construct $K$ as follows: firstly, choose a base
        vertex $v=v_0$. For each generator $x$ of $G$, add a 
        directed edge $e_x$ beginning 
        and ending at the base vertex, and label the edge with $x$. 
        For each defining relator $r$ of $S$, add an oriented 2-disc 
        $D_r$ such that the boundary of $D_r$ is labelled by $r$. 
        Thus far, we have constructed the associated 2-complex to the
        group $G$. 

        Now attach an edge $e=e_0$ from the base vertex $v = v_0$ to $v_1$. 
        This process is iterated -- we attach an edge $e_i$ from 
        $v_i $ to $v_{i+1}$ for each $i\in\N$. 
        For each generator $y$ of the subgroup $H$, attach an edge
        beginning and ending at $v_i$,  $e_{i,y}$. For each defining 
        relator $s$ of $H$, attach a 2-cell $D_{i,s}$ whose boundary is 
        labelled by the relator $s$.
        Also attach a 2-cell $C_{i,y}$ to the cycle $e_{i,y} e_i e_{i+1,y}^{-1} e_i^{-1}$.
        This relator is some sort of commutativity relator, it essentially
        says that the edge $e_{i,y}$ should be homotopic to $e_{i+1,y}$.
        
        This 2-complex defines a groupoid $\tilde{G}$ ( namely, the 
        edge-path homotopy
        groupoid ). While the groupoid is not finitely generated, it does 
        have a local finiteness property -- there are only finitely many 
        edges at each vertex. Moreover, there are only finitely many homotopy
         classes of generators. Indeed, the fundamental group of 
        this complex is $G$.

        We define some notation for the generators of $\tilde{G}$.   
        For each edge $e_x$, where $x\in X$,    denote by $\tilde{x}$ 
        the corresponding generator of $\tilde{G}$.
        For each edge $e_{j,y_i}$, where 
        $y_i\in Y, j\in\N$, let $y_{i,j}$ be the corresponding generator of 
        $\tilde{G}$. The union of $Y_{i,j}$ will be denoted by $\tilde{Y}$.
        For each edge $e_i, i\in\N$, let $z_i$ be the corresponding generator in
        $\tilde{G}$. Define $\tilde{Z}$ to be the union of the $z_i$.
        It is sometimes convenient to have a notation for a word of the 
        form $z_i, z_{i+1}, \dots , z_j$ where $j>i$ or 
        $z_i, z_{i-1}, \dots , z_j$ where $i<j$. 
        We denote both of these words by $z_{i\dots j} $
        For each relator $r_i\in R$, there is a relator 
        for $\tilde{G}$
        $\tilde{r_i}$ given by the disk $D_{r_i}$. The set of such relators 
        will be denoted by $\tilde{R}$. Similarly, for each relator $s_i$ of 
        $G$, for each disc $D_{s_i,j}$, there is a relator of $\tilde{G}$ 
        denoted by $\tilde{s}_{i,j}$. The set of all 
        $\tilde{s}_{i,j}$ will be denoted by $\tilde{S}$.
        There is also a relator $c_{i, y_j}$
        for each disc $C_{i, y_j}$ where $i\in\N, y_j \in Y$. Denote the 
        set of the relators $C_{i, y_j}$  by $C$. We use the term \emph{vertical relators}
        to describe relators in $C$.  A simply connected union of 2-cells corresponding
        to vertical relators is called a \emph{vertical disc}.
        The relators  in $\tilde{S}$ are referred to
        as \emph{horizontal relators}. A disc corresponding with a union of horizontal
        relators is called a \emph{horizontal disc}.
        To summarize, the presentation for $\tilde{G}$ is 
        $$
                \left\langle 
                        \tilde{X}, \tilde{Y}, \tilde{Z} 
                \right\vert  \left.
                        \tilde{R},\tilde{S},\tilde{C} 
                \right\rangle
        $$
        There is a $G$ action on the groupoid $\tilde{G}$. To construct this 
        action, we view $\tilde{G}$ as the homotopy groupoid of the complex
        $K$. Given $g\in G$, there is a path $\alpha$ beginning and ending at $v_0$ 
        representing $g$. If $\gamma$ is a path beginning at $v_0$ 
        ( ie $\gamma$ represents an element of $\tilde{G}$ ), define 
        $g\cdot \gamma = \alpha\gamma$. If 
        $\gamma$ begins at some other vertex $v_k$, then we can conjugate 
        $g$ by the path $\zeta$ that goes from $v_k$ to $v_0$ ( there is
        only one such path, up to homotopy ), then define 
        $g\cdot \gamma = \zeta\alpha\zeta^{-1}\gamma$. It is clear that $\cdot$ defines
        an action on $\tilde{G}$.
        to $\gamma$.

        The advantage of using a groupoid approach is that it is easier to 
        define and compute isoperimetric inequalities, which will serve as
        the main tool used to show that the space is word-hyperbolic.
        We will need to study the geometry of the Cayley graph of this groupoid
        which we will denote by $\tilde{G} = \tilde{G}_H$.
        The Cayley graph of $\tilde{G}$ is just the 1-skeleton of the universal
        cover of the 2-complex $K$ we constructed.

        We now define a metric on $\tilde{G}$.   Let 
        $\psi, \omega$ be positive real numbers. The edges $\tilde{x_i}$
        will be assigned a length of $1$, the edges $\tilde{y}_{i,j}$ have length 
        $\psi^{-j}$. The edges $\tilde{z}_i$ have length $\omega\log(\psi)$.

        The purpose of this metric is to emulate
        the geometry of a Horosphere (or cusp, depending on whether you 
        are looking at $K$ or $\tilde{G}$.) While the metric we use closely
        mirrors $\hyp^n$ horosphere geometry, the area we will define later
        deviates substantially from that of the hyperbolic horosphere.

\subsection{The Cusp Complex}

        Inside $K$, there is a sub-complex $K_H$ consisting of the edges $e_i$ 
        and $e_{i,y}$ 
        and 2-cells $C_{i,y}$ and  $D_{i,s}$.
        Let $\tilde{H}$ be the edge-path homotopy groupoid of $K_H$.
        Denote by $\Gamma_{\tilde{ H }}$ the universal cover of $K_H$.
        We call this space the \emph{cusp complex} of $H$.
        We will use the notation $\Gamma^{(1)}_{\tilde{H}} $ to 
        denote the 1-skeleton of the cusp complex.
        $\Gamma^{(1)}_{\tilde{H}}$ is the Cayley graph of the groupoid generated by 
        $\tilde{Y}$ and $\tilde{Z}$. As such, it can be given a metric by assigning
        lengths to the generators in these sets. We do this in the obvious way -- 
        by assigning the generators their lengths in the $\tilde{G}$-metric.

          The cusp complex depends only on $K_H$, which in turn 
        only depends on the presentation for $H$. It does not depend on $G$.
        Note that $\Gamma_{\tilde{H}}$ is analogous to a horoball.
        The aim of the following discussion will be to prove that $\Gamma^{(1)}_{\tilde{H}}$
        is a $\delta$-hyperbolic metric space.

        Let $\Gamma^{(1)}_H$ be the image of the embedding of the Cayley graph 
        of $H$ in $\Gamma^{(1)}_G$ induced by the inclusion $H\rightarrow G$.

        Given a vertex $v= v_0$ in $\Gamma^{(1)}_H$, there is a sequence
        of vertices $v_0, v_1, v_2 \dots $ in $\Gamma^{(1)}_{\tilde{H}}$, where 
        $v_{i+1} = v_i z_i$. 
        Given a vertex $v_i$, we call   $i$,  the \emph{depth} of 
        the vertex $v_i$.       
        We say that a path is 
        respectively \emph{descending}, \emph{level}, or
        \emph{ascending} if the corresponding sequence of vertices is respectively
        of descending, constant or ascending depth For example, the path
        $(v_0, v_1, v_2 \dots, v_n)$ is descending, and its inverse is ascending. 
        We also define 
        strictly descending and strictly ascending paths, 
        which are (resp) descending
        and ascending paths with no level sub-paths.
        There is a natural map $\pi_H : \Gamma^{(1)}_{\tilde{H}} \rightarrow \Gamma^{(1)}_H$ 
        given by mapping each vertex $v_i$ to $v=v_0$ and mapping each edge
        between vertices $v_i$ and $u_i$ to the edge between $v$ and $u$.
        We call the image under $\pi_H$ of a set the \emph{shadow} of that set.
        Two paths/sets are said to be $H$-parallel if they differ by some map 
                $v_i \mapsto v_{i+k}$ for some $k\in\N$. ( Alternatively, $A$ and $B$
                are $H$-parallel if $A = B z_{i\dots j}$ for some $i,j\in \N$. )
        
\begin{lemma}
                Any geodesic $\alpha$ in $\Gamma^{(1)}_{\tilde{H}}$ can be decomposed
                into $\alpha_1 \alpha_2 \alpha_3$, where
                        \begin{itemize}
                                \item $\alpha_1$ is a strictly descending path, 
                                \item $\alpha_2$ is a level path, 
                                                and the shadow of $\alpha_2$ is a geodesic between 
                                                the shadow of the endpoints of $\alpha_2$.
                                \item $\alpha_3$ is a strictly ascending path.
                        \end{itemize}
\end{lemma}
\begin{proof}
        First, we show that no level edge can precede a descending edge. 
        Assume a path $\alpha$ has a descending edge $( v_i , v_{i+1} )$.
        Let $(v_i , v_{i+1} )$ be the first in a sequence of descending
        edges. Note that if $\alpha$ has an ascending edge that is
        followed immediately by a descending edge, $\alpha$ can easily be
        shortened by canceling those two edges hence $\alpha$ is not geodesic. 
        Let $(u_i , v_i )$ be the last level edge preceding $(v_i , v_{i+1} )$.
        Then we can shorten $\alpha$ by replacing the edges 
        $((u_i, v_i) , (v_i, v_{i+1} ))$ with the shorter 
        path $((u_i, u_{i+1} ) , ( u_{i+1} , v_{i+1} ))$
        so $\alpha$ is not a geodesic.

        We then define $\alpha_1$ to be the maximal initial descending sub-path
        of $\alpha$. $\alpha_1$ is followed by a (possibly empty) level 
        sub-path. Let $\alpha_2$ be a maximal such sub-path. 
        Let $\alpha_3$ be the remaining sub-path of $\alpha$.   
        We can apply the
        argument above to the path $\alpha^{-1}$ to show that it is not possible
        for an ascending edge to be followed by a level edge. From this, we 
        conclude that $\alpha_3$ is strictly ascending.
        The shadow of $\alpha_2$ is a geodesic in $H$, because the restriction 
        of the metric on $\Gamma_{\tilde{H}}$ to a given depth is the same ( up to scaling )
        as the metric on $H$.
\end{proof}     

        We now try to formulate the depth of a path $\alpha$. Let $i$ be the 
        depth of the first vertex of $\alpha$. Let $k$ be the depth of 
        the last vertex. Assume without loss of generality that $i\geq k$.
        Let $L$ be the length of the shadow of $\alpha$.
        Decompose $\alpha$ into $\alpha_1\alpha_2\alpha_3$ where $\alpha_1$
        is strictly descending, $\alpha_2$ is level, and $\alpha_3$ is strictly
        ascending. Let $D$ be the depth of $\alpha_2$. Then the length of
        $\alpha_1$ is $\omega\log(\psi) ( D-i ) $, the length of $\alpha_3$ is 
        $\omega \log(\psi)( D-k )$, 
        and the length
        of $\alpha_2$ is $\psi^{-D} L$. So 
        $$\length ( \alpha  ) = \omega \log(\psi) ( 2D - i - k ) + \psi^{-D} L$$

        One approach to the problem is to extend the domain to $\R$ and
        use differential calculus.
        To minimize the length of $\alpha$, we differentiate the length 
        function with respect to $D$ and set the derivative equal to $0$ :

        $$
                0 = 
                \frac{\partial}{\partial D } ( \length ( \alpha ) ) = 
                2\omega\log(\psi) - \log(\psi) L  \psi^{-D} 
        $$

        Setting the derivative equal to $0$ and solving, we get 

\begin{gather*} 
        \log(\psi) L \psi^{-D} = 2\omega\log(\psi) \\
        \frac{L  }{2\omega } = \psi^{D} \\
                D =  \log_\psi \left( \frac { L } { 2\omega } \right)   
\end{gather*}

        We say that $\log_\psi ( \frac{L}{2\omega} )  $  is the \emph{optimal depth} of the 
        path $\pi_H ( \alpha )$. It is a relative minimum since the second 
        derivative $\frac{\partial^2}{\partial D^2 } ( \length ( \alpha ) ) = 
        L \log(\psi)^2 \psi^{-D} $ is greater than $0$. The shape of $\alpha$ dictates
        that $D \geq i$. So if the optimum depth lies outside this 
        domain, the minimum value of the function corresponds to the endpoint 
        of the domain, ie $D = i$.
        Since we are interested in the minimum value at integer points of $D$, 
        we observe that the value of $D$ that minimizes the length of $\alpha$
        over the integer points of the domain 
        is one of the natural numbers within distance $1$ of $D$. 
        This establishes the following lemma:

\begin{lemma}
        Let $u$ and $v$ be vertices in $\Gamma_{\tilde{H}}$. Let $i$ and $k$ be respectively
        the depths of $u$ and $v$, and suppose $i\geq k$. Let $\alpha$ be a geodesic 
        between $u$ and $v$. Then 
        \begin{itemize}
                \item If  the optimal depth of $\alpha$ is less than $i$,
                                then the maximum depth of any point on $\alpha$ is $i$. 
                \item If the optimal depth of $\alpha$ is $C > i $, then the 
                maximum depth $D$ of any point on $\alpha$ satisfies the inequality
                $$ \vert D - C \vert < 1 $$
        \end{itemize}

\end{lemma}
                
From this, we can also deduce the length of the level sub-path of $\alpha$: 

\begin{lemma}
        Let $\alpha = \alpha_1 \alpha_2 \alpha_3$ be a geodesic in $\Gamma^{(1)}_{\tilde{H}}$ 
        where $\alpha_1$ is strictly descending, $\alpha_2$ is level, and $\alpha_3$
        is strictly ascending. Then the length of $\alpha_2$ is at most $2\omega\psi$.
\end{lemma}
\begin{proof}
        Let $L$ be the length of the shadow $\pi_H( \alpha )$. Then the optimal 
        depth of $\alpha$ is $\log_\psi ( \frac { L } { 2\omega } )  $. 
        So the depth of $\alpha$
        is at least $\log_\psi ( \frac { L } { 2\omega } ) - 1$, 
        and the length of $\alpha_2$ is 
        $\psi^{ - [ \text{Depth} ( \alpha ) ] }L$.
        \begin{gather*}
        \psi^{ - [ \text{ Depth } ( \alpha ) ] } L \leq \\
        \psi^{ - [ \log_\psi( \frac{L}{2\omega} ) - 1  ] } L = \\
        \psi^{  [ - \log_\psi( \frac{L}{2\omega} ) + 1  ] } L = \\
        \frac { \psi L } {  \frac {L}{2\omega} } = \\
        2\psi\omega
        \end{gather*} 
\end{proof}

\begin{theorem}
        The metric space $\Gamma^{(1)}_{\tilde{H}}$ is $\delta$-hyperbolic. 
\end{theorem}
\begin{proof}
        Let $a,b$ and $c$ be vertices in $\Gamma^{(1)}_{\tilde{H}}$. 
        Let $\alpha, \beta,\xi$ be geodesics joining the points
        $(b,c) , (a,c), (a,b)$ respectively.
        Assume without loss of 
        generality that $\length(\alpha) \geq \length(\beta) \geq \length(\xi)$
        By the triangle inequality, the distance 
        $\length(\alpha) \leq \length(\beta) + \length(\xi) \leq 2 \length(\beta)$. 
        Assume that $\alpha = \alpha_1 \alpha_2 \alpha_3$,
        $\beta = \beta_1 \beta_2 \beta_3 $, and $\xi = \xi_1 \xi_2 \xi_3$  are
        decompositions into a descending, level, and ascending segment of
        $\alpha$, $\beta$, and $\xi$.

        So the optimal depth of $\alpha$ is at most $\log_\psi (2)$ more than the optimal
        depth of $\beta$. Taking into account possible rounding errors, 
        the depth of the level segment $\alpha_2$ of $\alpha$ exceeds that of the
        level segment $\beta_2$ by at most $\log_\psi (2) +2$.
        So the paths $\alpha_1$ and $\beta_1$ leave the vertex $c$ and travel 
        along the same line until $\beta_1$ terminates. At this stage,
        $\alpha_1$ descends a distance of no more than $(\log_\psi (2)+2) \omega \log(\psi) $ 
        then $\alpha_2$ travels a distance of at most $2\omega\psi$. $\beta_2$ 
        travels at most $2\omega\psi$.
        So $\alpha_1\alpha_2$ and $\beta_1\beta_2$ stay within a 
        $4\omega\psi + (\log_\psi (2) +2 ) \omega \log(\psi )$ neighborhood of each other.

        Let $\alpha_3'$ and $\beta_3'$ be maximal segments of $\alpha_3$
        and $\beta_3$ that don't intersect $\xi$. Note that any point on
        $\alpha_3'$ can be reached from $\beta_3'$ by traveling along 
        a geodesic parallel to $\xi_2$ ( hence of at length at most 
        $2\omega\psi$ ), and then possibly descending by distance of no more than
        $(\log_\psi(2) + 2)\omega \log(\psi)$.
        So we've established that any point on $\beta$ is either on 
        $\xi$ or within distance $4\omega\psi +(\log_\psi (2) +2 ) \omega \log(\psi ) $ 
        of $\alpha$. The same
        is true reversing the roles of $\alpha$ and $\beta$. 
        It is clear that any point on $\xi$ is either on $\alpha$, $\beta$
        or $\xi_2$ ( since $\xi_1 \subset \alpha, \xi_3 \subset \beta$ ). 
        So $\xi$ stays within distance $2\omega\psi$ of $\alpha \cup \beta$.
        Hence $\Gamma^{(1)}_{\tilde{H}}$ is $\delta$-hyperbolic with 
        $\delta=4\omega\psi + (\log_\psi (2) +2 ) \omega \log(\psi ) $.
\end{proof}

\subsubsection{Boundary Points of Cusps}\label{boundary_points_of_cusps}

Sometimes it is useful to consider boundary points of the cusp complex.
These are analogous to parabolic fixed points in $\hyp^n$.
This is done as follows: let $A$ be the set of points at depth $0$ 
in a cusp complex $\Gamma_{\tilde{H}}$. Then $\Gamma_{\tilde{H}}$  embeds in
the space $A\times [0,\infty )$. 
and that  induces a homeomorphism
into $A\times [0,1) \subset [0,1]/ \sim$ where the relation $\sim$ is
the equivalence generated by the relation $(a,x) \sim (b,y)$  
if $x = 1$ and $y = 1$. By identifying a cusp with its image under the
embedding into $A\times [0,1]/\sim $, and taking the union of that 
with the point $(a,1)$ (note that the choice of $a$ does not matter
because of the equivalence relation), we have obtained a way to add 
a boundary to a cusp. We call the point $(a,1)$ the \emph{boundary}
of the cusp. 

Note that when we take the union of a cusp and its boundary, we no longer have
a metric space. However, it is possible to salvage a cell complex. This is done
as follows: for each 1-cell of depth $0$ in $\Gamma_{\tilde{H}}$, take the union
of all cells beneath it (that is, all one cells parallel to it, all vertical 
1-cells beneath its endpoints, and all 2-cells bounded by the 1-cells
we've just described). Doing this, we obtain an infinite \emph{vertical strip} -- a 2-cell
with infinite sides. The boundary point of the cusp is the essential ingredient
to turn this space into a cell-complex: after we adjoin the boundary point, all
of the infinite strips are triangles with one horizontal edge, and two
``infinite'' vertical edges which meet the boundary point of the cusp.
We call this new complex the  \emph{cusp (cell) complex with boundary} 
 This complex is useful for counting area, because 
each 2-cell has an area of $1$  whereas in the cusp complex, cells can have
arbitrarily small area, and the granularity presents a substantial obstruction
to proving facts about area, because it makes it impossible to determine area
by naively counting 2-cells.
Area here is defined in the 
obvious way: the boundary points added have zero area. Area is computed elsewhere
via the pull-back of the area by the inverse map from the image of
the cusp complex in the \emph{cusp complex with boundary} to the cusp complex.

There are a number of different competing metrics one can use on the cusp complex
or various subsets of the complex. There are problems with the fact that we have
added points at infinity. While we don't get a nice induced path metric, we do
get a distance function that is defined on pairs of points outside the boundary,
and this is quite useful at times. Another metric is the edge-path metric on vertices.
Another metric yet is a ``2-cell metric'' which we will introduce later in the
discussion.

Sometimes, we will need to 
\emph{cut cells in the coset complex with boundary in half.}
This is done as follows: given a triangular 2-cell
in a cusp such that one vertex is a boundary vertex, 
one can divide it into the square and a triangle by subdividing
it as follows: let $T'$ be the triangle, let $T$ be its 
preimage in the coset complex without boundary. Let $e$
be the unique (horizontal) edge contained in $T$ of depth one.
Then $e$ divides $T'$ into two components (a parallelogram and an infinite
triangle).

\subsection{Isoperimetric Bounds For $\tilde{G}$ }

        The purpose of this section will be to establish a linear upper bound on
        the isoperimetric function for $\tilde{G}$.
        Before we proceed, we need to define such a function -- we need to 
        use a variation of the ``classical'' definition, because the edges of the 
        graph are scaled.

                For example, consider the embedding of the cusp complex
                $\tilde{H}$
                of $H = \Z\oplus\Z$ into $\hyp^3$.
                if we choose the embedding carefully, 
                the path metric on $\tilde{H}$ induced by the embedding 
                coincides with the metric on $\Gamma^{(1)}_{\tilde{H}}$.
                However, it is clear that the area of the 2-cells decreases 
                as the length of the edges decrease.

        We define the area function as follows:
        cells $C_{i,y}$ have an area $\psi^{-i} \omega$, cells 
        $D_{i,s}$ have area $\psi^{-i}$.

        To put it more formally, 
        we define a mass function $\mass$ on the defining relators of $\tilde{G}$.
        For each $r\in R$, $\mass ( r ) = 1 $. For each 
        $s\in S_i$, $\mass ( s ) = \psi^{-i }$. 
        For each $c\in c_{i,y}$, $\mass(c) = \psi^{-i} \omega$.

        Note that our choice may seem a little odd. 
        For example, one might expect
        that if the vertical relators were scaled by $\psi^{-i}$, then the horizontal
        relators should be scaled by $\psi^{-2i}$. This would more closely resemble
        the geometry of a hyperbolic cusp, because very small hyperbolic polygons look 
        like Euclidean polygons. Therefore if we treat each horizontal relator as a regular
        polygon and each vertical relator as a trapezium, then this would appear to be the right
        scaling. 
        
        We make our odd-looking choice in the name of pragmatism -- the space
        enjoys the property that we can always reduce the area of a disc in a cusp
        by ``pulling'' the disc towards the boundary point, obtaining something
        analogous to a cone over a boundary point in $\hyp^n$. The benefit of
        this is that it greatly simplifies the process of estimating areas, since
        we are able to eliminate horizontal discs from the computation.

        Given any loop $\alpha$, we define 
        $$
                A(\alpha) = 
                \inf  
                \left\{  A\left\vert 
                          A = \sum_{1,\dots, n  } \right.
                        \mass ( R_i )  , 
                \alpha \sim \Pi_{1, \dots , n} ( p_i R_i p_i^{-1} ) 
                \right\}
        $$
        where $\sim$ is the relation of free equality in $\tilde{G}$ 
         \myfootnote{
                From a geometric standpoint, 
                $\gamma \sim \eta$ if and only if $\eta$ and $\gamma$ are homotopic in the 
                Cayley graph of $\tilde{G}$.
        }

        Given a groupoid $\tilde{G}$ and a chosen set of weighted generators and weighted
        relators, 
        a function $f:\N \mapsto \N$ is said to be an isoperimetric function for 
        $\tilde{G}$ if all loops of length less than $n$ have area no more than $f(n)$.
        We will show that $\tilde{G}$ has a linear isoperimetric function.

        \myfootnote{ Linearity of the isoperimetric function is invariant under a change of
                generators in $G$, which induces an isomorphism of $\tilde{G}$, 
                though the linearity constants themselves may change. Note that one can't say the
        same about arbitrary isometries of the groupoid. This is further complicated by
        the fact that $\tilde{G}$ contains additional structure to the groupoid (namely,
        the area and length weights). 
        }

        \begin{theorem}\label{linear_iso_cuspoff}
                The weighted groupoid $(\tilde{G}, \mass, H )$ has a linear isoperimetric function.     
        \end{theorem}

                The proof will use similar logic to the proof that 
                geometrically finite hyperbolic
                groups have a linear relative isoperimetric function 
                with respect to the cusp groups.
                Given a loop in the Cayley graph of $\tilde{G}$, we will 
                use two classes of reduction moves -- coset reduction
                and length reduction. Using these moves, we can reduce 
                any loop to a relative quasi-geodesic in linear time. Then
                we bound the $\tilde{G}$-length of any loop that is a 
                relative quasi-geodesic.

        \begin{definition}
                A word in $\tilde{G}$ is \emph{coset reduced} if any sub-word 
                in the generators $y_{i,j}$ and $z_i$ is a geodesic with respect
                to the $\tilde{H}$-metric.
        \end{definition}

        \begin{lemma}[Coset Reduction]
                There exists a constant $K\in\R$ such that the following is true:
                any word $\alpha$ of length less than or equal to $n$ 
                in $\tilde{H}$ is freely equal to a word
                $\beta\gamma^{-1}$ where $\beta$ is a loop with area no more than 
                $nK + K$ and $\gamma$ is a $\tilde{H}$-geodesic.
        \end{lemma}
        \begin{proof}
                Let $\gamma$ be the geodesic with the same endpoints as $\alpha$.
                Let  $i$ be the depth of $\alpha(0)$, and let $j$ be the depth
                of the other endpoint of $\alpha$.
                Let $A$ be the area of the shadow of the loop $\alpha \gamma^{-1}$.
                Choose an integer $D$ such that 
        $$ \log_\psi A  \leq D \leq   \log_\psi ( A )  + 1 $$ 
                and $D \geq i, D\geq j$.
                Let $\alpha'$ and $\gamma'$ be paths parallel to the shadows of
                $\alpha$ and $\gamma$
                (respectively) 
                of depth $D$.
                Let $\xi = z_{i\dots D} \alpha' \gamma^{'-1} z_{i\dots D}^{-1}$                 
                The isoperimetric function $f$ and the scaling of the area function
                imply that $\xi$ has an area of no more than $1$.
                Let $\zeta_1$ be the loop 
                $\alpha z_{j \dots D} \alpha^{'-1} z^{-1}_{i \dots D} $.
                Let $\zeta_2$ be the loop
                $z_{i \dots D} \gamma' z^{-1}_{j \dots D}\gamma^{-1} $.
                Then $\alpha \gamma^{-1}$ is homotopic to
                $\zeta_1 \xi \zeta_2$. 
                \begin{gather*}
                \zeta_1 \xi \zeta_2 \\
                =       
                        ( \alpha z_{j \dots D} \alpha^{'-1} z^{-1}_{i \dots D} ) 
                        ( z_{i\dots D} \alpha' \gamma^{'-1} z_{i\dots D}^{-1} )
                        ( z_{i \dots D} \gamma' z^{-1}_{j \dots D}\gamma^{-1} ) 
                        \\
                =
                        ( \alpha z_{j \dots D} )( \alpha^{'-1} z^{-1}_{i \dots D}  
                         z_{i\dots D} \alpha' ) ( \gamma^{'-1} z_{i\dots D}^{-1} 
                         z_{i \dots D} \gamma'  ) ( z^{-1}_{j \dots D}\gamma^{-1} ) 
                        \\
                \sim 
                        ( \alpha z_{j \dots D}  ) 
                        (  z^{-1}_{j \dots D}\gamma^{-1} ) \\
                \sim \alpha \gamma^{-1}
                \end{gather*}   
                The loop $\zeta_1$ has area at most 
                $\length ( \alpha ) \sum_0^{D}  \omega \psi^{-n} $. Similarly,
                the loop $\zeta_2$ has area at most
                $\length ( \gamma ) \sum_0^{D}  \omega \psi^{-n} $.
                The area of $\alpha \gamma^{-1}$ is equal to 
                the area of $\zeta_1 \xi \zeta_2$ which is at most 
                $$ \omega ( \length ( \alpha ) + \length ( \gamma) ) 
                \left[ 
                        \sum_{n=0}^{\infty} (\psi^{-n}) 
                        \right]
                        + 1
                        \leq
                        \omega ( \length ( \alpha ) + \length ( \gamma ) ) 
        \left( \frac{\psi }{\psi - 1 }\right) + 1$$
        \end{proof}

        We can define a pseudo metric on $\Gamma_{\tilde{G}}$ by assigning 
        a length of $0$ to all the edges in $\Gamma_{\tilde{H}}$. This space is the same
        ( modulo changes on sets of diameter $0$ ) as the electric Cayley graph.
        We have notions of relative geodesics, which are geodesics in this 
        pseudo metric. We will call this metric  the \emph{electric metric}. 
        Geodesics in this space are called \emph{electric geodesics}.
        The path length with respect to this metric is called 
        \emph{electric length}.

        If $\alpha$ is a coset reduced word that projects to the 
        identity morphism in $\tilde{G}$, we can employ a \emph{length reduction}
        algorithm similar to that exhibited previously ( in the discussion about
        relatively linear isoperimetric inequalities ).
        Note that in this context, length reduction doesn't reduce the 
        length-proper, it reduces electric length.


\begin{lemma}[Electric Length Reduction]
        There is a constant $E\in \R$ such that for any coset reduced word
        $\alpha$ that is not a $k$-local electric geodesic, and projects to the
        identity morphism in $\tilde{G}$, 
        there is a word $\alpha' \alpha''$ homotopic to $\alpha$ such that 
        $\alpha'$ is a loop of area at most $E$, and $\alpha''$ is a word that
        has electric length strictly less than $\alpha$.
\end{lemma}

Before we begin the proof, we need to make some preliminary observations.
        For any 
        $L\in\R$, there exists a function $\tau: \R\rightarrow \R$ such that
        a non-backtracking loop of $\tilde{G}$-length no more than $L$ 
        has $\tilde{G}$-length no more than $\tau(L)$.
        This claim follows immediately \ref{rel_iso_exists}.
        Note also that there exists a function $h:\R \rightarrow \R$
        such that any loop containing a point of depth $0$ of 
        $\tilde{G}$ length $t$ has area less than $h(t)$.
        So if the electric length of a non backtracking
        loop is less than $t$, then the area is no more than $h(t)$.

\begin{proof}
        First, we address the non-backtracking case.
        Decompose $\alpha$ as $\alpha = \eta_1 \xi \eta_2$ where
        $\eta$ is a minimal sub-word of $\alpha$ that is not 
        an electric geodesic.
        Then $\alpha$ penetrates no more than $k$ cosets.
        Assume the endpoints of $\xi$ are of depth $0$. 
        Let $\xi'$ be a electric geodesic 
        that has the same endpoints of $\xi$ and such that 
        the sub-word between two points in the any $H$ coset
        is a $H$-geodesic.
        By minimality of $\xi$, the cycle $\xi \xi^{ ' -1 } $
        is a non-backtracking cycle of electric length
        no more than $2k$.
        So its $\tilde{G}$-length is no more than $\tau(2k)$ 
        and its area is no more than $h(\tau(2k))$.
        So the result holds with $\alpha' = \eta_1 \xi \xi^{'-1} \eta^{-1}$
        and $\alpha'' = \eta_1 \xi' \eta_2$.

        Now consider the backtracking case.
        In this case, let $\alpha = \eta_1 \eta_2 \eta_3 $ where 
        $\eta_2$ does not backtrack, $\eta_2$ begins and ends in the same 
        $H$-coset, and $\eta_2$ is maximal with this property.
        There are two possibilities to consider -- $\eta_2$ either 
        is or is not a $k$-local geodesic in the electric
        metric. If it is not, choose a minimal sub-word, $\eta$ of $\eta_2$
        that is not an electric geodesic and apply 
        the technique used in the non-backtracking case.
        Otherwise choose a $\tilde{H}$-geodesic $y$ with the same endpoints
        as $\eta_2$, and let $\alpha' = \eta_1 \eta_2 y^{-1} \eta_1^{-1}$,
        and $\alpha'' = \eta_1 y \eta_3$.
        $\eta_2 y^{-1}$ is a non-backtracking cycle of 
        electric length less than $2k$. So the 
        $\tilde{G} $ length is no more than $\tau(2k)$,
        and the area is no more than $h(\tau(2k))$.

\end{proof}

We now proceed with the proof of \ref{linear_iso_cuspoff}:

\begin{proof}[Proof of \ref{linear_iso_cuspoff}]
        Our aim will be to reduce a loop $\alpha$ to a $k$-local electric
        geodesic.
        First, perform a coset reduction on $\alpha$, which contributes
        an area of no more than 
        $2 \omega \length ( \alpha ) \left( \frac{\psi}{\psi-1} \right) + 1$.
        Then we perform a series of length reduction moves, each which 
        only contributes a constant amount of area, and reduces electric 
        length. The important point here is that the length reduction 
        moves can be performed in such a way as to ensure that the resulting
        words are coset reduced.
        The end result is a simple closed non-backtracking loop,
        which has bounded electric length and hence bounded length-proper,
        and bounded area. This shows that we have a linear bound on area --
        the initial coset reduction contributes a linear amount of area,
        and the number of length reductions (contributing constant area)
        is no more than the electric length of $\alpha$, and the final reduction
        of the resulting loop also contributes constant area.
\end{proof}


\subsection{Linear Isoperimetric Inequality Implies $\delta$-Hyperbolicity}

Here, we exploit the linear isoperimetric inequality to show that the cusped off
Cayley graph is $\delta$-hyperbolic. The result itself is hardly surprising --
``hyperbolicity is equivalent to a linear area function'' is a well known 
slogan. However, the proof is nontrivial. The argument here is a generalized 
version of an argument of Gersten and Short in appendix 2 of \cite{gersten}.
In this paper, Gersten and Short prove that a linear isoperimetric inequality in a group 
implies that the group is $\delta$-hyperbolic.

Let $\Gamma = \Gamma^{(2)}$ be the cusp 2-complex. 
When we wish to refer to the cusped off Cayley
graph, we will use the notation $\Gamma^{(1)}$ (since it is indeed the 1-skeleton
of $\Gamma$), and when we wish to refer to
the set of vertices in this complex, we use the notation $\Gamma^{(0)}$.
Let $\rho$ be the maximum word length of any defining relator 
(since this includes the
vertical relators, $\rho \geq 4$).
Let $\bar\Gamma$ be the union of $\Gamma$ and its boundary points. Let $X$ 
be the cusp complex with boundary (which is identical to $\bar\Gamma$ with some 
cells merged). So we view $\Gamma$ as a subset of $\bar\Gamma$, while $X$ is 
the same set of points, but a different CW structure. 
For the purposes of this argument, most of the work will take place in the 
space $X$. However, since the ultimate goal is to prove that $\Gamma$,
 is $\delta$-hyperbolic, we need to understand the $\Gamma$ metric
and area in the context of the space $X$.
Denote by $d_\Gamma: X\times X \rightarrow \R$ the 
restriction of the $\Gamma^{(1)}$ metric to $X^{(1)} - \bdy X$.
More precisely, there is an embedding, 
$f:\Gamma\rightarrow X$. Now use $d_{\Gamma}(a,b) = d(f^{-1}(a), f^{-1}(b))$
where the metric on the right hand side is the usual edge-path metric.
The map $f^{-1}$ is well defined on all but the boundary points.

We will need to introduce a construction called a \emph{disc diagram}.
Since our argument is largely inspired by that in \cite{gersten}, we will 
use a similar definition to theirs.

\begin{definition}
        A \emph{disc diagram}, $(D,h)$, for an edge cycle $\gamma$ 
        in a CW 2-complex $\Gamma$
        is a simply connected CW 2-complex with boundary, $D$ and a continuous map $h : D \rightarrow \Gamma$  with
        the properties that:
   \begin{itemize}
   \item        $h(\bdy D ) = \gamma$
        \item $h$  maps the interiors of $n$-cells in $D$ to interiors of 
                $n$-cells in $h(D)$ by homeomorphisms.
        \end{itemize}
        If the 2-cells in $\Gamma$ are assigned weightings, there is an area
        function on the set of disc diagrams:
        $$ A(D) = \sum_\sigma M(h(\sigma) ) $$
        where $\sigma$ ranges over all 2-cells in $D$, and $M$ is the mass 
        function on $\Gamma$ whose domain is the set of 2-cells in $\Gamma$. 

        Note that nothing in this scheme requires that $f$ be injective, it is
        possible that it is not. This is important, because it is not always
        possible to find an embedded disc whose boundary is $\gamma$ (an obvious 
        example is a power of a relator. Since $\gamma$ is not an embedded circle,
        the restriction of $h$ to $\bdy D$ is not injective either.)

        A \emph{minimal disc diagram} for $\gamma$ is a disc diagram $D$ 
        with the property that  no other disc diagram for $\gamma$ has less area.
        In general, it is not true that minimal disc diagrams always exist,
        though they do exist in the case where disc diagrams exist 
        the weightings have finitely 
        many values. 
   We define an $\epsilon$-minimal disc diagram for a loop as a disc diagram whose 
        area exceeds the area of the loop by at most $\epsilon$. For any $\epsilon >0$,
        if there exists a disc diagram for a loop $\gamma$, then 
        there exists an $\epsilon$-minimal disc diagram for $\gamma$.

\end{definition}

For each $r\in\R$, let $T_r$ be a geodesic triangle in $\Gamma^{(1)}$, 
with vertices $x=x(r),y=y(r) $ and $z=z(r)$ such that there is some point $w$ on 
the side $\overline{xy}$ that is of distance greater $2r$ from any point on 
$\overline{xz}\cup \overline{yz}$.  Let $K$ be the isoperimetric constant.

\begin{figure}
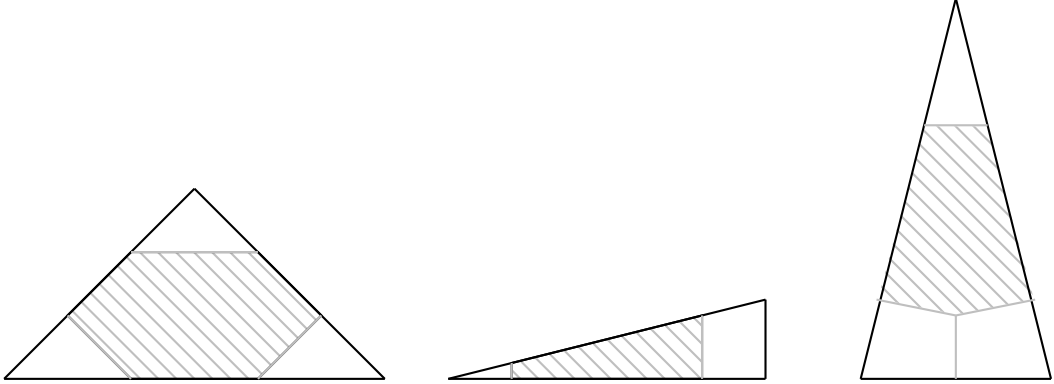
\label{iso_implies_hyp1_fig}
\caption{Truncating a triangle. The three diagrams correspond with case 1 (left), case 2 (center)
and case 3. The sides of length $4\epsilon$ are colored light gray. The black sides are
those of the original triangle.}
 
\begin{center}
 
\psset{xunit=2pc}
\psset{yunit=2pc}

\pspicture(-3,0)(15,6)

\psline[hatchcolor=lightgray,fillstyle=vlines](-1,0)(-2,1)(-1,2)(1,2)(2,1)(1,0)(-1,0)

\psline(-3,0)(3,0)
\psline(-3,0)(0,3)
\psline(3,0)(0,3)

\psline[linecolor=lightgray](-1,0)(-2,1)
\psline[linecolor=lightgray](1,0)(2,1)
\psline[linecolor=lightgray](-1,2)(1,2)

\rput(4,0)
{

\psline[hatchcolor=lightgray,fillstyle=vlines](1,0)(1,0.25)(4,1)(4,0)(1,0)

\psline(0,0)(5,0)
\psline(5,0)(5,1.25)
\psline(0,0)(5,1.25)

\psline[linecolor=lightgray](1,0)(1,0.25)
\psline[linecolor=lightgray](4,0)(4,1)
}

\rput(12,0)
{

\psline[linecolor=white,hatchcolor=lightgray,fillstyle=vlines](-1.25,1.25)(-0.5,4)(0.5,4)(1.25,1.25)(0,1)(-1.25,1.25)

\psline(-1.5,0)(1.5,0)
\psline(-1.5,0)(0,6)
\psline(1.5,0)(0,6)


\psline[linecolor=lightgray](-1.25,1.25)(0,1)
\psline[linecolor=lightgray](1.25,1.25)(0,1)
\psline[linecolor=lightgray](0,0)(0,1)
\psline[linecolor=lightgray](-0.5,4)(0.5,4)

}
\endpspicture
\end{center}
 
\end{figure}

Let $\epsilon > 16$ be a constant such that $\epsilon > 2\rho$.
$r$ can be chosen arbitrarily large, so given $\epsilon$,
we can choose $r$ so that $r>6\epsilon$.
Cut off the corners of $T_r$ in such a way that the remaining segments are maximal with the 
property that each truncated segment is distance at least $4\epsilon$ from the 
other two segments. The end result is as described in fig \ref{iso_implies_hyp1_fig}:
\begin{enumerate}
\item   A non-degenerate hexagon with three non-adjacent sides of length no less than $4\epsilon$.
\item   A non-degenerate quadrilateral with two opposite sides of length no less than $4\epsilon$.
\item   A degenerate hexagon.
\end{enumerate}

We will use $\alpha$, $\beta$, and $\gamma$ denote the interior 
(in $\Gamma^{(1)}$) of the intersection of each side
of the boundary of $T_r$ with the new polygon, and let $\alpha$ be the 
side $\overline{xy}$, so $w\in \alpha$. Note that it's possible 
that $\gamma$ is empty (note that at least 2 sides are non-empty, so we
                can assume that $\alpha$ and $\beta$ are nonempty). For 
convenience, we will define the length of
an empty side to be $0$.

Let $(D,h_D)$ be a minimal disc diagram for 
the path around the boundary
of the truncated triangle in $\Gamma$. Let $\length ( \bdy D )$ denote the 
$\Gamma$-length of this path. 

We need to modify $D$ somewhat. The idea of the argument is that we 
modify $D$ to obtain a space $D'$ where all but a uniformly 
bounded amount of area is contributed by vertical relators and relators
in $R - S$. Then we can simplify the problem by working with a disc diagram
in $X$. 
We need to deal with the fact that $\bdy D$ and $D$ are not unions of cells in $X$.
The disc diagram $D'$ will be a diagram in $\bar{\Gamma}$, in particular
the image of the map $h'_{D'}$ will contain boundary points of $\Gamma$.
We show that 
the modifications to $D$ can be done in such a way as to
only increase its area by at most a bounded 
multiplier. Hence the derived disc $D'$ should have an area bound above by 
a linear multiple of the perimeter of $D$, 
though the actual linearity constant itself will be different.

To begin, we need to modify $\bdy D$. Do this as follows: 
First, we need to deal with corners that are of depth greater than $0$. This is very
simple to do: at each corner of depth greater than $0$, insert a vertical path that
goes from the corner to a vertex of depth 0 and back, namely $z_{0 \dots k} z_{0\dots k}^{-1}$.
This ensures that no corner has nonzero depth.
We now need to deal with the faces of $\bdy D$.
Each sub-path of
$\bdy D$ that is contained inside some cusp, contains a point in that cusp 
of depth greater than $0$, and contains no corners
of $D$ consists of a descending vertical segment,
followed by a horizontal segment, followed by an ascending vertical segment.
Modify this path by deleting the horizontal segment, and extend the vertical 
segments to infinity, and attach the boundary point for that cusp. Note that 
while the length of the new path is undefined, the increase in area is no 
more than $1$, because the two paths are spanned by a $2$-disc 
                whose boundary is an ideal triangle with an
           area no greater than $1$. The spanning disc for $\bdy D$
can be extended to a spanning disc for our modified curve by gluing these ideal 
triangular discs
to $h_D(\bdy D)$, and gluing corresponding triangular discs to the $h_D$ pre-images 
of the attachment points. 
We will denote the new disc by $D_1$.
We then extend $h_D$ in the obvious way to obtain
a map $h_{D_1}$. 
Since the path length of the segments that were ``pushed down'' is at least $1$, and
the area is at most $1$, the effect of this operation is to increase the linearity
constant by no more than $1$. So $h_{D_1}$ has area no more than a factor of $K+1$ 
greater than the length of $\bdy D$.

We then construct a pair $(D',h')$ using several ``pushdown'' operations on horizontal
discs in $D_1$. These operations are performed as follows: first, let $A$ be 
a component of $h_{D_1}(D_1)\cap \Lambda\subset \Gamma$ where $\Lambda$ is a cusp. 
Given a horizontal disc, $R$, it can be pushed down as 
follows (see \ref{fig_pushdown}): 
\begin{figure}
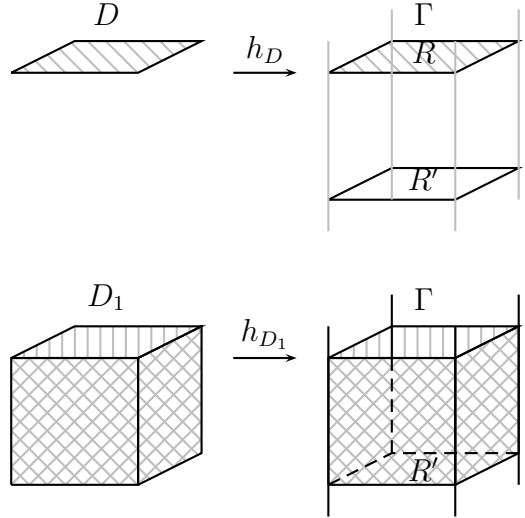
\label{fig_pushdown}
\caption{A pushdown operation on a horizontal 2-cell in $D$}
\begin{center}
\psset{xunit=1pc}
\psset{yunit=1pc}

\pspicture(0,-14)(16,6)


\rput[b](3,1.5){$D$}
\psline[hatchcolor=lightgray,fillstyle=vlines](0,0)(4,0)(6,1)(2,1)(0,0)


\rput(10,0)
{
\rput[b](3,1.5){$\Gamma$}

\psline[hatchcolor=lightgray,fillstyle=vlines](0,0)(4,0)(6,1)(2,1)(0,0)
\rput(0,-4)
{
\psline(0,0)(4,0)(6,1)(2,1)(0,0)
\rput[b](3,0.3){$R'$}
}
\rput[b](3,0.3){$R$}

\psline[linecolor=lightgray](0,-5)(0,1)
\psline[linecolor=lightgray](4,-5)(4,1)
\psline[linecolor=lightgray](2,-4)(2,2)
\psline[linecolor=lightgray](6,-4)(6,2)

}


\rput[b](8,0.3){$h_D$}
\psline{->}(7,0)(9,0)


\rput(0,-13)
{
\rput[b](8,4.3){$h_{D_1}$}
\psline{->}(7,4)(9,4)

\rput[b](3,5.5){$D_1$}

\psline[hatchcolor=lightgray,fillstyle=vlines,hatchangle=0](0,4)(4,4)(6,5)(2,5)(0,4)
\psline[hatchcolor=lightgray,fillstyle=crosshatch](0,0)(4,0)(4,4)(0,4)(0,0)
\psline[hatchcolor=lightgray,fillstyle=crosshatch](4,0)(6,1)(6,5)(4,4)(4,0)
}

\rput(10,-13)
{
\rput[b](3,5.5){$\Gamma$}

\psline[hatchcolor=lightgray,fillstyle=vlines,hatchangle=0](0,4)(4,4)(6,5)(2,5)(0,4)
\psline[hatchcolor=lightgray,fillstyle=crosshatch](0,0)(4,0)(4,4)(0,4)(0,0)
\psline[hatchcolor=lightgray,fillstyle=crosshatch](4,0)(6,1)(6,5)(4,4)(4,0)

\psline[linestyle=dashed](0,0)(2,1)(6,1)
\psline[linestyle=dashed](2,1)(2,4)

\rput(3,0.5){$R'$}

\psline(0,-1)(0,5)
\psline(4,-1)(4,5)
\psline(2,4)(2,6)
\psline(6,0)(6,6)

}

\endpspicture
\end{center}
\end{figure}

Our goal will be to obtain a recursively defined sequence of discs 
$h_i:D_i\rightarrow \bar{\Gamma}$.
In $D_i$, we replace $h^{-1}(R)$ with a prism whose base is $h_{D_i}^{-1}(R)$, then 
remove $h^{-1}(R)$ to obtain a ``prism box'' $B$ whose boundary is $\bdy h_{D_i}^{-1}(R)$.
Let $R'$ be the parallel copy of $R$ whose depth is $\depth ( R ) + 1 $.
There is a unique face opposite to $\bdy  h_{D_i}^{-1} (R)$ in $B$. Map this 
to $R'$. Now extend in the obvious way -- the sides adjacent to $h^{-1}(R)$ 
are mapped to the vertical discs that contain an edge in $\bdy R$ and $\bdy R'$.
If the resulting diagram admits reduction, then we perform it.
First, the net effect of replacing $R$ with $R'$ is that the area is decreased
by $$\area ( R ) - \area ( R' ) = \psi^{- \depth ( R ) } - \psi^{ - \depth ( R' ) } 
= \psi^{ - \depth ( R' ) } ( \psi - 1 )$$
But we also must consider the effect of adding the vertical discs which increases
the area by up to 
$\rho \omega \psi^{- \depth ( R ) } = \rho \omega \psi \psi^{-\depth (R')}$.  
So the net effect is that
area is reduced by at least
$( \psi - 1 - \rho\omega\psi  ) \psi^{ - \depth ( R' ) } $. 
This always corresponds to a net decrease if $\psi > \frac{1}{1-\rho \omega }$ and
$\omega < 1/\rho$.
If we sum 
these increases where $\depth ( R' ) $ ranges between $\depth ( R ) + 1 $ 
and $\infty$, we get that the area is reduced by at least
$$( \psi- 1 - \rho\omega\psi   ) \sum_{n=1}^\infty \psi^{-n} = 
( \psi- 1 - \rho\omega\psi   ) \frac{1}{\psi - 1} $$
So ``pushing down'' decreases area if the above conditions on $\psi$ and $\omega$ 
are met.

This immediately demonstrates an important fact:
\begin{proposition}
        For any loop in $\Gamma$,
        for all $n$, there exists an $\epsilon >0$ such that no $\epsilon$-minimal 
        disc for the loop contains a horizontal disc in a cusp of depth less than $n$.
\end{proposition}
\begin{proof}
        This is immediate from the previous result: if 
   $2\epsilon   = (\psi - 1 - \rho \omega \psi   ) \psi^{-n}$
        then no $\epsilon$-minimal disc can contain a disc in a cusp of depth less than
        $n$, otherwise we could reduce area by at least $2\epsilon$ with a push down 
        operation which implies that the disc is not $\epsilon$-minimal.
\end{proof}

We wish to obtain a disc $D' $ by pushing down ``infinitely many'' times. 
This means doing a fairly extreme form of a pushdown -- instead of 
replacing a horizontal disc with the vertical discs immediately below
its boundary, and adding the horizontal disc immediately below, we
replace the horizontal disc by the union of the infinite strip 
below its boundary, and the ideal point corresponding to that
cusp.
Before we do this, we need the following definition:

\begin{definition}
A 2-cell path is a finite sequence $\sigma_i$ of 2-cells 
such that $\sigma_i \cap \sigma_{i+1}$ includes a 1-cell.
A cellular path from $a$ to $b$  $\in \Gamma$ with $a\neq b$ is a 
2-cell path $\sigma_1, \dots, \sigma_n$ where $a\in\sigma_1$ and 
$b\in\sigma_n$.
If the path is of minimal length, we say that $n$ is the 2-cell distance.
For any $a$,  we consider the empty path to be a 
cellular path from $a$ to itself. (hence the 2-cell distance between a point and
                itself is 0)
2-cell distance defines a pseudo-metric -- it clearly satisfies all the 
metric space axioms with the exception that it may not be defined
on all pairs of points, even in a connected space. ( for example, 
in a complex consisting of two 2-cells that meet at a single vertex,
we will see pairs for which the 2-cell distance is undefined ). If the disc in question has the property
that any two points can be connected by a 2-cell path, then 2-cell
distance is a metric on that disc. 
The relation on pairs of points, defined by the proposition ``$a$ and $b$
are connected by a 2-cell path'' is an equivalence relation, so it makes
sense to speak of ``2-cell path connectedness'' and ``2-cell path components'' 
(which are equivalence classes).
\end{definition}

Suppose we
wish to push down several discs, and moreover, we also wish to 
cancel adjacent pairs of vertical discs. 
We can do this by simultaneously modifying
the domain and range, but we need to make sure that the domain is a disc
after the modification. This works nicely if we wish to push down a region $R$
that is homeomorphic to a disc, because we can remove the interior of
the disc, and cone the boundary off to the ideal point of $\Lambda$, and the
resulting space is still a disc, and we can adjust $h$ so that the interior
of $h^{-1}(R)$ maps to the new pushed down region.
This strategy doesn't work however if the region is
an annulus (because coning over the boundary of an annulus does not 
produce another annulus). We will show that such problematic annular regions 
do not exist, and that we only need to push down regions that are homeomorphic
to a 2-disc.

For a given cusp $\Lambda$, we can 
decompose the 2-cells in $h^{-1}(h(D)\cap \Lambda)$ into 2-cell path components.
We show that any such region that contains a 2-cell is a disc.
Each 2-cell path component is either 
\begin{enumerate}
\item a point on an edge that is not in the boundary of any 2-cell in, $ h^{-1}(h(D)\cap \Lambda)$
\item the closure of a union of 2-cells. 
\end{enumerate}
In case 1, which case we ignore the component.
In case 2, 
the boundaries of the 2-cell components of $h^{-1}(h(D)\cap \Lambda)$ are all 
circles. Note that 2-cell path components have the property that any 1-cell
contained in the component must be in the boundary of some 2-cell in the same 
component, because given an edge $e$ and two distinct interior points
$a,b\subset e$, a 2-cell path from $a$ to $b$ consists of a cell whose boundary
includes $e$.
There is one outermost circle, it makes sense to speak of ``outermost'' since
$h^{-1}(h(D)\cap \Lambda)$ is a subset of a disc.
Take the union of $U$ 2-cells adjacent to one of the inner circles $\alpha$. Note
that each 1-cell in $\alpha$ does indeed lie adjacent to some 2-cell. Let $B$
be the disc in $D$ bounded by $\alpha$. $\alpha$ is of depth $0$, and $U$ contains
no horizontal 2-cells of depth $0$ (from $\epsilon$-minimality), so all the cells 
in $U$ are vertical. If $\psi >2 $, then the area of the vertical cells 
beneath $h(U)$ is less than the area of the cells in $U$, the difference is
$\omega \left( 1 - \frac { 1 } { \psi - 1 } \right)$.
So if $\epsilon < \omega \left( 1 - \frac { 1 } { \psi - 1 } \right)$, $D$ cannot
be $\epsilon$-minimal,
because we can replace $h(U\cup B)$ with the union of the vertical cells directly beneath
$h(U)$ and the boundary point of $\Lambda$, saving at least $\epsilon$ units of area.
So each 2-cell path component of $h^{-1}(h(D)\cap \Lambda)$ is a disc.
We push each such component down to infinity, by performing the same process that we
use to do an infinite pushdown on a 2-cell -- we replace it with a 
disc consisting of the union of the infinite vertical strip below its boundary, and
the boundary point of the cusp $\Lambda$. Now modify $h$ so that the 
2-cell path components of $h^{-1}(h(D)\cap \Lambda)$ are mapped to these new pushed down
discs, to obtain a new diagram $(h',D')$. Observe that $(h',D')$ has the property that
its image does not include any horizontal 2-cells in any cusp.

From $h'(D')$, we obtain a disc diagram $g: E\rightarrow X$ by merging 
vertical discs:
for all pairs $u_0, v_0$ of adjacent vertices of depth $0$ in
a cusp $\Lambda$, perform the following 
procedure: let $C_{i}$ be the unique vertical $2-cell$ whose boundary
contains $u_i, v_i, u_{i+1}$, and $v_{i+1}$.
Then merge the cells to obtain $C = \cup_{i=1}^{\infty} C_i$.
Remove all horizontal discs, and all vertices of depth greater than $1$.
Recall that the new vertical 
disc $C$ has an area of $\omega  \sum_{n=1}^\infty \psi^{-n} = \frac{\omega}{\psi-1}$, 
but are is of infinite $\Gamma$-length. Note also that the finite vertical disc
containing $u_0, v_0, u_1, v_1$ has area $\omega$.

We make some observations:
\begin{itemize}
\item
Observe that there is a uniform bound on edge-path length of the perimeter 
of polygons in
$E$, since the new polygons obtained by merging vertical relators are
infinite triangles.
\item
Also observe that there is \emph{no} uniform bound on the number of
polygons adjacent to each vertex. This is because an ideal vertex could be
adjacent to an arbitrarily large number of triangles.

\item Edge path length in the disc $E$ is not at all related to edge path
length in $\Gamma$ 
because the edges of the ideal triangle correspond to infinite paths in
$\Gamma$. 

\item Least area discs exist in $E$, because the set of weights
        for $2$-cells is the finite set $\left\{ 1, \omega, \frac{\omega}{\psi -1} \right\}$
        
\end{itemize}

For the purpose of the rest of this argument, we will set $\omega = 1 / \psi$.
Note that we require the following conditions on $\psi$:

\begin{itemize}
\item $\psi >2$
\item $D'$ is $\omega ( 1 - \frac{1}{\psi - 1} )$-minimal
\end{itemize}

There  the following proofs will require slightly stronger conditions:
\begin{itemize}
\item $\log ( \psi ) > 1$ (so $\psi > e$)
\end{itemize}

Let $d_E$ be 
the 2-cell metric on the space $E$.
We need to prove the following:

\begin{lemma}
In the space $E$, if $a$ and $b$ 
are non-ideal vertices, then the following holds:

$$  d_{\Gamma} ( f\circ g(a), f\circ g(b) ) \leq \rho d_E ( a,b )$$ 
\end{lemma}
\begin{proof}
We argue this as follows:
Let $\sigma_1, \dots , \sigma_n$ be a 2-cell path in $E$ from $a$ to
$b$. Then connect $f\circ g(a)$ and $f\circ g(b)$ as follows: 

For each $\sigma_i$ with $i>1$, choose a non-ideal vertex $a_i$ with 
        $a_i \in \sigma_i\cap \sigma_{i-1}$. Such a vertex exists because 
        each edge contains two distinct vertices, and only one of those
        vertices may be ideal.
        Join each $a_i$ to $a_{i+1}$
        by a path $\eta_i$ of minimal edge-path length in $\bdy\sigma_i$ 
        that contains no ideal vertices. It is possible to avoid ideal
        vertices because no
        polygon may contain more than one ideal vertex.
        The path $\eta_i$ 
        cannot have length greater than $\rho$. This is because 
        if $\sigma_i$ does not include an ideal vertex, then a path
        $\eta_i$ in $\bdy \sigma_i$ contains no more than $\rho$ one-cells,
        and each of these cells maps via $f\circ g$ to a one-cell of length $1$ in
        $\Gamma$. On the other hand, if $\sigma_i$ includes an ideal vertex,
        then let $\eta_i$ be the path between the other two vertices
        of $\sigma_i$. Then $\eta_i$ maps via $f\circ g$ to a path of length $1$.

        Now let $\eta_0$ be a path from $a$ to $a_1$ and let $\eta_n$
        be a path from $a_n$ to $b$.

Let $\eta$ be the concatenation of the $\eta_i$. Then there's
a path $g(\eta)$ whose length is the same as $\eta$ given
by $g(\eta)(t) = g ( \eta(t) )$. Since $g(\eta)$ contains 
no ideal vertices, we can further compose the path with $f$
to obtain a path whose length is the same.
So $f\circ g(\eta)$ has length of no more than $\rho d_E (a,b)$.
\end{proof}

Now we prove the following lemmas:
\begin{lemma} \label{iso_implies_hyp_l1}
        Let $\theta$ be either one of $\alpha$, $\beta$ or $\gamma$. 
        Then there exists a neighborhood of $\theta$, $N(\theta)$ 
        with the following properties:
        \begin{itemize}
                \item   
                        $$
                                2 (\rho -1 )\text{area} ( N(\theta) ) > 
                                \ell(\theta) \epsilon/{\rho^2 } - C_1
                        $$
                        where $C_1 = C_1(\epsilon )$ is a constant that depends only
                        on $\epsilon$.
                \item
                        All the points in $N(\theta)$ are within 2-cell distance
                        $\lfloor \epsilon / \rho \rfloor +1$ of $\theta$. 
        \end{itemize}

\end{lemma}

\begin{lemma} \label{iso_implies_hyp_2}

        There exists a constant $C_2\in \N$, depending only on $\epsilon$,
        such that 
        $$ \text{area} (D) > (\length ( \alpha ) + \length ( \beta ) + \length ( \gamma ))\epsilon / \rho^2  - C_2(\epsilon) + 2r  \rho $$

\end{lemma}

We need some definitions.

\begin{definition}
        The \emph{star } of a sub-complex $E' \subset E$ is
        the union of  all 2-cells $\sigma$ such that $\sigma\cap E'$
        contains a one-cell or a 2-cell.
        The \emph{corners} of $E$ are the endpoints of the sides
        $\alpha$, $\beta$ and $\gamma$.

\end{definition}

\begin{proof}[ Proof of \ref{iso_implies_hyp_l1} ]

        First, we need to subdivide the disc $E$. 
        In particular, we need to deal with ideal vertices on the boundary of 
        $E$.    
        For each ideal vertex $v_i$ in 
        $\theta$, let $\xi_i$ and $\xi_{i+1}$ be the adjacent one-cells in $\theta$
        (unless $v_i $ is a corner of $E$, in which case, we just have
        one such edge). Let $\eta_i$ be the other edges.
        
        For each $v_i$  we cut each two cell adjacent to $v_i$ in half as described previously(\ref{boundary_points_of_cusps}). The newly added 1-cells are denoted
        $\tau_{i,j}$ and the number of such one cells about each vertex $v_i$
        will be denoted by $k_i$.

        Now we count the number of the 2-cells adjacent to 
        these ideal vertices and compare with the length of the corresponding path
        in $\Gamma$. 
        Then a geodesic between any two points on 
        $f^{-1}(\tau_{i,j})$ and $f^{-1}(\tau_{i,j'})$
        will have a penetration depth of  
        no more than $\log_\psi ( k_i / \omega ) +  1$, and hence length of no more 
        than $2\log_\psi (k_i / \omega )+2+2\psi\omega$. In particular, the ratio between 
        the length of such a geodesic, and the number of 2-cells
        adjacent to the vertex $v_i$ is at least 
        $$
                \frac{k_i}{ 2\log_\psi (k_i / \omega )+2+2\psi\omega } 
                = \frac{k_i}{ 2\log_\psi (k_i \psi )+ 4 } 
                $$
        This function is increasing with respect to $k_i$ if $k_i > 0$, and 
        $\log(\psi)>1$ (one can demonstrate this by differentiating).
        If $k_i = 2$, then $\log_\psi(k_i) < 1$, so 
        $\log_\psi ( k_i \psi ) < 2$, so
        $$\frac {k_i}{ 2\log_\psi (k_i \psi )+ 4 } \geq \frac{2}{8} = 1/4 \geq 1/\rho $$

                Let $N_0$ be the set of 2-cells adjacent to the ideal points in 
                $\theta$.
                Remove $N_0$ from $E_1$ to obtain a disc 
                $E_1 = \overline{E - N_0}\subset E$ such 
                that $$\area ( E_1 - E ) \geq \sum k_i$$

        We will continue with a similar strategy -- remove cells from $E_m$, and demonstrate
   that $E - E_m$ has area that is directly proportional to the $\Gamma$-distance between
        the endpoints of $\theta$.
        Since we've already found a slice of area corresponding with the parabolic sub-words of $\theta$,
        it remains to count 2-cells near $\theta\cap E_1$. So we proceed to do this. 
        Let $\theta_1 = \theta\cap E_1$. Each 1-cell in $\theta_1$ is adjacent to
        at least a single 2-cell. It is possible that there is some double-counting, 
        but no cell can be counted more than $\rho$ times. Let $L$ be the $\Gamma$-length
        of $\theta$. Let $L'$ be the $\Gamma$-length of the portion of $\theta$ that 
        penetrates a cusp corresponding to one of the boundary points we have removed 
        from $E$. 
        Then there are at least $\left\lfloor \frac{(L-L')}{\rho} \right\rfloor$ 
        2-cells adjacent to $\theta \cap E_1$. 
        Denote these cells by $N_1$. Let $E_2 = \overline{E_1 - N_1}$. Let $\theta_1$ be some path in 
        $\overline{(E_2 - E_1)} \cap E_1$.  
        The endpoints of $\theta_1$ are no less than distance $L - 2\rho$ apart (since the endpoints
        are within distance $\rho$ of the endpoints of $\theta$.
        So there are at least $L - 2\rho - L'$ 1-cells in $\theta_1$, excluding the cells $\tau_{i,j}$,
        and they have at least $$ \left\lfloor \frac{ L - 2\rho - L'  }{ \rho } \right\rfloor $$
        adjacent 2-cells. Denote these 2-cells by $N_2$. 

        Iterating in this manner, we attain $E_{i} = E_{i-1} - N_{i-1}$,
        and $N_i$ has at least $$  \frac { L - 2(i-1) \rho - L' } { \rho }  $$
        2-cells. $\theta_i$ contains at least $ L - L' - 2i\rho $
        cells.

        We can iterate the process $\left\lfloor \frac{\epsilon}{\rho}\right\rfloor + 1$
        times and then count the sum of the $N_i$. Note that it's possible
   that some of the $N_i$ are negative, but that's not a problem (we are
   trying to get a lower bound, so the fact that some of the $N_i$ are
        negative merely makes our bound less optimal).   
   Then the sum of the $N_i$ is:
   \begin{gather*}
        \sum_{1, \lfloor \epsilon / \rho \rfloor +1 } 
                 \frac{ L - L' - 2i\rho } { \rho }  \geq \\
         \frac { L - L' } { \rho }  
         \frac{ \epsilon }{\rho } 
        -  \left\lfloor \frac{ \epsilon }{\rho} \right\rfloor
                \left(\left\lfloor \frac{\epsilon}{\rho} \right\rfloor + 1 \right)              
        \end{gather*}

        Let $C_1(\epsilon) =  \left\lfloor \frac{ \epsilon }{\rho} \right\rfloor
                \left(\left\lfloor \frac{\epsilon}{\rho} \right\rfloor + 1 \right)$

        Note that initially, we removed an area of at least
        $L'\epsilon / \rho$, so after counting
        this we have an area of at least 
        $$ \frac{ L \epsilon}{ \rho^2 } - C_1 $$
        The multiplicative constant $\psi$ on the left hand side of the 
        inequality in the statement of \ref{iso_implies_hyp_l1}  is needed because 
        the hypothesis refers area in $D\subset \Gamma^{(2)}$, 
        while we've proved a result
        about area in $E \subset X$. Recall that 
        the area in $X$  is between the area in $\Gamma^{(2)}$  and $\psi$ times that
        area (because of the subdivision of cells near ideal points)

        To prove the other assertion, observe that each 2-cell in $N_i$ is within 2-cell 
        distance 1 of some 2-cell in $N_{i-1}$, since each $N_i$ is in the star of $N_{i-1}$.
        It follows that $N_{ \left\lfloor \epsilon / \rho \right\rfloor + 1 }$ is within 
        2-cell distance $\lfloor \epsilon / \rho \rfloor $ of $N_1$, and 2-cell distance 
        $\lfloor { \epsilon / \rho } \rfloor + 1 $ of $\theta$.

\end{proof}

\begin{proof}[Proof of \ref{iso_implies_hyp_2}]
        First, we note that the neighborhoods of $\alpha, \beta$ and $\gamma$ 
obtained in the above argument are disjoint. They are disjoint because 
they stay within a 2-cell distance of no more than $2\epsilon$ of 
the side that was used in their construction. But the nonempty sides 
of $\alpha, \beta$ and $\gamma$ triangle are sufficiently 
far apart that the 2-cells in each neighborhood do not overlap.

Let $w$ be a point in $\alpha$ that is a distance of at least $2r$
from any other side. Let $\phi'$ be the set of 1-cells in 
$N_\alpha \cap \overline{( \text{int} (D)  - N_\alpha )}$.
Note that all the points in $\phi'$ are 2-cell distance 
of at least $r-2\epsilon$ from $N_\beta \cup N_\gamma$.         
Let $\phi = B \cap \phi'$ where $B$ refers to all points within 
2-cell distance $ r/\rho  $ of $w$ (hence the $\Gamma^{(1)}$ distance 
                from $f^{-1}(w)$ must be no more than $r$ )
        So $\phi$ contains an arc whose endpoints $p_1$ and 
        $p_2$ have the property that 
        $d_{\Gamma^{(1)}} (f^{-1}(p_1), f^{-1}(p_2)) \geq 2r - 4\epsilon$. 
So the 2-cell distance between them is no less than $ \frac{2r-4\epsilon }{\rho}$, and
so the neighborhood of this arc contains at least that many 2-cells.

So the result holds with $C_2(\epsilon) = 4\epsilon / \rho + 3 C_1(\epsilon )$ 
In case 2 and 3,
        $4 /\epsilon / \rho + 2C_1 $ is sufficient, but in case 1, we need $4\epsilon / \rho + 3C_1$.

\end{proof}

We can now use this result to obtain a contradiction to the hypothesis 
of a linear isoperimetric inequality. 

Suppose that $K$ is the isoperimetric constant in the space $\Gamma^{(2)}$. 
Then compensating for the distortions, the isoperimetric constant 
for the space $X$ is no more than $\psi K$.
By hypothesis, the inequality 
$$ A(D) < 12 K\epsilon  + (\length ( \alpha ) + \length ( \beta ) + \length ( \gamma ) )  K$$ 
should hold. 
In terms of area in $X$, this inequality is:
$$ A_X(D) < 12 K'\epsilon  + (\length ( \alpha ) + \length ( \beta ) + \length ( \gamma ) )  K'$$ 
where $K'$ is $\psi K$. Suppose the previous lemma holds.
Let $\epsilon = \max ( K'\rho^2, 16, 2\rho )$. 
Let $r > \max ( 6\epsilon , \rho/2 ( C_2 + 12 K^{'2} \rho^2 ) )$.
The fact that $r \geq \rho/2 ( C_2 + 12K^{'2} \rho^2 ) $ contradicts 
the conclusion of the lemma:

First combine the inequalities to obtain:

$$(\length ( \alpha ) + \length ( \beta ) + \length ( \gamma ) ) \epsilon / \rho^2 - C_2 + 2r/\rho \leq (\length(\alpha )+ \length ( \beta ) + \length ( \gamma ) ) K' + 12 K'\epsilon$$
Since $\epsilon \geq K'\rho^2$, we obtain:
$$ - C_2 + 2r/\rho \leq  12 K^{'2} \rho^2 $$
But substituting $r$  into the expression on the left hand side, we get:
$$ -C_2 + 2r/\rho > -C_2 + 2 ( \rho /2 ( C_2 + 12 K^{'2} \rho^2 ) ) / \rho  
 =   12 K^{'2}/ \rho^2   $$
 But this is contradicts the above inequality.

 \section{An Automaticity Theorem}
\setcounter{footnote}{0}

In the previous section, we introduced the cusped-off Cayley graph.
This construction,  in addition to being a pretty geometric object 
is a useful means of generalizing Epstein's theorem that asserts that
geometrically finite hyperbolic groups are biautomatic. The theorem
depends heavily on the fact that geometrically finite hyperbolic groups
can be shown to act ``nicely'' (properly discontinuously, cocompactly, by isometries)
on a subset of $\hyp^n$ known as the \emph{neutered space}. The obvious 
difficulty with reproducing such an 
argument for a Cayley graph with similar properties is obtaining an ambient
Gromov-hyperbolic space in which to embed the Cayley graph. This is
the problem that the cusped-off Cayley graph addresses.

The aim here is to prove the following result:

\begin{theorem}\label{automatic2}
        Let $G$ be a group that is hyperbolic relative to a finite set of 
        biautomatic subgroups each with a prefix-closed normal
        form $H_1, \dots, H_n$. Then $G$ is biautomatic.
\end{theorem}

A critical subtlety in correctly proving this result is to 
prove that one has obtained a \emph{synchronous}
biautomatic structure. A weaker version of this result where the conclusion
is that $G$ is asynchronously biautomatic is a known 
result \myfootnote{Benson Farb has proven it in unpublished work}

Epstein resolves this problem by showing that one can obtain geodesic normal
forms with the (asynchronous) fellow traveler property, and the fact that
they are geodesic implies that the stronger synchronous fellow traveler 
property follows. The main difficulty in taking this approach to the proof
is in demonstrating that a geodesic $\gamma$ in the neutered space 
stays close to the union of the horospheres on the boundary of the 
neutered space and the $\hyp^n$ geodesic $\gamma'$ between the endpoints 
of $\gamma$. It then follows that the geodesics in the neutered space
enjoy the asynchronous fellow traveler property. The proof is completed 
by mapping the Cayley graph of an appropriate groupoid quasi-isometrically 
into a lattice in the neutered space and observing that statements about
fellow traveling in the neutered space apply to the groupoid.

In Epstein's setting, translating between the discrete and continuous 
settings involves a lot of work. In our setting, we have no such issues --
however, working in a discrete setting comes with problems of it's own,
namely finding an ambient space. On the other hand, having found the cusped off Cayley
graph, a lot of work is done. 
The key result we will need to prove \ref{automatic2} is as follows:

\begin{theorem}\label{theorem2a}
        Let $G$ be a group that is $\delta$-hyperbolic relative to 
        subgroups $H_1, \dots, H_m$, with the $\delta$-bounded coset 
        penetration property.
        There exists a constant $N(\delta)\in \N$ such that for all $n>N$,
        there exists a constant 
        $C(n) \in \N$ such that 
        the following is true:
        Let $x,y\in \Gamma_n$. Let $\alpha$ be a $\Gamma$-geodesic
        between $x$ and $y$. Let $\beta\subset \Gamma_n \Subset \Gamma$ be the 
        $\Gamma_n$-geodesic
        from $x$ to $y$. Let $\gamma \subset \Gamma$ be a path obtained 
        from $\beta$ by replacing all sub-paths of depth $n$ with $\Gamma$-geodesics.
        Then $\gamma$ lies in a $C$-neighborhood of $\alpha$.
\end{theorem}

Before we proceed, we quote the following result:

\begin{lemma}[$k$-local geodesics are quasi-geodesics]\label{k-local}
Let $M$ be a $\delta$-hyperbolic space. Then there exists $N(\delta)\in \N$ 
such that for all $k>N$, there exist $\lambda(k), \epsilon(k)$ such that
any $k$-local geodesic in $M$ is a $(\lambda,\epsilon)$ quasi-geodesic.

\end{lemma}
\begin{proof}

Having noted this result, it suffices to show that 
the path $\gamma$ is in fact a quasi-geodesic in $\Gamma$.

We proceed as follows: let $n = 2k(\delta)$ where $\delta$ is the 
hyperbolicity  constant of the space $\Gamma$, and $k$ is a 
constant obtained from \ref{k-local}. 

Now consider the path $\gamma$. We establish the $k$-local property
by observing sub-paths of $\gamma$. Consider a sub-path, $\eta$. There
are two possibilities, $\eta$ either does or does not contain a point
of depth greater than or equal to $k$. 

We start with the case where $\eta$  does contain 
such a point.  In this case, $\eta$ lies entirely within some cusp. 
First, we consider the sub-case where $\eta$ contains no point of depth 
greater than $2k$. Then $\eta$ must travel vertically, horizontally
(at optimal depth),
then vertically. (either of the two vertical segments may be empty, 
but not both) Any other behavior can be adjusted to yield a shorter
path in $\Gamma_{2k}$, because a segment $uvw$ where 
$v$ is horizontal, $w$ and $u$ are vertical
and $uw$ is a geodesic, then one of $uwv$ or $vuw$ is a shorter 
path than $uvw$.
By a similar argument, if $\eta$ contains a point of depth greater than
$2k$, then it contains a sub-word $\eta'$ that is a $\Gamma$-geodesic, 
and similar logic implies that $\eta'$ must travel vertically, then 
horizontally, then vertically.
Either way, $\eta$ is a geodesic.

We now address the case where $\eta$ does not include any point 
of depth $k$ or more. 
We first demonstrate that in this case, if $\eta$ is a 
$\Gamma_{2k}$-geodesic, then
it is also a $\Gamma$ geodesic.
We do this by demonstrating the contrapositive. 
Suppose that $\eta$ is not a $\Gamma$-geodesic. Then take a geodesic 
$\eta'$ with the same points as $\eta$. For all $x$ in $\eta'$,
$\text{depth}(x) \leq d ( x,\eta(0) ) + \text{depth} (\eta(0)) \leq k + k \leq 2k$.  
So $\eta' \subset \Gamma_{2k}$, hence $\eta$ is not a $\Gamma_{2k}$ geodesic.

We've shown that $\gamma$ is indeed a $k$-local geodesic, so by
\ref{k-local}, the proof is complete.

\end{proof}

To finish the proof of \ref{automatic2}, we need to apply the following 
result of \cite{epstein1}.
First, we provide some context. Let $\tilde G$ be a groupoid. Let $A$ and
$B$ be two different finite ordered weighted sets of generators. Let
$\Gamma ( \tilde G, A)$ and $\Gamma (\tilde G, B)$ respectively be the
Cayley graphs of $\tilde G$ with the generating sets $A$ and $B$. Let 
$\Theta$ be a connected subgraph of $\Gamma(\tilde G,B)$.

\begin{lemma}
        Let $\Theta$ be as above. We suppose that the identity map on vertices
        is a quasi isometry between the metric induced from the path metric
        of $\Theta$ and the path metric of $\Gamma(\tilde G,B)$ (or equivalently, 
                        of $\Gamma(\tilde G,A)$). Let $V$ be a finite state automaton, and
        let $L(V)$ be the language accepted by $V$. Let $L(V)$ be
        prefix closed and consist entirely of certain strings which can be traced
        out entirely within $\Theta$; these strings will not in general be labels 
        on paths starting at the same point. Let $L\subset L(V)$ be the set of
        strings representing paths which are geodesic for the path metric of 
        $\Theta$. We suppose that, for every pair $(v',v'')$ of vertices of 
        $\Gamma(\tilde G,B)$, there is a path from $v$ to $v'$ labelled by an element
        of $L$. We also suppose that there is a number $k$  with the 
        following property:

        \begin{quote}
        Let $w_1$ and $w_2$ be paths in $\Theta$ labelled by elements of $L$. 
        Let their initial points be the vertices $w_1'$ and $w_2'$, 
        and their final points be the vertices $w_1''$ and $w_2''$. 
        Then, in the uniform metric induced by the path metric
        of $\Theta$, their distance apart is $k(d(w_1',w_2') + d(w_1'',w_2'')+1)$
        at the most.
        \end{quote}
        
        Under the above hypothesis, $(L,B)$ is a biautomatic structure for $\tilde G$.

\end{lemma}

\begin{proof}[proof of \ref{automatic2}]
To apply the above lemma, we need to first explain how it applies in our
context. First, $\Theta$ will be our clipped Cayley graph ($\Gamma_{2k}$), 
as will $\Gamma(\tilde G, B)$. We will set $L(V)$ to be the set of 
words that
follow the prefix closed biautomatic normal form on the cusp subgroups.

The property established in \ref{theorem2a} immediately implies the required 
2-sided fellow traveler property for this lemma.
Hence the lemma applies, and we're done.
\end{proof}

\renewcommand\baselinestretch{1}


\begin{thebibliography}{99}
\bibitem[BH1]{bridson1} Martin R.Bridson, Andr\'e Haefliger
\emph{Metric Spaces of Non-Positive Curvature}
Springer-Verlag 1999\\
\bibitem[BH2]{bridson2} Martin R.Bridson, Andr\'e Haefliger
\emph{Metric Spaces of Non-Positive Curvature}
Preprint\\
\bibitem[Bo]{bowditch} Brian H. Bowditch
\emph{Relatively Hyperbolic Groups}
1999 Preprint\\
\bibitem[Brn]{brown1} Kenneth S. Brown
\emph{Cohomology of Groups}, 
1982 Springer-Verlag New York, Inc  \\
\bibitem[F]{farb1} Benson Farb
Relatively Hyperbolic Groups
1998 Preprint\\
\bibitem[E]{epstein1} David B.A. Epstein, J.W. Cannon, D.F. Holt, S.V.F. Levy, M.S. Paterson, W.P. Thurston
\emph{Word Processing in Groups}, 
1992 Jones and Bartlett Publishers, Inc\\
\bibitem[GS]{gersten} Steve Gersten and Hamish Short
\emph{Automatic Structures on Small Cancellation Groups}
1990 Invent. math. 102, 305-334\\
\bibitem[Gr]{greenberg1} Marvin J. Greenberg
\emph{Lectures on Algebraic Topology}
1967 W.A. Benjamin, Inc\\
\bibitem[H]{heintze} E.Heintze, H.-C. Im Hof
\emph{Geometry of Horospheres}
Journal of Differential Geometry
12(1977), 481-491\\
\bibitem[Kl]{klingenberg} Wilhelm Klingenberg
\emph{Riemannian Geometry}
1982 de Gruyter Berlin New York\\
\bibitem[MKS]{mks} Magnus, Karrus and Solitar
\emph{Combinatorial Group theory}\\
\bibitem[Mn]{menasco} Menasco
\emph{Closed Incompressible Surfaces in Alternating Knot and Link Complements}
Topology 1984, 37-44\\
\bibitem[M]{mosher1} Lee Mosher
\emph{Central Quotients of Biautomatic Groups} 
Comm. Math. Helv. 72 no.1 (1997) 16--29\\
\bibitem[NR1]{neumann1} Walter D. Neumann, Lawrence Reeves
\emph{Regular Cocycles and Biautomatic Structures} 
internat. J. Alg. Comp. 6 (1996), 313-324\\
\bibitem[NR2]{neumann3} Walter D. Neumann, Lawrence Reeves
\emph{Central Extensions of Word Hyperbolic Groups} 
Annals of Math. 145 (1997), 183-192\\
\bibitem[NS]{neumann2} Walter D. Neumann, Michael Shapiro
\emph{ Automatic Structures, Rational Growth, and Geometrically Finite Hyperbolic Groups } 
1994 Preprint\\
\bibitem[R]{ratcliffe} John G. Ratcliffe
\emph{Foundations of Hyperbolic Manifolds}
1994 Springer-Verlag New York, Inc\\
\ifthenelse{\boolean{verbose}}
{}
{
\bibitem[Re]{rebbechi} Donovan Rebbechi
\emph{On Boundary-Trivial Central Extensions of Geometrically Finite Hyperbolic Groups }
1999 Preprint ( Rutgers PhD Thesis )
}
\bibitem[Se]{serre} Jean-Pierre Serre
\emph{Trees}
1980 Springer-Verlag New York, Inc\\
\bibitem[Sp]{spanier1} Edwin H.Spanier
\emph{Algebraic Topology} 
1966 Springer-Verlag New York, Inc \\
\bibitem[T1]{thurston1} William P. Thurston
\emph{The Geometry and Topology of Three-Manifolds} 
1997 Preprint
\end{thebibliography}
\end{document}